\documentclass[a4paper]{article}

\usepackage[utf8]{inputenc}  % for older LaTeX distributions
\usepackage{a4wide}
\usepackage{amsthm,amssymb}
\usepackage{mathtools}       % imports and enhances amsmath
\usepackage{graphicx}
    \graphicspath{{./}{figures/}}
\usepackage{bm}              % for boldface math
\usepackage{booktabs} % for tables
\usepackage[UKenglish]{babel}
\usepackage{xcolor} % for nonstandard colors
\usepackage[colorlinks=true,linkcolor=black,citecolor=black]{hyperref}
\usepackage{cite}
\usepackage{url}
\usepackage{float}
\usepackage[shortlabels]{enumitem}

\usepackage{subcaption}
\usepackage{caption}

\newcommand*{\cO}{\mathcal O}
\newcommand*{\C}{\mathcal C}

\newcommand*{\R}{\mathbb R}

\renewcommand{\(}{\begin{equation}}
\renewcommand{\)}{\end{equation}}

\DeclareMathOperator\erfi{erfi}

\newtheorem{theorem}{Theorem}[section]
\newtheorem{proposition}[theorem]{Proposition}
\newtheorem{lemma}[theorem]{Lemma}

\theoremstyle{remark}\newtheorem{remark}[theorem]{Remark}
\usepackage[most]{tcolorbox}
\usepackage{xcolor}

\usepackage[most]{tcolorbox}
\usepackage{xcolor}
\usepackage{xparse}

\numberwithin{equation}{section}

\begin{document}
\title{A particle method for the space homogeneous Lenard--Balescu equation}

\author{
José Antonio Carrillo\thanks{
Mathematical Institute, University of Oxford, Oxford, United Kingdom.
Email: \texttt{carrillo@maths.ox.ac.uk}
}
\and
Timon S. Gutleb\thanks{
School of Computer Science, University of Leeds, Leeds, United Kingdom.
Email: \texttt{T.S.Gutleb@leeds.ac.uk}
}
\and
Andrea Medaglia\thanks{
Department of Mathematics and Computer Science, University of Ferrara, Ferrara, Italy.
Email: \texttt{andrea.medaglia@unife.it}
}
\and
Yijia Tang\thanks{
Scientific Computing Center, Karlsruhe Institute of Technology, Karlsruhe, Germany.
Email: \texttt{yijia.tang@kit.edu}
}
}

\date{}

\maketitle

\begin{abstract}
We propose a structure preserving particle method to numerically approximate the spatially homogeneous Lenard--Balescu equation for collisional plasmas. The Lenard--Balescu operator is an improvement of the Landau equation with a dielectric function taking the dynamic screening effect into account. Our particle method is based on the variational formulation of the Lenard--Balescu equation under suitable regularization analogous to the strategy followed for the Landau collision operator in \cite{CARRILLO2020}. The method guarantees the conservation laws and the entropy dissipation at the semi-discrete level. The novel difficulty of the nonlocal nonlinear dependence of the dielectric function is tackled by rewriting it as integrals in terms of orthogonal and parallel components to the frequency variable. This equivalent form reduces the evaluation of the dielectric function to numerically computing an integral of the directional derivative of the solution on a hyperplane orthogonal to the frequency variable and its one-dimensional Hilbert transform.
An efficient implementation of the numerical approximation for the Lenard--Balescu kernel is derived by exploiting this expression of the dielectric function. We finally validate our method by means of several numerical tests, further allowing us to explore qualitative properties of the Lenard--Balescu equation.

\medskip


\end{abstract}

%\tableofcontents
\section{Introduction}
Collisional relaxation in weakly coupled plasmas \cite{krall1973principles,LifshitzPitaevskii} is usually modelled at the kinetic level by the Landau equation, which describes the cumulative effect of many small-angle binary interactions. In regimes where collective mean-field effects cannot be neglected, the corresponding collisional dynamics is more accurately described by the Lenard--Balescu equation. Its collision operator modifies the Landau operator through a dielectric response, thereby incorporating dynamical screening effects generated by the plasma itself. This additional response is the main feature of the model, but it is also the main source of analytical and numerical difficulty: the collision kernel is nonlinear and nonlocal through its dependence on the solution, it is no longer a convolution kernel in velocity, and its definition involves both a resonant Dirac constraint and a principal-value singularity.

In this paper we consider the spatially homogeneous Lenard--Balescu equation in dimension $d\geq 2$,
\begin{equation}\label{eq:LBeq}
    \partial_t f(v,t)=Q_{LB}(f)(v)
    :=\nabla_v\cdot \int_{\R^d} B(v,v-v_*;\nabla f)
    \bigl(f_*\nabla_v f-f\nabla_{v_*}f_*\bigr)\,\mathrm{d}v_*,
\end{equation}
where $f=f(v,t)$ are $f_*=f(v_*,t)$ are the one-particle distribution functions depending on the velocities $v,v_*\in\R^{d}$ and the time $t>0$. In \eqref{eq:LBeq}, $B=B(v,v-v_*;\nabla f)$ is the Lenard--Balescu kernel, that takes the form
\begin{equation}\label{eq:LBkernel}
    B(v,v-v_*;\nabla f)
    =\int_{\R^d} (k\otimes k)\,\pi \widehat V(k)^2
    \frac{\delta(k\cdot(v-v_*))}{|\varepsilon(k,k\cdot v;\nabla f)|^2}\,\mathrm{d}k,
\end{equation}
where $\widehat V=\widehat V(k)$ is the Fourier transform of the pair-interaction potential, with $k$ wave vector, and $\varepsilon$ is the dielectric function
\begin{equation}\label{eq:dielectric}
    \varepsilon(k,k\cdot v;\nabla f)
    =1+\widehat V(k)\lim_{\eta\to0^+}\int_{\R^d}
    \frac{k\cdot \nabla f(u)}{k\cdot(v-u)-\mathrm{i}\eta}\,\mathrm{d}u.
\end{equation}
By the Sokhotski--Plemelj formula, we obtain
\begin{equation}\label{eq:plem_intro}
    \varepsilon(k,k\cdot v;\nabla f)
    =1+\widehat V(k)\,\mathrm{p.v.}\!\int_{\R^d}
    \frac{k\cdot \nabla f(u)}{k\cdot(v-u)}\,\mathrm{d}u + \mathrm{i}\pi \widehat V(k)\int_{\R^d}
    (k\cdot \nabla f(u))\,\delta(k\cdot(v-u))\,\mathrm{d}u
\end{equation}
where p.v. stands for the principal value integral.
The non-vanishing of $\varepsilon$ is related to Penrose stability \cite{duerinckx2023}, while the factor $|\varepsilon|^{-2}$ in the kernel accounts for collective screening.

The Lenard--Balescu equation was introduced independently by Guernsey \cite{guernsey1960,guernsey1962}, Lenard \cite{lenard1960}, and Balescu \cite{balescu1960,balescu1963}, through derivations based on the BBGKY hierarchy and on Prigogine's theory of irreversibility. Rigorous derivations and validity results remain substantially more difficult than for the Landau equation; recent progress includes the works \cite{duerinckx2021,duerinckxsaintraymond2021}. On the PDE side, the mathematical theory has also developed only recently, including well-posedness and convergence results for smooth interactions \cite{duerinckx2023}, regular Lenard--Balescu equations in spatially inhomogeneous periodic settings \cite{JungNguyen2026}, and thermalization results in simplified hierarchical models \cite{DuerinckxLeBihan2025}. The linearized equation around Maxwellian equilibrium was studied in \cite{merchant1973}, and later coercivity estimates and decay rates were established in \cite{strain2007}.

Numerical work on the Lenard--Balescu equation is much more limited than for the Landau equation, mainly because the collision kernel requires the evaluation of the dielectric response and of the resonant integral in \eqref{eq:LBkernel}. Early contributions include the numerical integration of the Lenard--Balescu equation for isotropic distributions in \cite{Dolinsky1965}, the one-dimensional electrostatic study in \cite{RicciLapenta2002}, and the analysis of the dielectric function in the Lenard--Balescu--Poisson setting in \cite{jasperse1986}. More recently, Laguerre and adaptive spectral methods were developed for the quantum Lenard--Balescu equation in \cite{ScullardEtAl2016,ScullardEtAl2020}. Related generalized Lenard--Balescu computations have also been used for electron--ion equilibration, transport coefficients, and dense-plasma applications, involving quadrature, dynamical screening, local-field corrections, and comparisons with molecular dynamics \cite{VorbergerGericke2009,BenedictEtAl2012,ScullardEtAl2018,FuWangLiKangZhang2015,WhitleyEtAl2015}. Beyond plasma physics, numerical and semi-analytical implementations of Lenard--Balescu-type equations have been developed for long-range interacting systems, especially in stellar dynamics and related models \cite{Heyvaerts2010,Chavanis2012,FouvryEtAl2015A,FouvryEtAl2015B,FouvryEtAl2017,BenettiMarcos2017,HeyvaertsEtAl2017,HamiltonEtAl2018,FouvryBarOrChavanis2019,RouleEtAl2025}. These works show that Lenard--Balescu operators have been evaluated numerically in several long-range settings, but they do not provide a deterministic particle discretization of the nonlinear plasma Lenard--Balescu equation preserving the collision invariants and entropy structure.

At a formal level, the Lenard--Balescu operator has the same fundamental thermodynamic structure as the Landau operator. It conserves mass, momentum, and kinetic energy, dissipates the Boltzmann entropy, and admits Maxwellian distributions as equilibria. In particular, writing the operator in logarithmic form gives
\begin{equation}\label{eq:LB_log}
    Q_{LB}(f)=\nabla_v\cdot\int_{\R^d}B(v,v-v_*;\nabla f)
    ff_*\bigl(\nabla_v\log f-\nabla_{v_*}\log f_*\bigr)\,\mathrm{d}v_*.
\end{equation}
For suitable test functions $\phi$, the corresponding weak formulation is
\begin{equation*}%\label{eq:LB_weak_intro}
    \int_{\R^d}Q_{LB}(f)\phi\,\mathrm{d}v
    =-\frac12\iint_{\R^{2d}}
    (\nabla_v\phi-\nabla_{v_*}\phi_*)\cdot B(v,v-v_*;\nabla f) \, ff_*
    (\nabla_v\log f-\nabla_{v_*}\log f_*)\,\mathrm{d}v\,\mathrm{d}v_*.
\end{equation*}
Taking $\phi=1,v,|v|^2$ gives the conservation laws, while $\phi=\log f$ gives entropy dissipation. Preserving this structure at the discrete level is the central requirement in the construction below.
Indeed, the numerical strategy of the present work is inspired by deterministic particle methods for collisional kinetic equations preserving gradient flow structures \cite{carrillo2019}. A blob method for the homogeneous Landau equation was introduced in  \cite{CARRILLO2020}: the entropy is regularised by convolution, the regularised equation is evaluated on empirical measures, and the resulting particle system preserves mass, momentum, energy, and entropy dissipation at the semi-discrete level. This blob-method strategy has since been further analysed in connection with the regularised Landau equation and its gradient-flow structure \cite{carrillo2023convergence,carrillo2024landau}. It has also been extended in several directions, including random-batch accelerations \cite{CarrilloJinTang}, multispecies Landau systems \cite{CarrilloHuVanFleet2024}, the Vlasov--Maxwell--Landau system \cite{bailo2024collisional}, and Landau equations with uncertainty \cite{Bailo2025UQ}. These developments show that particle discretisations can retain the structural properties of collisional kinetic equations, provided that the regularisation is designed at the level of the entropy and not merely at the level of the flux.

A useful point of comparison is obtained by replacing the dynamical dielectric function with a static response leading to the kernel
\begin{equation}\label{eq:static_kernel}
    B_{\mathrm{static}}(v-v_*)=
    \int_{\R^d}(k\otimes k)\,\pi
    \frac{\widehat V(k)^2}{\varepsilon_{\mathrm{static}}(k)^2}\,
    \delta(k\cdot(v-v_*))\,\mathrm{d}k \quad \mbox{with } \,\varepsilon_{\mathrm{static}}(k)=1+\alpha\widehat V(k), 
\end{equation}
with $\alpha>0$. When $\alpha=0$ and $\widehat V$ is isotropic, the dielectric factor is replaced by $1$ and this kernel reduces to a Landau-type kernel
\begin{equation}\label{eq:LB_to_Landau}
    B_L(v-v_*)=\frac{L}{|v-v_*|}
    \left(I_d-\frac{(v-v_*)\otimes(v-v_*)}{|v-v_*|^2}\right),
\end{equation}
with $I_d$ identity matrix in dimension $d$, and $L>0$ a suitable constant related to the regularity of $V$, see \cite{duerinckx2023} for a detailed discussion.
In this perspective, the Lenard--Balescu equation can be viewed as a dynamically screened version of the Landau equation. 

The purpose of this paper is to construct a deterministic particle method for the spatially homogeneous Lenard--Balescu equation. The starting point is the gradient-flow structure of the collision operator, which allows us to regularise the equation at the level of the Boltzmann entropy, following the strategy developed for the Landau equation. This regularisation makes the equation meaningful on empirical measures and leads to a closed system of ODEs for the particle velocities. The main additional difficulty, compared with the Landau case, is that the Lenard--Balescu kernel is not fixed: it depends on the particle approximation itself through the dielectric function. We show that, despite this extra nonlinear and nonlocal dependence, the resulting semi-discrete method preserves the fundamental structure of the continuous equation, namely conservation of mass, momentum, and energy, together with entropy dissipation. 

We also discuss practical strategies for evaluating the Lenard--Balescu kernel, including static screening, radial reductions, and fully dynamical dielectric functions, as well as reduced-cost approximations based on random batching. The efficient approximation of the nonlocal nonlinear dependence of the dielectric function is tackled by rewriting \eqref{eq:dielectric} as integrals in terms of orthogonal and parallel components to the frequency variable. This change of variables applied to the expression \eqref{eq:plem_intro} involving principal-value integrals allow us to reduce the evaluation of the dielectric function to numerically computing an integral of the directional derivative of the solution on a hyperplane orthogonal to the frequency variable and its one-dimensional Hilbert transform. Therefore, numerical methods for Hilbert transforms and plasma dispersion functions become relevant for its practical evaluation \cite{king_hilbert_2009,olver2011computing,Weideman1994}.

The paper is organised as follows. Section~\ref{sec:2} recalls the gradient-flow formulation of the Lenard--Balescu equation and introduces the entropy regularisation. Section~\ref{sec:3} derives the deterministic particle system and proves the semi-discrete conservation and entropy-dissipation properties. Section~\ref{sec:4} is devoted to the computation of the dielectric function, 
including radial, and fully dynamical cases. Section~\ref{sec:5} discusses implementation aspects and computational cost. Section~\ref{sec:6} presents numerical simulations validating the kernel evaluation, the conservation laws, and the entropy decay of the proposed method.

\section{Gradient flow formulation and regularized Lenard--Balescu equation}
\label{sec:2}
We now recall the variational structure on which the particle method will be based. The idea follows the deterministic particle method introduced for the homogeneous Landau equation in \cite{CARRILLO2020} and further developed in subsequent works on the Landau equation and related gradient-flow formulations \cite{carrillo2022boltzmann,carrillo2023convergence,carrillo2024landau}. The main point is not to regularize the collision flux directly. Instead, we regularize the entropy variable that drives the equation. This is the mechanism that makes empirical measures admissible and, at the same time, preserves the conservative and dissipative structure of the collision operator.

We denote the Boltzmann entropy by
\[
E(f)=\int_{\mathbb{R}^d}f\log f \,\mathrm{d}v,
\]
so that its variation reads
\[
\frac{\delta E}{\delta f}=\log f+1.
\]
Using the logarithmic form \eqref{eq:LB_log}, the Lenard--Balescu equation \eqref{eq:LBeq} can be written as the nonlinear continuity equation
\begin{equation} \label{eq:gradient_flow}
	\partial_t f = -\nabla \cdot (U(f)f),
\end{equation}
where the velocity field is
\begin{equation} \label{eq:velocity_field}
U(f)=-\int_{\mathbb{R}^d} B(v, v-v_*; \nabla f)
\left(\nabla \frac{\delta E}{\delta f}
- \nabla_* \frac{\delta E_*}{\delta f_*}\right)f_*\,\mathrm{d}v_*.
\end{equation}
Here and in the following, we use the short notation
\[
f=f(v), \qquad f_*=f(v_*), \qquad
\nabla=\nabla_v, \qquad \nabla_*=\nabla_{v_*}.
\]
Equation \eqref{eq:gradient_flow} with \eqref{eq:velocity_field} is the Lenard--Balescu analogue of the gradient-flow formulation of the Landau equation. The velocity field is nonlocal and depends on the solution through the kernel $B$, but the dissipated functional is still the Boltzmann entropy.

We next introduce a mollifier $\psi_{\epsilon}$, with $\epsilon>0$, satisfying
\[
\psi_{\epsilon}(v)\geq 0, \qquad
\psi_{\epsilon}(-v)=\psi_{\epsilon}(v), \qquad
\int_{\R^d}\psi_\epsilon(v)\,\mathrm{d}v = 1.
\]
\begin{remark}
Several choices of mollifier are possible. For instance, one may use Gaussian kernels, as in \cite{CARRILLO2020}, or compactly supported kernels such as tensor-product B-splines \cite{bailo2024collisional}. From the analytical point of view, these choices play the same role, provided that the mollifier is nonnegative, symmetric, normalized, and sufficiently regular. From the numerical point of view, however, the choice is more relevant: the size of the support, the regularity order, and the computational cost of evaluating the mollifier affect the implementation. These aspects will be discussed in the subsection \ref{subsec:mollifier}.
\end{remark}
There are two natural ways to regularize the entropy functional. The first one is the symmetric regularization
\begin{equation}
\label{eq:regular_E}
	E_{\epsilon}(f)=\int_{\mathbb{R}^d} (f* \psi_{\epsilon})
        \log (f * \psi_{\epsilon}) \,\mathrm{d}v,
\end{equation}
the second one is the non-symmetric regularization
\begin{equation}
\label{eq:regular_E2}
	E_{\epsilon}(f)=\int_{\mathbb{R}^d} f
        \log (f * \psi_{\epsilon}) \,\mathrm{d}v,
\end{equation}
where only the logarithmic term is regularized. In what follows, $E_\epsilon$ denotes either one of these regularized entropies.

\begin{remark}
From the analytical point of view, the two regularizations play the same role: both replace the singular entropy variable $\log f$ with a smooth quantity that can be evaluated on empirical measures. The difference becomes relevant at the numerical level. When the symmetric regularization \eqref{eq:regular_E} is evaluated on a particle approximation, the entropy variable still involves an integral in velocity space and therefore requires a quadrature. By contrast, the non-symmetric regularization \eqref{eq:regular_E2} leads directly to a fully discrete expression once the empirical density is inserted. These two alternatives will be discussed in more detail in the next sections.
\end{remark}

Replacing $E$ by the regularized entropy $E_\epsilon$ in the velocity field gives the regularized Lenard--Balescu equation
\begin{equation}\label{eq:regular_LB}
	\partial_t f = Q_{LB}^{\epsilon}(f)
        =-\nabla \cdot (U_{\epsilon}(f)f),
\end{equation}
with regularized velocity field
\begin{equation*}
%\label{eq:regular_U}
U_{\epsilon}(f)=-\int_{\mathbb{R}^d} B_\epsilon(v, v-v_*; \nabla f)
\left(\nabla \frac{\delta E_{\epsilon}}{\delta f}
- \nabla_* \frac{\delta E_{\epsilon,*}}{\delta f_*}\right)f_*\,\mathrm{d}v_*.
\end{equation*}
Here $B_\epsilon$ denotes the Lenard--Balescu kernel given by
\begin{equation}
\label{eq:regular_kernel}
B_\epsilon(v,v-v_*;\nabla f)
=
\int_{\mathbb R^d}
(k\otimes k)\,\pi \widehat V(k)^2
\frac{\delta(k\cdot(v-v_*))}
{|\varepsilon_\epsilon(k,k\cdot v;\nabla f)|^2}\,\mathrm{d}k,
\end{equation}
where the regularized dielectric function is
\begin{equation*}
%\label{eq:regular_dielectric}
\varepsilon_\epsilon(k,k\cdot v;\nabla f)
=
1+\widehat V(k)\lim_{\eta\to 0^+}
\int_{\mathbb R^d}
\frac{k\cdot \nabla (\psi_{\epsilon}*f)(u)}
{k\cdot(v-u)-\mathrm{i}\eta}
\,\mathrm{d}u.
\end{equation*}
Thus the regularization enters both in the entropy driving the transport field and in the dielectric response defining the Lenard--Balescu kernel.

We next make explicit the structural properties of the kernel that will be used in the weak formulation.
\begin{lemma}[Structural properties of the regularized kernel]
%\label{lem:regularized_kernel_properties}

Assume that the dielectric function $\varepsilon_\epsilon$ does not vanish and the integral \eqref{eq:regular_kernel} is finite. Then, for all $v
\ne v_*\in\mathbb R^d$, the matrix $B_\epsilon(v,v-v_*;\nabla f)$ is symmetric and positive semidefinite. Moreover, the kernel satisfies that 
\begin{equation}
\label{eq:kernel_exchange_symmetry}
B_\epsilon(v,v-v_*;\nabla f)
=
B_\epsilon(v_*,v_*-v;\nabla f),
\end{equation}
and
\begin{equation}
\label{eq:kernel_null_direction}
B_\epsilon(v,v-v_*;\nabla f)(v-v_*)=0.
\end{equation}
\end{lemma}
\begin{proof}
The symmetry as a matrix follows from the identity
\[
(k\otimes k)^T=k\otimes k.
\]
Moreover, for any $\xi\in\mathbb R^d$, using \eqref{eq:regular_kernel} gives
\[
\begin{aligned}
\xi\cdot B_\epsilon(v,v-v_*;\nabla f)\xi
&=
\int_{\mathbb R^d}
\pi \widehat V(k)^2
\frac{(k\cdot \xi)^2\delta(k\cdot(v-v_*))}
{|\varepsilon_\epsilon(k,k\cdot v;\nabla f)|^2}
\,\mathrm{d}k
\geq 0.
\end{aligned}
\]
Hence $B_\epsilon(v,v-v_*;\nabla f)$ is positive semidefinite.

We now prove the symmetry in the exchange $v\to v_*$. By definition, we have
\[
B_\epsilon(v_*,v_*-v;\nabla f)
=
\int_{\mathbb R^d}
(k\otimes k)\,\pi \widehat V(k)^2
\frac{\delta(k\cdot(v_*-v))}
{|\varepsilon_\epsilon(k,k\cdot v_*;\nabla f)|^2}
\,\mathrm{d}k.
\]
Since $\delta(k\cdot(v_*-v))=\delta(k\cdot(v-v_*))$,
the integral is restricted to the set $\{k\cdot(v-v_*)=0\}$. On this set,
$k\cdot v=k\cdot v_*$ holds, and therefore we infer that
\[
\varepsilon_\epsilon(k,k\cdot v;\nabla f)
=
\varepsilon_\epsilon(k,k\cdot v_*;\nabla f).
\]
Thus the two integrals defining $B_\epsilon(v,v-v_*;\nabla f)$ and $B_\epsilon(v_*,v_*-v;\nabla f)$ coincide, which proves \eqref{eq:kernel_exchange_symmetry}. Finally, \eqref{eq:kernel_null_direction} is obtained by
\[
\begin{aligned}
B_\epsilon(v,v-v_*;\nabla f)(v-v_*)
&=
\int_{\mathbb R^d}
k\,\bigl(k\cdot(v-v_*)\bigr)\,\pi \widehat V(k)^2
\frac{\delta(k\cdot(v-v_*))}
{|\varepsilon_\epsilon(k,k\cdot v;\nabla f)|^2}
\,\mathrm{d}k
=0,
\end{aligned}
\]
because $x\delta(x)=0$ in the sense of distributions.
\end{proof}
We can now write the weak form of the regularized operator. Let $\phi=\phi(v)$ be a smooth test function. From the divergence form of $Q_{LB}^{\epsilon}$, an integration by parts in $v$ gives
\[
\begin{aligned}
\int_{\mathbb R^d} Q_{LB}^{\epsilon}(f)\phi\,\mathrm{d}v
=
-\iint_{\mathbb R^d\times\mathbb R^d}
\nabla\phi\cdot B_\epsilon(v,v-v_*;\nabla f)ff_*
\left(
\nabla \frac{\delta E_{\epsilon}}{\delta f}
-
\nabla_* \frac{\delta E_{\epsilon,*}}{\delta f_*}
\right)
\,\mathrm{d}v\,\mathrm{d}v_* .
\end{aligned}
\]
Using the exchange symmetry \eqref{eq:kernel_exchange_symmetry} and exchanging the variables $v$ and $v_*$, we obtain the symmetrized weak form
\begin{equation}
\label{eq:regular_weak_form}
\begin{split}
\int_{\mathbb R^d} Q_{LB}^{\epsilon}(f)\phi\,\mathrm{d}v
=
-\frac12
\iint_{\mathbb R^d\times \mathbb R^d}
&\bigl(\nabla\phi-\nabla_*\phi_*\bigr)
\cdot B_\epsilon(v,v-v_*;\nabla f) ff_* \\
&\times
\left(
\nabla \frac{\delta E_{\epsilon}}{\delta f}
-
\nabla_* \frac{\delta E_{\epsilon,*}}{\delta f_*}
\right)
\,\mathrm{d}v\,\mathrm{d}v_* .
\end{split}
\end{equation}

We show the resulting structural properties in the following proposition.

\begin{proposition}[Conservation laws and entropy dissipation]
%\label{prop:regularized_structure}
Let $f$ be a smooth, positive, sufficiently decaying solution of the regularized Lenard--Balescu equation \eqref{eq:regular_LB}. Assume that the kernel $B_\epsilon$ is symmetric, positive semidefinite, and satisfies
\[
B_\epsilon(v,v-v_*;\nabla f)(v-v_*)=0.
\]
Then the regularized equation preserves mass, momentum, and kinetic energy:
\[
\frac{\mathrm{d}}{\mathrm{d}t}\int_{\mathbb R^d} f\,\mathrm{d}v=0,
\qquad
\frac{\mathrm{d}}{\mathrm{d}t}\int_{\mathbb R^d} v f\,\mathrm{d}v=\bm{0},
\qquad
\frac{\mathrm{d}}{\mathrm{d}t}\int_{\mathbb R^d}\frac12 |v|^2 f\,\mathrm{d}v=0.
\]
Moreover, the regularized entropy is dissipated along the flow:
\[
\frac{\mathrm{d}}{\mathrm{d}t}E_\epsilon(f)
=
-D_\epsilon(f)\leq 0,
\]
where
\begin{equation*}
% \label{eq:regular_dissipation}
D_\epsilon(f)
=
\frac12
\iint_{\mathbb R^d\times\mathbb R^d}
\left(
\nabla \frac{\delta E_{\epsilon}}{\delta f}
-
\nabla_* \frac{\delta E_{\epsilon,*}}{\delta f_*}
\right)
\cdot B_\epsilon(v,v-v_*;\nabla f)
\left(
\nabla \frac{\delta E_{\epsilon}}{\delta f}
-
\nabla_* \frac{\delta E_{\epsilon,*}}{\delta f_*}
\right)
ff_*\,\mathrm{d}v\,\mathrm{d}v_* .
\end{equation*}
\end{proposition}

\begin{proof}
The conservation laws follow from the weak formulation \eqref{eq:regular_weak_form}. Taking $\phi=1$, we have 
$
\nabla\phi-\nabla_*\phi_*=0, 
$
hence mass is conserved. Taking $\phi=v_i$, $i=1,\dots,d$, gives 
$
\nabla\phi-\nabla_*\phi_*=e_i-e_i=0, 
$
with $e_i$ vector of the direction $i$, hence momentum is conserved. Finally, taking $\phi=\frac12|v|^2$, we obtain
$
\nabla\phi-\nabla_*\phi_*=v-v_*,
$
and the energy contribution vanishes due to \eqref{eq:kernel_null_direction}.
For the entropy identity, we use the chain rule and the equation
\[
\frac{\mathrm{d}}{\mathrm{d}t}E_\epsilon(f)
=
\int_{\mathbb R^d}
\frac{\delta E_\epsilon}{\delta f}
Q_{LB}^{\epsilon}(f)\,\mathrm{d}v.
\]
Applying the weak formulation \eqref{eq:regular_weak_form} with
$
\phi=\frac{\delta E_\epsilon}{\delta f}
$
gives exactly
\[
\frac{\mathrm{d}}{\mathrm{d}t}E_\epsilon(f)
=
-D_\epsilon(f).
\]
Since $B_\epsilon$ is positive semidefinite and $f f_*\geq 0$, we have $D_\epsilon(f)\geq 0$, and therefore
\[
\frac{\mathrm{d}}{\mathrm{d}t}E_\epsilon(f)\leq 0.
\]
\end{proof}
The proposition shows that the entropy regularization preserves the metric structure of the Lenard--Balescu operator at the continuous regularized level. This is the key property that will later pass to the particle discretization: empirical measures are admissible because the entropy variable is regularized, while conservation of the collision invariants and entropy dissipation follow from the same symmetrized weak formulation. 

As a consequence of the entropy decay property, we can also characterize the first variation of the regularized entropy for a stationary state. This proof is similar to that in \cite{GUALDANI2017509,CARRILLO2020}.
\begin{proposition}[Characterization of the equilibrium]
    If $f$ is a equilibrium distribution of the regularized Lenard--Balescu equation \eqref{eq:regular_LB}, then the first variation of $E_{\epsilon}$ is quadratic in $v$. That is, there exist constants $a, c \in \R, b \in \R^d$ determined by the conserved quantities, satisfying 
$$
\frac{\delta E_{\epsilon}}{\delta f} = a + b\cdot v + \frac{c}{2}|v|^2.
$$
\end{proposition}
\begin{proof}
At equilibrium, $Q_{LB}^{\epsilon} (f) = 0$, i.e., $D_{\epsilon}(f) = 0$. Since $B_{\epsilon}$ is symmetric positive semidefinite, $D_{\epsilon} = 0$ holds if and only if
$$
B_\epsilon(v,v-v_*;\nabla f)
\left(
\nabla \frac{\delta E_{\epsilon}}{\delta f}
-
\nabla_* \frac{\delta E_{\epsilon,*}}{\delta f_*}
\right) = 0.
$$
$$
\left(
\nabla \frac{\delta E_{\epsilon}}{\delta f}
-
\nabla_* \frac{\delta E_{\epsilon,*}}{\delta f_*}
\right) \cdot
B_\epsilon(v,v-v_*;\nabla f)
\left(
\nabla \frac{\delta E_{\epsilon}}{\delta f}
-
\nabla_* \frac{\delta E_{\epsilon,*}}{\delta f_*}
\right) = 0.
$$
Refer to Substep 2.2 in the proof of Theorem 1 \cite{duerinckx2023}, this is equivalent to 
$$
\left(
\nabla \frac{\delta E_{\epsilon}}{\delta f}
-
\nabla_* \frac{\delta E_{\epsilon,*}}{\delta f_*}
\right) || (v-v_*)
$$
Thus, there exists $\lambda: \R^d\times \R^d\to \R$ such that
\begin{equation}
\label{eq:linear_dependent}
\left(
\nabla \frac{\delta E_{\epsilon}}{\delta f}
-
\nabla_* \frac{\delta E_{\epsilon,*}}{\delta f_*}
\right) = \lambda (v, v_*) (v-v_*).
\end{equation}
Taking the derivative of the $j$-th component of equation \eqref{eq:linear_dependent} with respect to $v_i$ yields
$$
\partial_{v_i}\partial_{v_j} \frac{\delta E_{\epsilon}}{\delta f} = \partial_{v_i}\lambda (v,v_*)(v_j-v_{*,j})+\lambda(v,v_*)\delta_{ij}.
$$
Let $v=v_*$, we deduce 
$$
\partial_{v_i}\partial_{v_j} \frac{\delta E_{\epsilon}}{\delta f} = \lambda(v,v)\delta_{ij}.
$$
Differentiating again with respect to $v_k$ yields
$$
\partial_{v_k}\partial_{v_i}\partial_{v_j} \frac{\delta E_{\epsilon}}{\delta f} = \partial_{v_k}\lambda(v,v)\delta_{ij}.
$$
Permutating indices $i,j,k \in \{1,\cdots,d\}$ gives the
following identity
$$\partial_{v_k}\lambda(v,v)\delta_{ij}= \partial_{v_i}\lambda(v,v)\delta_{kj}.$$
Take $i=j\neq k$, we obtain $$\partial_{v_k}\lambda(v,v) = 0, \quad  k = 1,\cdots, d.
$$
It's clear from \eqref{eq:linear_dependent} that $\lambda(v,v_*) = \lambda(v_*,v)$, so $\lambda$ is a constant, i.e. $\nabla^2\frac{\delta E_{\epsilon}}{\delta f}  := c I_d $. 
Hence $\frac{\delta E_{\epsilon}}{\delta f}$ is quadratic in $v$ with form
$$
\frac{\delta E_{\epsilon}}{\delta f} = a + b\cdot v + \frac{c}{2}|v|^2.
$$
\end{proof}
\begin{remark}
    Generally speaking, it is not possible to give an explicit characterisation of $f$ at equilibrium.
    One known result is  when we choose the symmetric regularization and the Gaussian mollifier \cite{CARRILLO2020}. In this case, the equilibrium $f$ is a Maxwellian distribution with conserved quantities determined by constants $(a,b,c)$ in an explicit one-to-one mapping. 
\end{remark}

\section{Particle method} 
\label{sec:3}

We now derive the particle approximation associated with the regularized Lenard--Balescu equation introduced in Section~\ref{sec:2}. The construction follows the same principle as the deterministic particle method for the Landau equation: the equation is first regularized at the level of the entropy variable, and only then evaluated on empirical measures.

We approximate $f$ by an empirical measure of the form
\begin{equation*}
 	f^N(v,t)=\sum_{i=1}^N w_i \delta(v-v_i(t)),
 	%\label{eq:empirical}
\end{equation*}
where $v_i(t)\in\mathbb R^d$ is the velocity of particle $i$, $w_i>0$ is its fixed weight, and $N$ is the total number of particles. The associated blob reconstruction is
\begin{equation*}
%\label{eq:blob_reconstruction}
\tilde{f}^N(v,t):=(f^N*\psi_{\epsilon})(v,t)
=\sum_{i=1}^N w_i \psi_{\epsilon}(v-v_i(t)).
\end{equation*}
This smooth reconstruction is used both to evaluate the regularized entropy functional and to visualize the particle approximation.

We now evaluate the two regularized entropies introduced above on the
empirical measure $f^N$. The two choices lead to different expressions
for the discrete term
$\nabla (\delta E_\epsilon^N/\delta f)$. In particular, the first
regularization is only semi-discrete in velocity: although $f^N$ is represented by
particles, the computation of the entropy variable still involves an
integral in velocity space. The second regularization, instead, is fully
discrete in velocity and can be evaluated entirely through finite sums over the
particles.
For the regularized entropy \eqref{eq:regular_E}, we have
\[
E_\epsilon^N:=E_\epsilon(f^N)
=
\int_{\mathbb R^d}\tilde f^N(v)\log \tilde f^N(v)\,\mathrm{d}v,
\]
and the corresponding gradient of the variation of the entropy satisfies
\[
\nabla \frac{\delta E_{\epsilon}^N}{\delta f}(v)
=
\int_{\mathbb R^d}
\nabla \psi_{\epsilon}(v-u)
\log\left(
\sum_{\ell=1}^N w_\ell\psi_{\epsilon}(u-v_\ell)
\right)\,\mathrm{d}u.
\]
Thus, even after replacing $f$ by the empirical measure $f^N$, the
evaluation of the gradient of the variation of the entropy requires the computation of an
additional integral over the velocity variable $u$. For the alternative regularization \eqref{eq:regular_E2}, we obtain
\[
E_\epsilon^N:=E_\epsilon(f^N)
=
\int_{\mathbb R^d} f^N(v)\log \tilde{f}^N(v)\,\mathrm{d}v
=
\sum_{i=1}^N w_i
\log\left(
\sum_{j=1}^N w_j\psi_\epsilon(v_i-v_j)
\right).
\]
In this case, we obtain
\[
\nabla \frac{\delta E_\epsilon^N}{\delta f}(v_i)
=
\sum_{j=1}^N
w_j\nabla\psi_\epsilon(v_i-v_j)
\left(
\frac{1}{\sum_{\ell=1}^N w_\ell\psi_\epsilon(v_i-v_\ell)}
+
\frac{1}{\sum_{\ell=1}^N w_\ell\psi_\epsilon(v_j-v_\ell)}
\right).
\]
Therefore, this second regularization yields a fully discrete formulation:
both the regularized entropy and the corresponding discrete gradient of the variation of the entropy are computed by finite sums over the particles, with no additional quadrature in velocity space.

The particle system is obtained by evaluating the regularized velocity field at the particle locations. Namely, for $i=1,\dots,N$, we set
\begin{equation}
\label{eq:particle_I}
\frac{\mathrm{d}v_i}{\mathrm{d}t}
=
U_{\epsilon}(f^N)(v_i)
=
-\sum_{j=1}^N w_j
B_\epsilon^N(v_i,v_i-v_j)
\left(
\nabla \frac{\delta E_{\epsilon}^N}{\delta f}(v_i)
-
\nabla \frac{\delta E_{\epsilon}^N}{\delta f}(v_j)
\right).
\end{equation}
Here the particle Lenard--Balescu kernel is given by
\begin{equation*}
%\label{eq:particle_kernel}
B_\epsilon^N(v_i,v_i-v_j)
=
\int_{\mathbb R^d}
(k\otimes k)\,\pi \widehat V(k)^2
\frac{\delta(k\cdot(v_i-v_j))}
{|\varepsilon_\epsilon^N(k,k\cdot v_i)|^2}
\,\mathrm{d}k,
\end{equation*}
where
\begin{multline*}
%\label{eq:particle_dielectric}
    \varepsilon_{\epsilon}^N(k,k\cdot v) =\varepsilon(k,k\cdot v;\nabla \tilde{f}^N) = 1+\widehat{V}(k)\lim_{\eta\to 0^+}
\int_{\mathbb R^d}
\frac{k\cdot \nabla \tilde{f}^N (u)}
{k\cdot(v-u)-\mathrm{i}\eta}
\,\mathrm{d}u\\
 =1+\widehat{V}(k)\sum_{\ell=1}^N w_\ell \,\mathrm{p.v.}\!
\int_{\mathbb R^d} 
\frac{k\cdot \nabla\psi_{\epsilon}(u-v_{\ell})}
{k\cdot(v-u)}
\,\mathrm{d}u +\mathrm{i}\pi \widehat{V}(k)\sum_{\ell=1}^N w_\ell \int_{\R^d} k\cdot  \nabla\psi_{\epsilon}(u-v_{\ell})\delta(k\cdot (v-u))\,\mathrm{d}u.
\end{multline*}
The diagonal terms $i=j$ are omitted from the particle sum, or equivalently their contributions are defined to be zero.

\begin{proposition}
Similar to the particle method for the Landau equation, the continuous in time particle method \eqref{eq:particle_I} preserves the structure properties, namely,
\begin{enumerate}
    \item[i)] Conservation of mass, momentum and energy
\begin{equation*}
   \frac{\mathrm{d}}{\mathrm{d}t} \sum_{i=1}^N w_i \phi(v_i) =0, \quad \text{for} \  \phi = 1, v, \frac12|v|^2.
\end{equation*}
\item[ii)] Dissipation of entropy: 
    \begin{equation*}
       \frac{\mathrm{d} E_{\epsilon}^N}{\mathrm{d} t}=-D_{\epsilon}^N \leq 0,
    \end{equation*}
with 
\[
D_\epsilon^N
=
\frac12
\sum_{i,j=1}^N
w_iw_j
\left(
\nabla \frac{\delta E_\epsilon^N}{\delta f}(v_i)
-
\nabla \frac{\delta E_\epsilon^N}{\delta f}(v_j)
\right)
\cdot
B_\epsilon^N(v_i,v_i-v_j)
\left(
\nabla \frac{\delta E_\epsilon^N}{\delta f}(v_i)
-
\nabla \frac{\delta E_\epsilon^N}{\delta f}(v_j)
\right).
\]
\end{enumerate}
\end{proposition}

\begin{proof}
We first consider the conservation laws. Along the particle
	system, for any smooth test function $\phi$, we have
	$$
	\begin{aligned}
		\frac{\mathrm{d}}{\mathrm{d}t}\sum_{i=1}^N w_i\phi(v_i)
		&=~\sum_{i=1}^N w_i \nabla \phi(v_i) \cdot
		U_{\epsilon}(f^N)(v_i)\\
	&=-\sum_{i=1}^N\sum_{j=1}^N w_i w_j \nabla\phi(v_i)\cdot B_\epsilon^N(v_i,v_i-v_j)\left(\nabla \frac{\delta E_{\epsilon}^N}{\delta f}(v_i)- \nabla \frac{\delta E_{\epsilon}^N}{\delta f}(v_j)\right)\\
	&=-\frac12\sum_{i=1}^N\sum_{j=1}^N w_i w_j\left(\nabla\phi(v_i)-\nabla\phi(v_j)\right) \cdot B_\epsilon^N(v_i,v_i-v_j)\left(\nabla \frac{\delta E_{\epsilon}^N}{\delta f}(v_i)- \nabla \frac{\delta E_{\epsilon}^N}{\delta f}(v_j)\right).
	\end{aligned}
	$$
	For $\phi(v)=1$ and $\phi(v)=v$, the difference
	$\nabla\phi(v_i)-\nabla\phi(v_j)$ is zero. For
	 $\phi(v)=\frac12|v|^2$, this difference is proportional
	to $v_i-v_j$. Hence the last expression vanishes, since
	$v_i-v_j \in \ker B_\epsilon^N (v_i,v_i-v_j)$. This proves the conservation of
	mass, momentum, and energy.

We now turn to the entropy dissipation. Take the symmetric regularization for example. Differentiating the
 regularized entropy along the particle trajectories gives
\begin{align*}
	\frac{\mathrm{d}E_{\epsilon}^N}{\mathrm{d}t}=&\int_{\mathbb{R}^d}\sum_{i=1}^N w_i\nabla\psi_{\epsilon}(v-v_i)\cdot\left(-\frac{\mathrm{d} v_i(t)}{\mathrm{d}t}\right)\log \tilde{f}^N(v)\mathrm{d}v\\
& +\int_{\mathbb{R}^d}\tilde{f}^N(v)
\frac{\sum_{\ell=1}^N w_\ell\nabla\psi_{\epsilon}(v-v_\ell) \cdot\left(-\frac{\mathrm{d} v_\ell(t)}{\mathrm{d}t}\right)}
{\tilde{f}^N(v)} \mathrm{d}v
=: I_1+I_2.
\end{align*}
The second term does not contribute. Indeed,
\begin{equation*}
    I_2=-\sum\limits_{\ell=1}^N\int_{\mathbb{R}^d} w_\ell\nabla\psi_{\epsilon}(v-v_\ell) \cdot\frac{\mathrm{d} v_\ell(t)}{\mathrm{d}t}\mathrm{d}v=\frac{\mathrm{d}}{\mathrm{d}t}\sum\limits_{\ell=1}^N w_\ell \int_{\mathbb{R}^d}\psi_{\epsilon}(v-v_\ell) \mathrm{d}v = 0,
%\label{eqn: I_2}
\end{equation*}
since $\int_{\mathbb{R}^d}\psi_{\epsilon}(v-v_\ell) \mathrm{d}v=1$.

It remains to identify the first term. By the definition of the discrete
entropy variable, we obtain
$$
 \begin{aligned}
 	I_1&=\sum_{i=1}^N w_i
 	\left(\int_{\mathbb{R}^d}-\nabla\psi_{\epsilon}(v-v_i)\log \tilde{f}^N(v)
 	\mathrm{d}v\right)\cdot\frac{\mathrm{d} v_i}{\mathrm{d}t}
 	=\sum_{i=1}^N w_i \nabla \frac{\delta E_{\epsilon}^N}{\delta f}(v_i)\cdot\frac{\mathrm{d} v_i}{\mathrm{d}t}\\
 	&=-\sum_{i=1}^N w_i \nabla \frac{\delta E_{\epsilon}^N}{\delta f}(v_i) \cdot \sum_{j=1}^N w_j  B_\epsilon^N(v_i,v_i-v_j)\left(\nabla \frac{\delta E_{\epsilon}^N}{\delta f}(v_i)- \nabla \frac{\delta E_{\epsilon}^N}{\delta f}(v_j)\right)\\
 	&=-\frac{1}{2}\sum\limits_{i=1}^N\sum_{j=1}^N w_i w_j \left(\nabla \frac{\delta E_{\epsilon}^N}{\delta f}(v_i)- \nabla \frac{\delta E_{\epsilon}^N}{\delta f}(v_j)\right)\cdot B_\epsilon^N(v_i,v_i-v_j)\left(\nabla \frac{\delta E_{\epsilon}^N}{\delta f}(v_i)- \nabla \frac{\delta E_{\epsilon}^N}{\delta f}(v_j)\right)\\
 	&=-D_{\epsilon}^N\leq 0.
 \end{aligned}
$$
The last inequality follows from the positive semidefiniteness of
$B_\epsilon^N(v_i,v_i-v_j)$. Hence the discrete entropy is non-increasing in time, and
the proof is complete.
\end{proof}
\begin{remark}
\label{remark:3.2}
The proof above is independent of the particular choice of the entropy
regularization. It only uses the values of the term
\[
\nabla \frac{\delta E_\epsilon^N}{\delta f}(v_i),
\qquad i=1,\dots,N,
\]
the fact that
\[
\frac{\mathrm{d}E_{\epsilon}^N}{\mathrm{d}t} 
=\sum_{i=1}^N w_i \nabla \frac{\delta E_{\epsilon}^N}{\delta f}(v_i)\cdot\frac{\mathrm{d} v_i}{\mathrm{d}t}
\]
and the structural properties of the kernel $B_\epsilon^N(v_i,v_i-v_j)$. Therefore,
the same argument applies both to the semi-discrete regularization, where
the entropy variable still involves an integral in velocity space, and to
the fully discrete regularization, where it is computed by finite sums over
the particles.
\end{remark}

\section{Computation of the dielectric function}
\label{sec:4}
A central additional difficulty in the numerical treatment of the Lenard--Balescu equation, compared with existing structure-preserving methods for the Landau equation, is the computation of the Lenard--Balescu kernel itself. In the Landau case, the interaction kernel is explicit and depends only on the relative velocity. In contrast, the Lenard--Balescu kernel contains the dielectric function $\varepsilon(\cdot)$, which depends on the distribution function through a singular integral. Its evaluation is therefore both analytically and numerically more delicate. For this reason, before discussing the full discretization of the operator, it is necessary to analyze in detail how the dielectric function can be computed. In this section, we first reformulate the dielectric function in terms of a one-dimensional Hilbert transform, which is a form suitable for numerical approximation, and then give one available reduction in radially symmetric case.

In the computation of the dielectric function \eqref{eq:dielectric}, the integral 
\begin{equation*}
%\label{eq:I(k,v)}
I(k,v):=\lim\limits_{\eta\to 0^+}\int_{\R^d} \frac{k \cdot \nabla f(u)}{k\cdot(v-u)-\mathrm{i}\eta} \,\mathrm{d}u,
\end{equation*}
is of crucial importance. 
Throughout this section, we assume that $k\neq 0$ and that $f$ is sufficiently regular and decays sufficiently fast at infinity so that all the integrals below are well defined and the order of integration can be exchanged. Under these assumption, the Sokhotski--Plemelj formula can be applied 
$$
\lim\limits_{\eta\to 0^+}\int_{\R^d} \frac{k \cdot \nabla f(u)}{k\cdot(v-u)-\mathrm{i}\eta} \,\mathrm{d}u =  \textrm{p.v.}\int_{\R^d} \frac{k\cdot\nabla_u f(u)}{k\cdot(v-u)} \,\mathrm{d}u + \mathrm{i} \pi \int_{\R^d} (k\cdot \nabla f(u)) \delta(k\cdot(v-u)) \,\mathrm{d}u
.$$
Hence it is enough to compute the real and imaginary parts separately.

In order to isolate the direction parallel to $k$, we denote $\hat{k}=\frac{k}{|k|}$ and consider the orthogonal decomposition
\begin{equation}
\left\{
\begin{array}{ll}
u = u_\parallel \hat{k} + u_\perp, \quad \textrm{with} \quad u_\parallel = u \cdot  \hat{k},\\
v = v_\parallel \hat{k} + v_\perp, \quad \textrm{with} \quad v_\parallel = v \cdot  \hat{k}.
\end{array}\right.
\label{eq:decomposition}
\end{equation}
With this notation, we can find the following identities
\begin{equation*}
k \cdot \nabla_u f(u) 
%=  |k| \hat{k}\cdot\nabla_u f(u) 
= |k| \hat{k}\cdot\nabla_u f(u_\parallel\hat{k} +u_{\perp} ) 
%= |k| \frac{\partial}{\partial x} f (x \hat{k} +u_{\perp}) 
\end{equation*}
and
\begin{equation*}
k\cdot(v-u) =  |k| (v_\parallel - u_\parallel).   
\end{equation*}
Thus, the real part can be rewritten as 
\begin{equation}
\label{eq:order_integration}
\Re I(k,v)=\textrm{p.v.}\int_{\R^d} \frac{k\cdot\nabla_u f(u)}{k\cdot(v-u)} \,\mathrm{d}u
  = \textrm{p.v.}\int_{\R} \int_{\R^{d-1}}\frac{ \hat{k} \cdot \nabla_u f(u_\parallel\hat{k} +u_{\perp} )}{v_\parallel-u_\parallel} \,\mathrm{d} u_\perp\,\mathrm{d}u_\parallel .
\end{equation} 
Define
\begin{equation}
\label{eq:general_G}
G_{\hat{k}}(x) = \int_{\R^{d-1}}\hat{k}\cdot\nabla_u f(x\hat{k} +u_{\perp} ) \,\mathrm{d}u_{\perp} =  \int_{\R^{d-1}} \frac{\partial}{\partial x} f (x \hat{k} +u_{\perp}) \,\mathrm{d}u_{\perp}.
\end{equation}
With the convention
\[
H[g](x)=\frac{1}{\pi}\,\mathrm{p.v.}\!\int_{\R}\frac{g(y)}{x-y}\,\mathrm{d}y,
\]
we deduce
\begin{equation}\label{eq:hilbertlast}
     \Re I(k,v)  =\pi H[G_{\hat{k}}](v_\parallel).
\end{equation}
This means that the real part can be represented as a one-dimensional Hilbert transform, after integrating the directional derivative of $f$ over the hyperplane orthogonal to $k$.

Similarly, using again the new variables introduced in \eqref{eq:decomposition}, we deduce the identity
$$
\delta (k\cdot(v-u)) = \frac{1}{|k|} \delta (v_\parallel-u_\parallel)
$$
and the imaginary part of $I(k,v)$ is given by
$$
\begin{aligned}
    \Im I(k,v)=\ & \pi \!\int_{\R^d}\! (k\cdot \nabla_u f(u)) \delta(k\cdot(v-u)) \,\mathrm{d}u
 = \pi\!\int_{\R^{d-1}}\!\!\int_{\R} \hat{k}\cdot \nabla_u f(u_{\parallel}\hat{k}+u_{\perp})  \delta(v_{\parallel}-u_{\parallel})\,\mathrm{d}u_{\parallel}\,\mathrm{d}u_\perp\\
 = \ & \pi\!\int_{\R^{d-1}}\!\! \hat{k}\cdot \nabla_u f(v_{\parallel}\hat{k}+u_{\perp}) \,\mathrm{d}u_\perp = \pi G_{\hat{k}}(v_{\parallel}).
\end{aligned}
$$
Combining the two parts, we obtain
\begin{equation}
\label{eq:dielectric_integral}
I(k,v)= \pi H[G_{\hat{k}}](v_\parallel)+ \mathrm{i} \pi G_{\hat{k}}(v_\parallel).
\end{equation}
Therefore, evaluation of the dielectric function $\varepsilon(k,k\cdot v;\nabla f)$ reduces to computing function $G_{\hat{k}}$ and its Hilbert transform numerically. 
Notably, for fixed direction $\hat{k}$, the dependence on $v$ enters only through the scalar
\[ v_\parallel = \frac{k\cdot v}{|k|}.\]

\paragraph{Radially symmetric case}

Now we consider the radially symmetric case. This assumption is relevant as it holds for both the equilibrium Maxwellian and the Bobylev–Krook–Wu (BKW) type distribution. Since the Lenard--Balescu operator is covariant under rotations, the homogeneous Lenard--Balescu equation preserves radial symmetry, provided that the interaction potential $\widehat{V}(k)$ is isotropic. Therefore, this simplification is particular useful when the initial condition is radially symmetric. 
Suppose that $f$ is radially symmetric, i.e., $f(u) = f(|u|)$, we now compute the integrand in \eqref{eq:general_G} as
$$
\frac{\partial}{\partial x} f (x \hat{k} +u_{\perp}) = f'\left(\sqrt{x^2+|u_{\perp}|^2}\right)\frac{ \hat{k}\cdot(x \hat{k} +u_{\perp})}{\sqrt{x^2+|u_{\perp}|^2}}= f'\left(\sqrt{x^2+|u_{\perp}|^2}\right)\frac{x}{\sqrt{x^2+|u_{\perp}|^2}}.
$$
Therefore, plugging this into \eqref{eq:general_G} we obtain 
$$
G_{\hat{k} }(x) = \int_{\R^{d-1}} f'\left(\sqrt{x^2+|u_{\perp}|^2}\right)\frac{x}{\sqrt{x^2+|u_{\perp}|^2}} \,\mathrm{d}u_{\perp}.
$$
Since $G_{\hat{k} }(x)$ does not depend on $\hat{k}$, then we simply denote it by $G(x)$. Using the polar coordinate transformation,
\begin{equation} 
\label{eq:radial_G}
G(x) 
= x S_{d-2}\int_0^{\infty} \frac{f'\left(\sqrt{x^2+y^2}\right)}{\sqrt{x^2+y^2}} y^{d-2}\,\mathrm{d}y
= x S_{d-2} \int_{|x|}^{\infty} f'(r)(r^2-x^2)^{\frac{d-3}{2}}\,\mathrm{d}r.
\end{equation}
Here $S_{d}= 2\pi^{\frac{d+1}{2}} /\Gamma(\frac{d+1}{2})$ denotes the surface area of the $d$-dimensional unit sphere. Particularly, in $3$D case, 
$$
G(x)= xS_1 \int_{|x|}^{\infty}f'(r)dr = -2\pi xf(|x|)
$$
provided that $\lim\limits_{|x|\to \infty} f = 0$.

To conclude, in the radial symmetric case, the function $G_{\hat{k} }$ admits a simplified representation independent of direction $\hat{k}$. Consequently, the result $I(k,v)$ in \eqref{eq:dielectric_integral} amounts to the Hilbert transform of $G$ and relies only on the active point $v_{\parallel}$, though the dielectric function still depends on $k$ due to the presence of $\widehat{V}(k)$.

\paragraph{BKW type distribution}

As a specific example, for the BKW type distribution
\begin{equation}
\label{eq:BKW}
    f (u) = \frac{1}{(2\pi T)^{d/2}}(a+b|u|^2) \exp\left(-\frac{|u|^2}{2T}\right),
\end{equation}
the above quantities can be computed explicitly by 
\[
G(x) = \frac{1}{\sqrt{2\pi T}}\left((3-d)b-\frac{a}{T}-\frac{b}{T} x^2\right)x\exp\left(-\frac{x^2}{2T}\right),
\]
and
\[
H[G](x) = G(x)\erfi\left(\frac{x}{\sqrt{2T}}\right) + \frac{1}{\pi}\left((d-2)b+\frac{a}{T}+\frac{b}{T}x^2 \right).
\]

\section{Implementation of the method}
\label{sec:5}
The particle discretization introduced above leads to a system of interacting particles in which the velocity of each particle depends, in principle, on all the others. Moreover, each pairwise interaction requires the evaluation of the Lenard--Balescu kernel, which itself involves an integral in Fourier space and the computation of the dielectric function. As a consequence, a direct implementation is computationally expensive and quickly becomes prohibitive, especially in high dimensions.

The purpose of this section is to describe the main implementation strategies used to reduce this cost while preserving the structure of the method. We first clarify the initialization of the particles. We then discuss the numerical evaluation of the Lenard--Balescu kernel and the approximation of the dielectric function, emphasizing how precomputation, suitable quadrature rules, and the decomposition of the Fourier integral can be used to accelerate the computation. Next, we describe additional implementation details, including the correction step used to enforce discrete energy conservation after time discretization, the choice of mollifiers used in the regularized entropy and the cell-list data structure used to reduce the cost of interactions. Finally, we introduce a random batch particle method, which replaces the full set of pairwise interactions by interactions within randomly selected smaller batches.

\subsection{Initialization}
\label{subsec:initialization}
We truncate the velocity domain to $\Omega_v=[-L_v,L_v]^d$ and divide it into $n_v$ uniform cells in each coordinate direction. We denote the cells by $C_{\boldsymbol{\alpha}}$ and their centres by $v_{\boldsymbol{\alpha}}$, where
\[
\boldsymbol{\alpha}\in\{1,\ldots,n_v\}^d.
\]
The cell centres form the Cartesian velocity grid, whose mesh size is
\[
h=\frac{2L_v}{n_v}.
\]
We denote by $N_{\rm pc}$ the average number of particles per cell, so that the total number of particles is $N=N_{\rm pc}n_v^d$.

Given an initial distribution $f_0(v)$, we define the mass fraction of each cell as
\[
p_{\boldsymbol{\alpha}}
=
\frac{\displaystyle\int_{C_{\boldsymbol{\alpha}}} f_0(v)\,\mathrm{d}v}
{\displaystyle\int_{\Omega_v} f_0(v)\,\mathrm{d}v}.
\]
We assign to $C_{\boldsymbol{\alpha}}$ an integer number of particles
$N_{\boldsymbol{\alpha}}\simeq Np_{\boldsymbol{\alpha}}$, with the rounding chosen so that
$\sum_{\boldsymbol{\alpha}}N_{\boldsymbol{\alpha}}=N$, and sample their velocities independently and uniformly within the cell. If
$M=\int_{\Omega_v}f_0(v)\mathrm{d}v$ is the total mass, each particle is assigned the weight
$w_i=M/N$. Thus, cells containing a larger fraction of the initial mass receive proportionally more particles.

\subsection{Calculation of Lenard--Balescu kernel}
\label{subsec:5.2}

In this subsection, we discuss the computation of the Lenard--Balescu kernel. For the moment, we assume that the dielectric function can be computed as in Section~\ref{sec:4} and focus only on the integral in Fourier space.

Due to the presence of the Dirac delta function, the integral only takes place on the $(d-1)$-dimensional resonance hyperplane
\[
\{k\in\R^d:\ k\cdot(v-v_*)=0\}.
\]
Motivated by this, we decompose
\[
\R^d
=
\operatorname{span}\{v-v_*\}
\oplus
(v-v_*)^\perp
\]
and split the $k$ variable into components parallel and orthogonal to $v-v_*$. Denote
\[
u=v-v_*,
\qquad
\Phi(k)
=
\frac{(k\otimes k)\,\pi\widehat V(k)^2}
{|\varepsilon_\epsilon(k,k\cdot v;\nabla f)|^2}.
\]
Then
\begin{equation*}
%\label{eq:B_epsilon}
B_\epsilon(u)
=
\int_{\R^d}
\Phi(k)\delta(k\cdot u)\,\mathrm{d}k.
\end{equation*}
For $u\neq 0$, let $R\in SO(d)$ be an orthogonal rotation matrix such that
\[
R\widehat u=e_1,
\qquad
\widehat u=\frac{u}{|u|},
\]
where $e_1$ is the first Cartesian basis vector. The coordinates of $k$ in the rotated reference frame are
\[
p=Rk=
\begin{pmatrix}
p_\parallel\\
p_\perp
\end{pmatrix},
\qquad
p_\parallel\in\R,
\quad
p_\perp\in\R^{d-1}.
\]
Since
\[
k=R^Tp,
\qquad
k\cdot u
=
(R^Tp)\cdot u
=
p\cdot Ru
=
|u|p_\parallel,
\]
and $|\det R^T|=1$, we obtain
\[
\begin{aligned}
B_\epsilon(u)
&=
\int_{\R^d}
\Phi(k)\delta(k\cdot u)\,\mathrm{d}k
\\
&=
\frac{1}{|u|}
\int_{\R^d}
\Phi(R^Tp)\delta(p_\parallel)\,\mathrm{d}p
\\
&=
\frac{1}{|u|}
\int_{\R^{d-1}}
\int_{\R}
\Phi\left(
R^T
\begin{pmatrix}
p_\parallel\\
p_\perp
\end{pmatrix}
\right)
\delta(p_\parallel)
\,\mathrm{d}p_\parallel\,\mathrm{d}p_\perp
\\
&=
\frac{1}{|u|}
\int_{\R^{d-1}}
\Phi\left(
R^T
\begin{pmatrix}
0\\
p_\perp
\end{pmatrix}
\right)
\,\mathrm{d}p_\perp.
\end{aligned}
\]
For the numerical evaluation of the integral, we truncate each transverse wave-number coordinate to $[-k_{\perp,\max},k_{\perp,\max}]$. No explicit cutoff at small $k$ is imposed. Let
$
\left\{
\left(p_\perp^{(m)},\omega_m\right)
\right\}_{m=1}^{N_k}
$
denote the resulting transverse quadrature nodes and weights. For a given relative velocity $u$, the corresponding nodes in the original Fourier coordinates are
\[
k^{(m)}(u)
=
R^T
\begin{pmatrix}
0\\
p_\perp^{(m)}
\end{pmatrix},
\qquad
m=1,\ldots,N_k.
\]
The Lenard--Balescu kernel is then approximated by
\[
B_\epsilon(u) \approx \frac{1}{|u|} \sum_{m=1}^{N_k} \omega_m \Phi\left(k^{(m)}(u)\right).
\]

In the numerical experiments, we use $n_k$ nodes in each of the $d-1$ transverse directions, so that
$
N_k=n_k^{d-1}
$. The nodes at each transverse direction are truncated to $k_{\perp,\max}$ and obtained using a sinh-mapped Gauss--Legendre quadrature, where $k_{\rm scale}$ controls their concentration near the origin. 
Thus, the cost of evaluating the kernel for a given velocity pair $(v,v_*)$ is
$
\mathcal O(N_k) = \mathcal O(n_k^{d-1})$.
The values of $n_k$, $k_{\perp,\max}$, and $k_{\rm scale}$ used in each experiment are reported in the corresponding numerical subsection.

For notational simplicity, we omit the dependence of $R$ on $u$. In two dimensions, for
\[
\widehat u=(\cos\theta,\sin\theta)^T,
\]
the rotation matrix can be chosen as
\[
R=
\begin{pmatrix}
\cos\theta & \sin\theta\\
-\sin\theta & \cos\theta
\end{pmatrix}.
\]
In three dimensions, for
\[
\widehat u
=
(\sin\theta\cos\varphi,
 \sin\theta\sin\varphi,
 \cos\theta)^T,
\]
it can be chosen as
\[
R=
\begin{pmatrix}
\sin\theta\cos\varphi
&
\sin\theta\sin\varphi
&
\cos\theta
\\
-\sin\varphi
&
\cos\varphi
&
0
\\
-\cos\theta\cos\varphi
&
-\cos\theta\sin\varphi
&
\sin\theta
\end{pmatrix}.
\]

In the general case, since $\varepsilon_\epsilon$ depends on the dynamic solution,
\[
\Phi\left(
R^T
\begin{pmatrix}
0\\
p_\perp^{(m)}
\end{pmatrix}
\right)
\]
must be computed online. However, in the static case, due to the rotational invariance of $|k|$,
\[
\Phi(k)
=
\frac{(k\otimes k)\,\pi\widehat V(k)^2}
{(1+\alpha\widehat V(k))^2}
=
R^T
\frac{(p\otimes p)\,\pi\widehat V(p)^2}
{(1+\alpha\widehat V(p))^2}
R.
\]
Therefore,
\[
\frac{(p\otimes p)\,\pi\widehat V(p)^2}
{(1+\alpha\widehat V(p))^2}
\]
can be precomputed at the Fourier quadrature nodes.

\subsection{Approximation of the dielectric function}

For the numerical evaluation of the dielectric function, we reconstruct a smooth distribution from the particles. Omitting the time dependence, we define
\[
\widetilde f^N(v)
=
(f^N*\psi_\epsilon)(v)
=
\sum_{i=1}^N w_i\psi_\epsilon(v-v_i).
\]
On the Cartesian velocity grid introduced in Subsection~\ref{subsec:initialization}, we set
\[
\widetilde f_{\boldsymbol{\alpha}}^N
=
\widetilde f^N(v_{\boldsymbol{\alpha}})
=
\sum_{i=1}^N
w_i\psi_\epsilon(v_{\boldsymbol{\alpha}}-v_i).
\]
We denote by $\nabla_h\widetilde f^N$ the numerical gradient of the reconstructed distribution on this grid and use interpolation whenever it must be evaluated away from the grid points.

\paragraph{General case}

For a fixed Fourier direction
\[
\widehat k=\frac{k}{|k|},
\]
let $\{x_\ell\}$ be the uniform grid, with spacing $h$, in the direction parallel to $\widehat k$. For each $x_\ell$, let
\[
\left\{
\left(u_{\perp,q}^{(\ell)},\omega_{\perp,q}^{(\ell)}\right)
\right\}_{q\in\mathcal J_\ell}
\]
denote the quadrature nodes and weights on the section of the velocity domain orthogonal to $\widehat k$. The discrete counterpart of \eqref{eq:general_G} is
\begin{equation}
\label{eq:discrete_general_G}
G_{\widehat k,h}^N(x_\ell)
=
\sum_{q\in\mathcal J_\ell}
\omega_{\perp,q}^{(\ell)}
\widehat k\cdot
\nabla_h\widetilde f^N
\left(
x_\ell\widehat k+u_{\perp,q}^{(\ell)}
\right),
\end{equation}
where the values of $\nabla_h\widetilde f^N$ at the rotated points are obtained by interpolation from the Cartesian velocity grid.

Once the values of $G_{\widehat k,h}^N$ have been computed, we evaluate their discrete Hilbert transform, denoted by
\[
H_h[G_{\widehat k,h}^N],
\]
using an FFT-based method. Since
\[
v_\parallel=\frac{k\cdot v}{|k|}
\]
does not generally coincide with one of the points $x_\ell$, both
$G_{\widehat k,h}^N$ and its Hilbert transform are interpolated at
$v_\parallel$. The discrete counterpart of \eqref{eq:dielectric_integral} is therefore
\begin{equation*}
%\label{eq:discrete_dielectric_integral}
I_h^N(k,v)
=
\pi H_h[G_{\widehat k,h}^N](v_\parallel)
+
\mathrm{i}\pi
G_{\widehat k,h}^N(v_\parallel).
\end{equation*}
Accordingly, the dielectric function is approximated by
\begin{equation*}
%\label{eq:discrete_dielectric}
\varepsilon_{\epsilon,h}^N(k,k\cdot v)
=
1+\widehat V(k)I_h^N(k,v).
\end{equation*}

Alternatively, one may exchange the order of the parallel and transverse integrations in \eqref{eq:order_integration}. For each fixed $u_\perp\in\widehat k^\perp$, define
\[
F_{u_\perp}(x) = \widehat k\cdot\nabla_u f(x\widehat k+u_\perp).
\]
Then
\[
\begin{aligned}
\Re I(k,v) &= \int_{\widehat k^\perp} \left( \operatorname{p.v.}\int_{\mathbb R} \frac{F_{u_\perp}(x)} {v_\parallel-x} \,\mathrm{d}x \right) \mathrm{d}u_\perp \\
&= \pi \int_{\widehat k^\perp} H[F_{u_\perp}](v_\parallel) \,\mathrm{d}u_\perp,
\end{aligned}
\]
while
\[
\Im I(k,v) = \pi \int_{\widehat k^\perp} F_{u_\perp}(v_\parallel) \,\mathrm{d}u_\perp.
\]
Equivalently,
\[
I(k,v) = \pi \int_{\widehat k^\perp} \left( H[F_{u_\perp}](v_\parallel) + \mathrm{i}F_{u_\perp}(v_\parallel) \right) \,\mathrm{d}u_\perp.
\]
Numerically, this alternative would require applying a separate discrete Hilbert transform at every transverse quadrature node before performing the transverse integration. This is less efficient in a parallel implementation, since it introduces several FFT computations inside an already parallel quadrature loop. We therefore follow \eqref{eq:hilbertlast}: for each direction $\widehat k$, we first complete the transverse quadrature in \eqref{eq:discrete_general_G} and then apply a single Hilbert transform to the resulting profile.

\paragraph{Radially symmetric case}

Under the radial-symmetry assumption, we first average the reconstructed values over spherical shells. Let
\[
r_q=\left(q+\frac12\right)h,
\qquad
q=0,\ldots,n_r-1,
\]
and define
\[
\mathcal S_q
=
\left\{
\boldsymbol{\alpha}:
qh\leq |v_{\boldsymbol{\alpha}}|<(q+1)h
\right\}.
\]
The discrete radial profile is
\[
f_{h,q}^N
=
\frac{1}{\#\mathcal S_q}
\sum_{\boldsymbol{\alpha}\in\mathcal S_q}
\widetilde f_{\boldsymbol{\alpha}}^N.
\]
We denote by $s_h^N(r)$ the cubic B-spline interpolant of these values. The derivative $(s_h^N)'(r)$ then approximates the radial derivative of the reconstructed distribution.

Let $\{(y_q,\omega_q^y)\}_{q=1}^{N_y}$ be a quadrature rule for the transverse radial variable in \eqref{eq:radial_G}. Its discrete counterpart is
\begin{equation*}
%\label{eq:discrete_radial_G}
G_h^N(x_\ell)
=
x_\ell S_{d-2}
\sum_{q=1}^{N_y}
\omega_q^y
\frac{
(s_h^N)'\left(\sqrt{x_\ell^2+y_q^2}\right)
}{
\sqrt{x_\ell^2+y_q^2}
}
y_q^{d-2}.
\end{equation*}
In contrast to the general case, $G_h^N$ is independent of the direction $\widehat k$ and therefore needs to be computed only once at each time step.

The profile is extended to a symmetric interval and zero-padded before applying the FFT-based discrete Hilbert transform, in order to reduce wrap-around errors. The corresponding approximation of \eqref{eq:dielectric_integral} is
\begin{equation*}
%\label{eq:discrete_radial_dielectric_integral}
I_{h,\mathrm{rad}}^N(k,v)
=
\pi H_h[G_h^N](v_\parallel)
+
\mathrm{i}\pi G_h^N(v_\parallel),
\qquad
v_\parallel=\frac{k\cdot v}{|k|}.
\end{equation*}
Thus, in the radially symmetric case,
\[
\varepsilon_{\epsilon,h}^N(k,k\cdot v)
=
1+\widehat V(k)I_{h,\mathrm{rad}}^N(k,v).
\]
The precomputed values of $G_h^N$ and $H_h[G_h^N]$ are subsequently interpolated at $v_\parallel$ during the evaluation of the Lenard--Balescu kernel.

\subsection{Numerical energy conservation}
In the fully discretized case, we need to further discretize the particle system \eqref{eq:particle_I} in time. Denote $t^n = n\Delta t$. For example, the explicit Euler method yields
\begin{equation}
\label{eq:forward_Euler}
\frac{v_i^{n+1}-v_i^n}{\Delta t}=-\sum_{j=1}^N w_j B_{\epsilon}^N(v_i^n, v_i^n-
 v_j^n)\left(\nabla \frac{\delta E_{\epsilon}^N}{\delta f}(v_i^n)- \nabla \frac{\delta E_{\epsilon}^N}{\delta f}(v_j^n)\right),\quad i=1,\cdots,N.
\end{equation}
Apparently, the total discrete mass $\sum_{i=1}^N w_i$, denoted by $M$, is always conserved. Simple computation shows that 
$$
\frac{\sum\limits_{i=1}^N w_i v_i^{n+1}- \sum\limits_{i=1}^N w_i v_i^n}{\Delta t}=-\sum_{i=1}^N\sum_{j=1}^N w_iw_j B_{\epsilon}^N(v_i^n,v_i^n- v_j^n)\left(\nabla \frac{\delta E_{\epsilon}^N}{\delta f}(v_i^n)- \nabla \frac{\delta E_{\epsilon}^N}{\delta f}(v_j^n)\right) = 0,
$$
which means the discrete momentum is conserved exactly. However, this does not hold for the discrete energy. Hence in the following, we apply the variance reduction technique to ensure energy conservation after time discretization.

Denote $\mu$ and $\sigma$ to be the target mean and variance determined by the initial data, i.e.
$$\mu = \frac{1}{M}\sum\limits_{i=1}^N w_i v_i^0, \quad \sigma = \frac{1}{M}\sum\limits_{i=1}^N w_i |v_i^0-\mu|^2.$$
At each time step, first do the forward Euler step as in \eqref{eq:forward_Euler} and denote the $i$-th particle velocity after this step as $v_i^*$. Then compute the current variance
$$\sigma^* = \frac{1}{M}\sum\limits_{i=1}^N w_i |v_i^*-\mu|^2. $$
and update the particles by
\begin{equation}
\label{eq:correction_step}
v_i^{n+1} = (v_i^*- \mu)\sqrt{\frac{\sigma}{\sigma^*}}+\mu,  \quad i=1,\cdots,N.
\end{equation}
Hence, we deduce 
$$
\frac{1}{M}\sum_{i=1}^N w_iv_i^{n+1} = \frac{1}{M}\sum_{i=1}^N w_i (v_i^*- \mu)\sqrt{\frac{\sigma}{\sigma^*}} + \mu = \mu,
$$
$$
\frac{1}{M}\sum_{i=1}^N w_i|v_i^{n+1}-\mu|^2 = \frac{1}{M}\sum_{i=1}^N w_i|v_i^*-\mu|^2 \frac{\sigma}{\sigma^*} =\sigma.
$$
This means that the mean and variance remain unchanged after the correction step \eqref{eq:correction_step}. It follows from $$\sum\limits_{i=1}^N w_i |v_i^{n+1}|^2 = \sum_{i=1}^N w_i|v_i^{n+1}-\mu|^2 - |\mu|^2\sum_{i=1}^N w_i + 2\mu \cdot\sum_{i=1}^N w_i v_i^{n+1} = (\sigma +|\mu|^2)M $$ 
that the discrete energy is also conserved.

\subsection{Choice of the mollifier and cell-list optimization}
\label{subsec:mollifier}
The choice of the mollifier affects both the regularisation error and the computational cost of the particle method. We consider either a Gaussian mollifier or a compactly supported tensor-product B-spline. Using $\epsilon$ as a length scale, the Gaussian mollifier is
\begin{equation*}
%\label{eq:gausskernel}
\psi_\epsilon^{\rm G}(v)
=
\frac{1}{(2\pi\epsilon^2)^{d/2}}
\exp\left(-\frac{|v|^2}{2\epsilon^2}\right).
\end{equation*}
This choice is smooth but has full support. Consequently, its exact evaluation requires contributions from all particle pairs. A finite cutoff may be introduced in practice, but this adds a further approximation and an additional numerical parameter.

For the compactly supported case, following \cite{bailo2024collisional}, we consider
\begin{equation*}
%\label{eq:scaled_spline}
S_{r,\epsilon}(v)
=
\frac{1}{\epsilon^d}
\prod_{\ell=1}^d
S_r\left(\frac{v_\ell}{\epsilon}\right),
\qquad r=0,1,2,
\end{equation*}
where
\[
S_0(z)=(1-|z|)_+,
\]
\[
S_1(z)
=
\begin{cases}
\dfrac{3}{4}-|z|^2,
& |z|\leq \dfrac{1}{2},\\[1ex]
\dfrac{1}{2}\left(\dfrac{3}{2}-|z|\right)^2,
& \dfrac{1}{2}<|z|\leq \dfrac{3}{2},\\[1ex]
0,
& |z|>\dfrac{3}{2},
\end{cases}
\]
and
\[
S_2(z)
=
\begin{cases}
\dfrac{1}{6}\left(4-6|z|^2+3|z|^3\right),
& |z|\leq 1,\\[1ex]
\dfrac{1}{6}\left(2-|z|\right)^3,
& 1<|z|\leq 2,\\[1ex]
0,
& |z|>2.
\end{cases}
\]
These choices correspond, respectively, to linear, quadratic, and cubic B-splines.

For the spline mollifiers, we link the regularisation scale to the velocity-grid resolution by taking
\[
\epsilon=h,
\]
where $h$ is the mesh size introduced in Subsection~5.1. Their supports are therefore
\[
\operatorname{supp}S_{0,h}=[-h,h]^d,
\qquad
\operatorname{supp}S_{1,h}
=
\left[-\frac{3h}{2},\frac{3h}{2}\right]^d,
\qquad
\operatorname{supp}S_{2,h}=[-2h,2h]^d.
\]
Thus, independently of the grid resolution, each particle interacts through the mollifier with particles contained in only a fixed number of neighbouring cells.

This compact support permits an exact cell-list optimization. At every time step, each particle is assigned to the cell $C_{\boldsymbol{\alpha}}$ containing its velocity. When evaluating the mollifier-based sums for a particle $v_i$, we inspect only the cell containing $v_i$ and the neighbouring cells that intersect
\[
v_i+\operatorname{supp}S_{r,h}.
\]
Particles in all other cells give exactly zero contribution and can therefore be discarded. The cell-list procedure does not modify the particle method; it only avoids evaluating particle pairs outside the support of the mollifier. If the number of particles in each cell remains uniformly bounded, the mollifier evaluation costs $\cO(N)$. This exact sparsity is not available for the Gaussian mollifier unless it is artificially truncated.

\subsection{Random batch particle method}
Recall the particle method
$$
\frac{\mathrm{d}v_i}{\mathrm{d}t}=-\sum_{j=1}^N w_j B_{\epsilon}^N(v_i,v_i-v_j)\left(\nabla \frac{\delta E_{\epsilon}^N}{\delta f}(v_i)- \nabla \frac{\delta E_{\epsilon}^N}{\delta f}(v_j)\right), \quad i=1,\cdots,N.
 $$
 At each time step, the cost is $\cO(N^2)$ due to pairwise particle interactions. Hence, we adopt the random batch idea \cite{JinLiLiu20,CarrilloJinTang} to reduce the computational cost.

 For simplicity, suppose $p$ divides $N$. At each time step, randomly divide the $N$ particles into $N/p$ batches of size $p$. Interactions then take place only inside each batch, namely,
\begin{equation*}
\frac{\mathrm{d}v_i}{\mathrm{d}t}=-\frac{N-1}{p-1}\sum_{j\in \C_q} w_j B_{\epsilon}^N(v_i,v_i-v_j)\left(\nabla \frac{\delta E_{\epsilon}^N}{\delta f}(v_i)- \nabla \frac{\delta E_{\epsilon}^N}{\delta f}(v_j)\right),\quad  i\in \C_q,\ q=1,\cdots,N/p.
%\label{eq:RBM_update}
\end{equation*}
In the implementation, the requested batch size determines a balanced partition. For a batch $\C_q$ of actual size $p_q$, the normalization is $(N-1)/(p_q-1)$.
In this case, the cost inside each small batch is $\cO(p^2)$ and since there are $N/p$ batches the overall cost is $\cO(pN)$. Furthermore, the random batch version still inherits the structure properties of the regularized Lenard--Balescu equation. This is due to fact that the random batch sum can be rewritten into a full batch sum with a symmetric indicator which identifies whether two particles are in the same batch, and the proof only relies on the structure properties of $B_{\epsilon}^N$ as mentioned in Remark \ref{remark:3.2}, see \cite{CarrilloJinTang} for details.

\section{Numerical simulations}
\label{sec:6}

In this section, we perform some numerical simulations to test our particle method. Firstly, we show the validity of the particle method in the simple static case by comparing the solution of the static Lenard--Balescu equation to the Landau equation in section \ref{subsec:6.1}. Then we test our method in general cases when the dielectric function is computed as mentioned in Section \ref{sec:4}, including the radially symmetric case (section \ref{subsec:6.2}) and the fully dynamic screening case (section \ref{subsec:6.4}). The solutions relax towards the limiting Maxwellians. Furthermore, the evolutions of the macroscopic quantities and the entropy functional show that the proposed method is structure preserving as expected. Last but not least, we also give convergence tests with respect to the velocity-grid resolution $n_v$ and the particle number $N$, observing approximately second-order and first-order rates, respectively, over the tested ranges in section \ref{subsec:6.3}.

In the following, we use two quantities to measure the `distance' between the simulated distribution $f(v)$ and a reference distribution $g(v)$. One is the `relative entropy' (also known as the Kullback–Leibler divergence) defined as
\[
S[f \| g] = \int f(v) \log\left( \frac{f(v)}{g(v)} \right) \, \mathrm{d}v.
\]
This quantity measures the `distance' between $f$ and $g$ in an information-theoretic sense (though note that it is not a metric). For distributions of equal mass, this quantity is nonnegative and vanishes if and only if $f=g$ almost everywhere. We will also discuss `relative error' in the form of the relative 2-distance, i.e.
\begin{align*}
    \text{err}_{\text{rel}}(f,g) = \frac{\| f(v) - g(v) \|_2}{\| g(v) \|_2}.
\end{align*}
For instance, in the tests where we check the relaxation towards the equilibrium, $g$ is chosen as the limiting Maxwellian $$\mathcal{M}_{\rho,u,T}(v)= \frac{\rho}{(2\pi T)^{d/2}} \exp\left(-\frac{|v-u|^2}{2T}\right)$$ with $\rho,u,T$ given by the initial data.

Unless otherwise stated, we use the non-symmetric entropy regularization and the cubic B spline $S_{2,\epsilon}$, with $\epsilon=h=2L_v/n_v$. The particle system is advanced using the explicit Euler method, followed by the energy correction~\eqref{eq:correction_step} at each time step. The velocity domain, particle number, time step, Fourier quadrature parameters, and batch size are specified for each experiment.

The companion Julia code is available under the open MIT license as LenardBalescu.jl~\cite{gutleb2026lenardbalescu}.

\subsection{Convergence of static Lenard--Balescu to Landau}
\label{subsec:6.1}

As mentioned in the introduction, when there is no screening effect, i.e., $\alpha=0$, and $\widehat V$ is isotropic, the static Lenard--Balescu kernel \eqref{eq:static_kernel} reduces to a Landau-type kernel \eqref{eq:LB_to_Landau}. Motivated by this, in this subsection, we investigate the convergence of the static Lenard--Balescu equation to the Landau equation in the vanishing screening limit $\alpha\to 0$. This special case can thus act as a nice sanity check of any Lenard--Balescu implementation as it can use essentially the same numerical framework as the dynamic Lenard--Balescu implementation, only differing in the dielectric function used.

Recall the simplified setting with static kernel \eqref{eq:static_kernel}.
More precisely, when $\alpha = 0$, one obtains the Landau-type kernel with Coulomb interactions \cite{duerinckx2023}
$$
B_{L}(v-v_*)
=
\frac{L}{|v-v_*|}
\left(
I_d-
\frac{(v-v_*)\otimes(v-v_*)}{|v-v_*|^2}
\right),
$$
where
$$
L
=
\frac{S_{d-2}}{(d-1)S_{d-1}}
\int_{\mathbb R^d}
|k|\,\pi \widehat V^2(k)\,\mathrm dk
=
\frac{S_{d-2}}{d-1}
\int_0^\infty
r^d\pi \widehat V^2(r)\,\mathrm dr.
$$
Note that this identity holds only when
$L<\infty$, which is equivalent to
$V\in\dot H^{1/2}(\mathbb R^d)$.

For the genuine Coulomb interaction, 
$$
\widehat V(k)= |k|^{-2},
$$
however, the above integral is divergent. Indeed,
$L\propto \int_0^\infty r^{d-4}\,\mathrm dr$. 
In $2$D,
$ L$ diverges at $0$, whereas in $3$D,
$L$ is logarithmically divergent both as $r\to0$ and as $r\to\infty$. Thus, for the true Coulomb potential, the formal limit $\alpha=0$ does not directly define a finite Landau coefficient $L$.

For this reason, a finite coefficient can be defined by restricting the $k$-integral to
$$
k_{\min}\le |k|\le k_{\max}.
$$
Accordingly, we define the cutoff coefficient
$$
L_{\text{cutoff}}
=
\frac{S_{d-2}}{d-1}
\int_{k_{\min}}^{k_{\max}}
r^d\pi\widehat V^2(r)\,\mathrm dr,
$$
which is always finite.

In the implementation, we instead calibrate the Landau coefficient from the same finite Fourier quadrature used for the static Lenard--Balescu kernel at $\alpha=0$, rather than evaluating $L_{\text{cutoff}}$ analytically. The quadrature avoids $k=0$; its smallest nonzero node acts as a numerical lower cutoff, without an explicitly imposed $k_{\min}$. Keeping this quadrature fixed gives a consistent comparison as $\alpha\to0$.

\begin{figure}[!ht]
  \centering

  % Optional tightening
  \vspace{-2mm}

  \begin{tabular}{@{}ccc@{}}
    % -------- Trajectory row --------
    \includegraphics[width=0.47\textwidth]{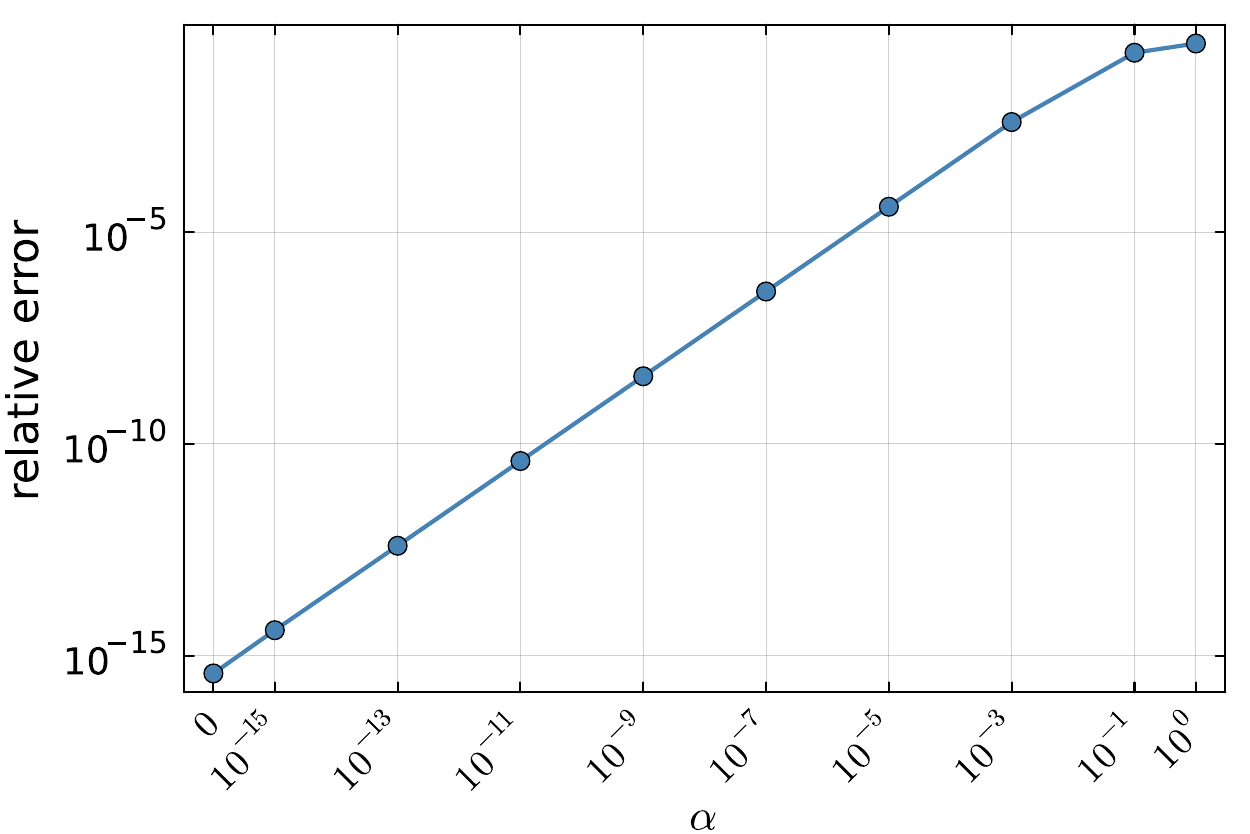} &
    \includegraphics[width=0.47\textwidth]{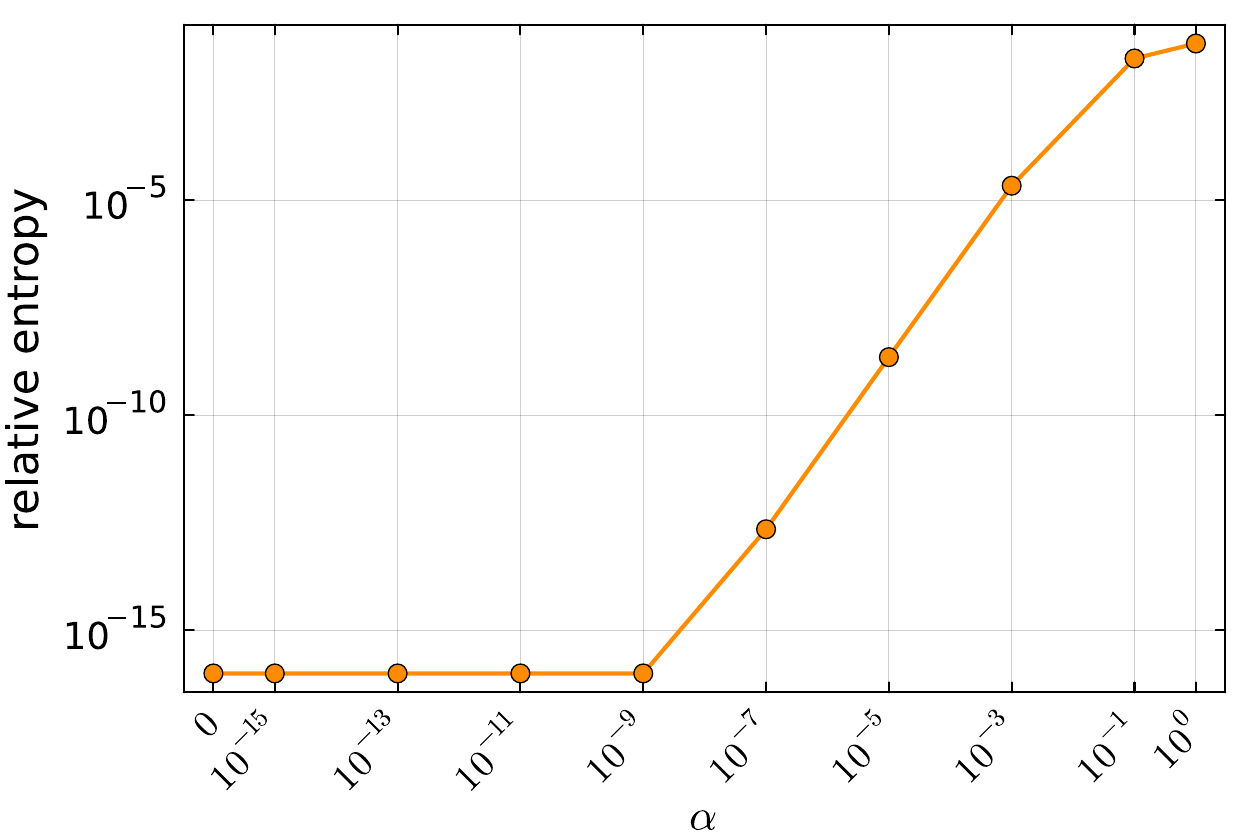} \\

    % -------- Endpoint row --------
    \includegraphics[width=0.35\textwidth]{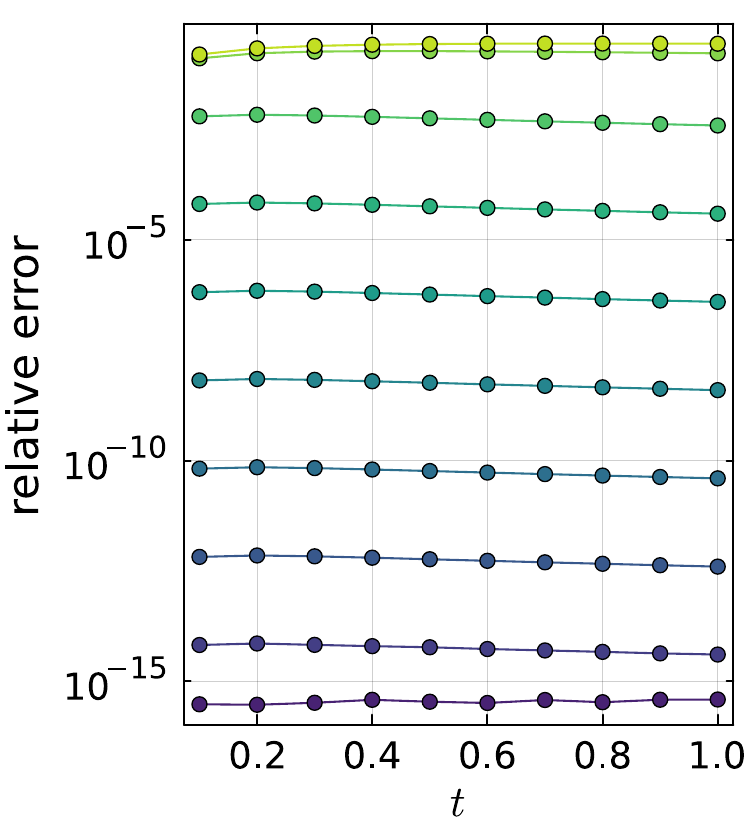} &
    \includegraphics[width=0.46\textwidth]{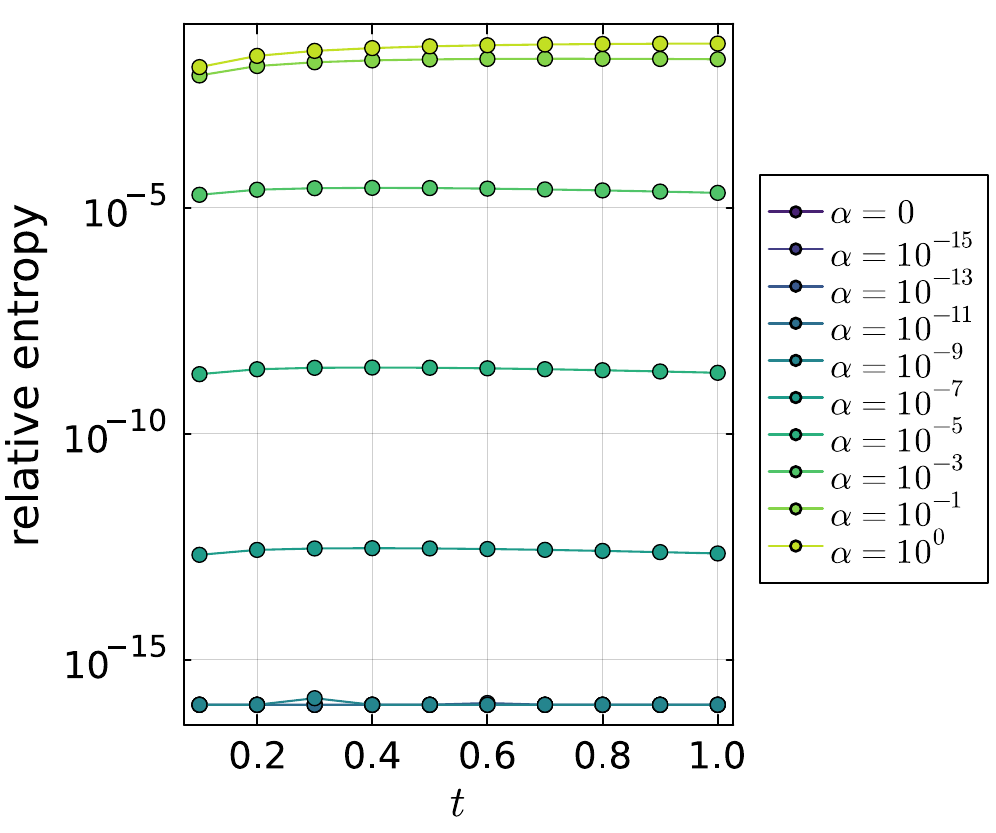} 
  \end{tabular}

  \vspace{-2mm}

  \caption{Convergence of the static Lenard--Balescu solution to the Landau solution as $\alpha\to0$. Top: final relative $L^2$ error (left) and relative entropy (right) versus $\alpha$. Bottom: their evolution in time.}
  \label{fig:statictolandau1}
\end{figure}

The numerical results in Figure \ref{fig:statictolandau1} confirm this convergence.
 We show a representative numerical experiment with Coulomb interactions starting from the $2$D BKW distribution \eqref{eq:BKW} with $a=0$, $b=1$, and $T=0.5$.
 We use $[-6,6]^2$, $n_v=40$, $N_{\rm pc}=3$, and $N=4800$, with full interactions, $\Delta t=0.01$, and final time $T_{\rm end}=1$. The sinh-mapped transverse quadrature has $n_k=50$, $k_{\perp,\max}=100$, and $k_{\rm scale}=1$; the parallel quadrature has width $0.02$ and order $3$. The smallest active wave number is approximately $0.165516$, giving a matched Landau coefficient $L_{\rm num}\approx2.992488$ in the implementation's time normalization. As $\alpha$ decreases, we can see the solution of static Lenard--Balescu converges to solution of Landau for the final time (Figure \ref{fig:statictolandau1} Top) as well as through evolution (Figure \ref{fig:statictolandau1} Bottom), as demonstrated by both the relative $L^2$ and the relative entropy measures.

\subsection{BKW type initial states}
\label{subsec:6.2}
We consider the Lenard--Balescu equation with BKW type initial distributions \eqref{eq:BKW} with the simplest choice of $\widehat{V}$, i.e., $\widehat{V} = 1$, as a benchmark test. We use these smooth, radially symmetric profiles as they provide a clean setting to assess relaxation towards equilibrium as well as the numerical conservation of mass, momentum, and energy. 

\begin{figure}[!ht]
  \centering

  % Optional tightening
  \vspace{-2mm}

  \begin{tabular}{@{}c@{\hspace{0.01\textwidth}}c@{\hspace{0.01\textwidth}}c@{}c@{}}
    $t=0.0$ & $t=0.25$ & $t=5$ & \\[1mm]
    % -------- Scatter row --------
    \includegraphics[width=0.30\textwidth]{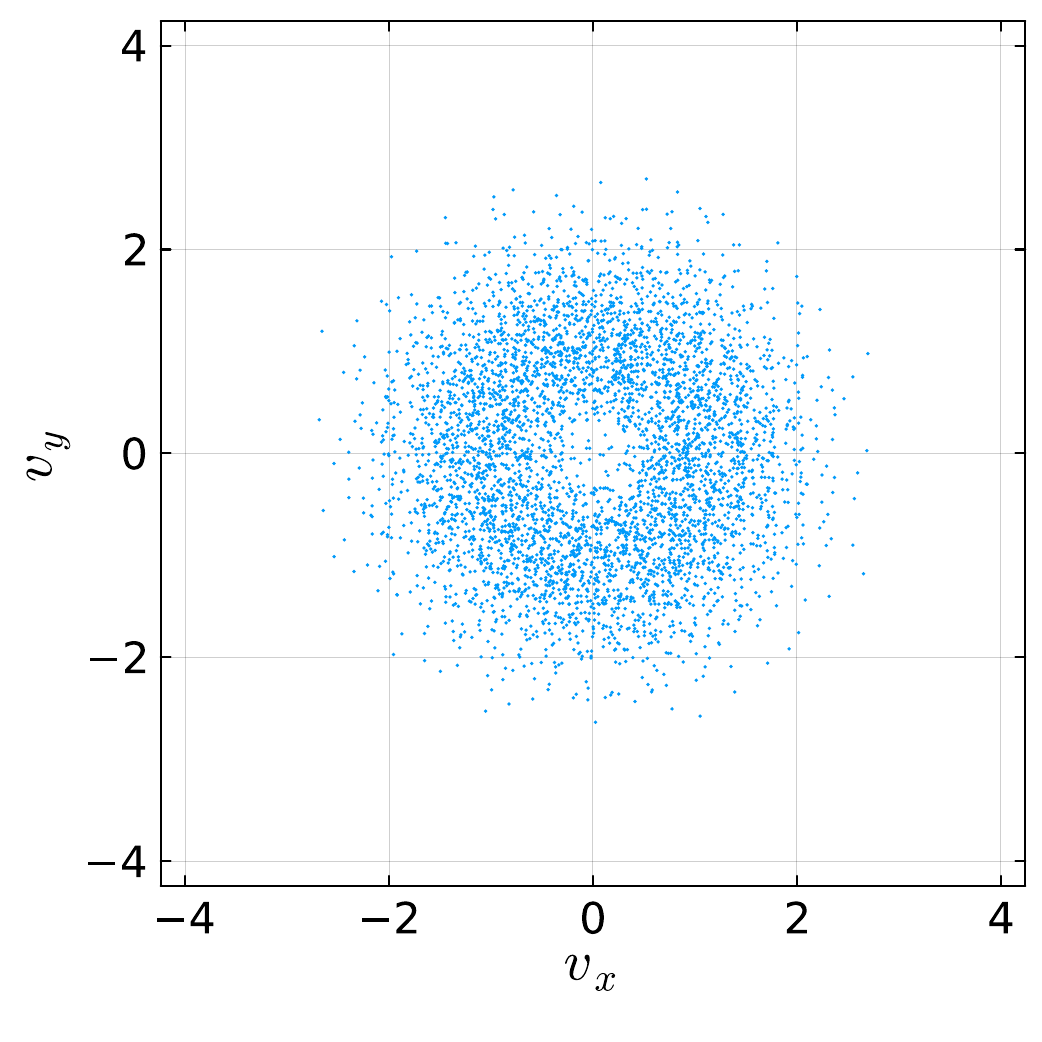} &
    \includegraphics[width=0.30\textwidth]{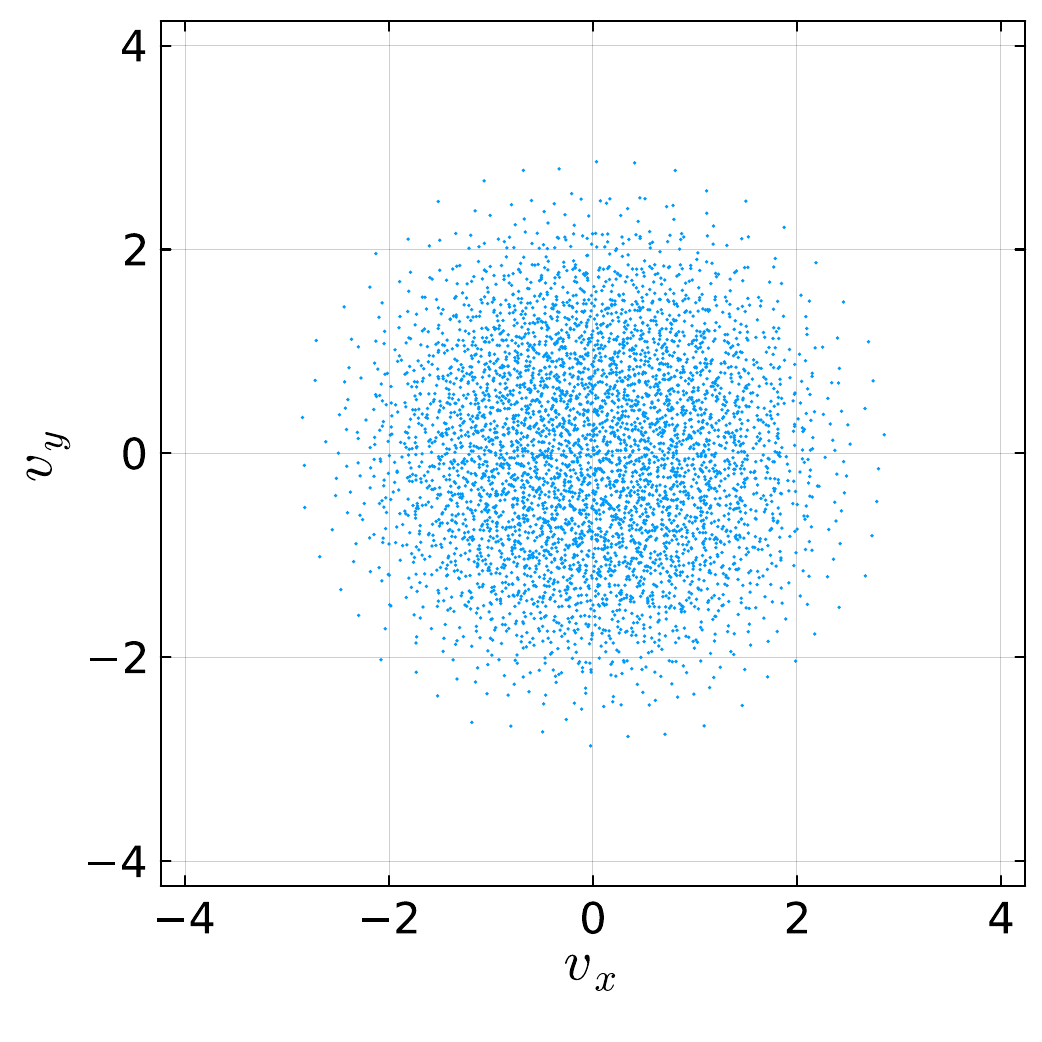} &
    \includegraphics[width=0.30\textwidth]{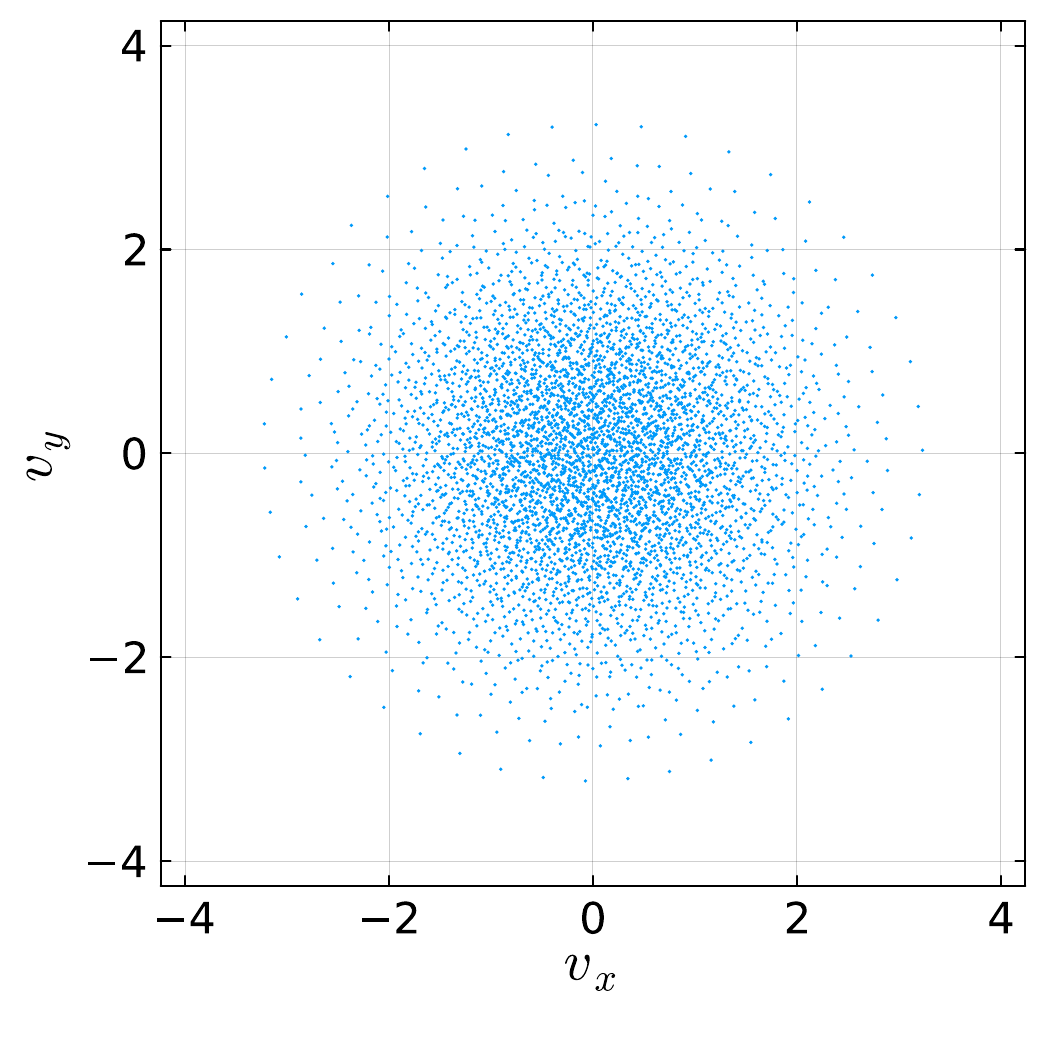} & \\[-2mm]

    % -------- Heatmap row --------
    % Match the axes boxes to the scatter panels, with a separate shared colorbar.
    \includegraphics[pagebox=mediabox,viewport=49 -21 573 503,clip,width=0.30\textwidth]{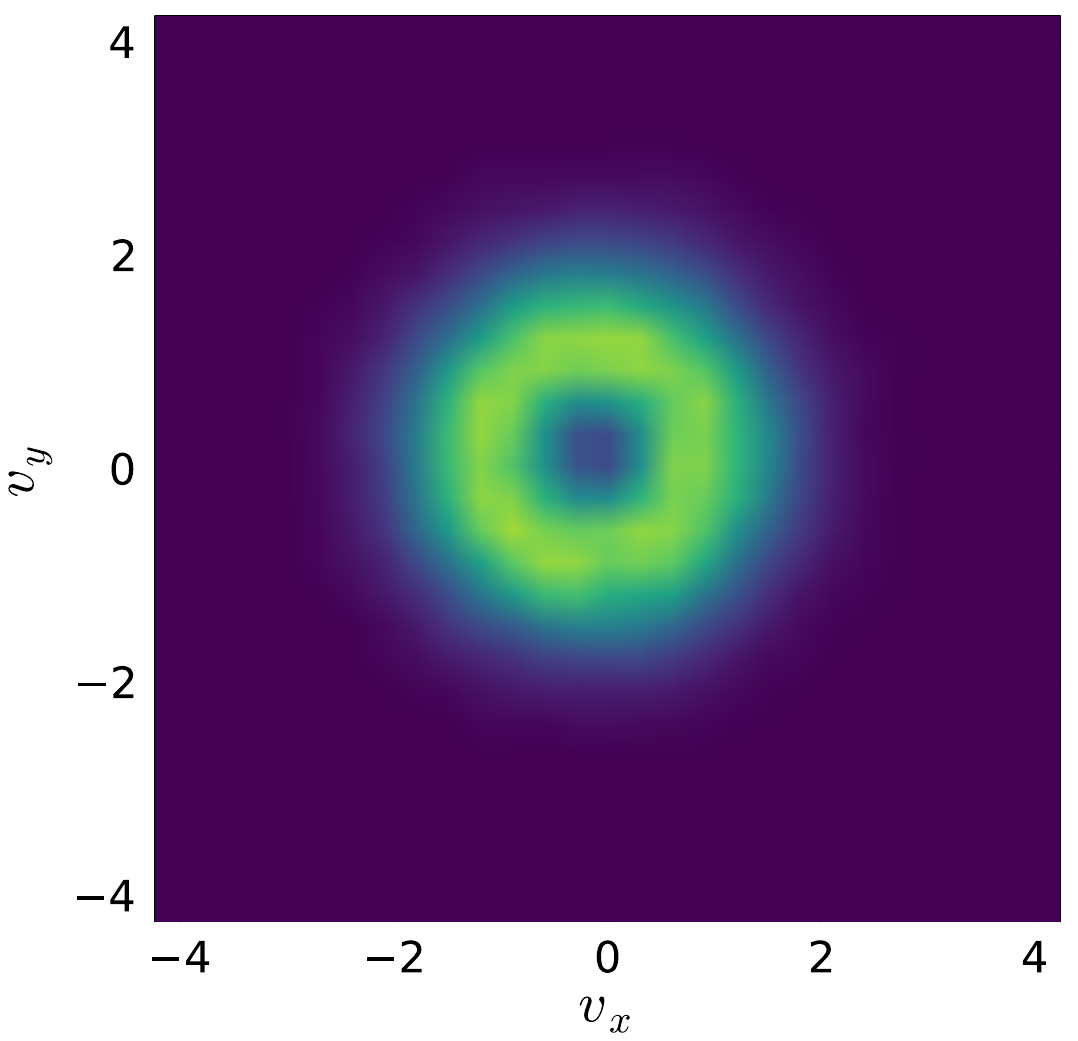} &
    \includegraphics[pagebox=mediabox,viewport=49 -21 573 503,clip,width=0.30\textwidth]{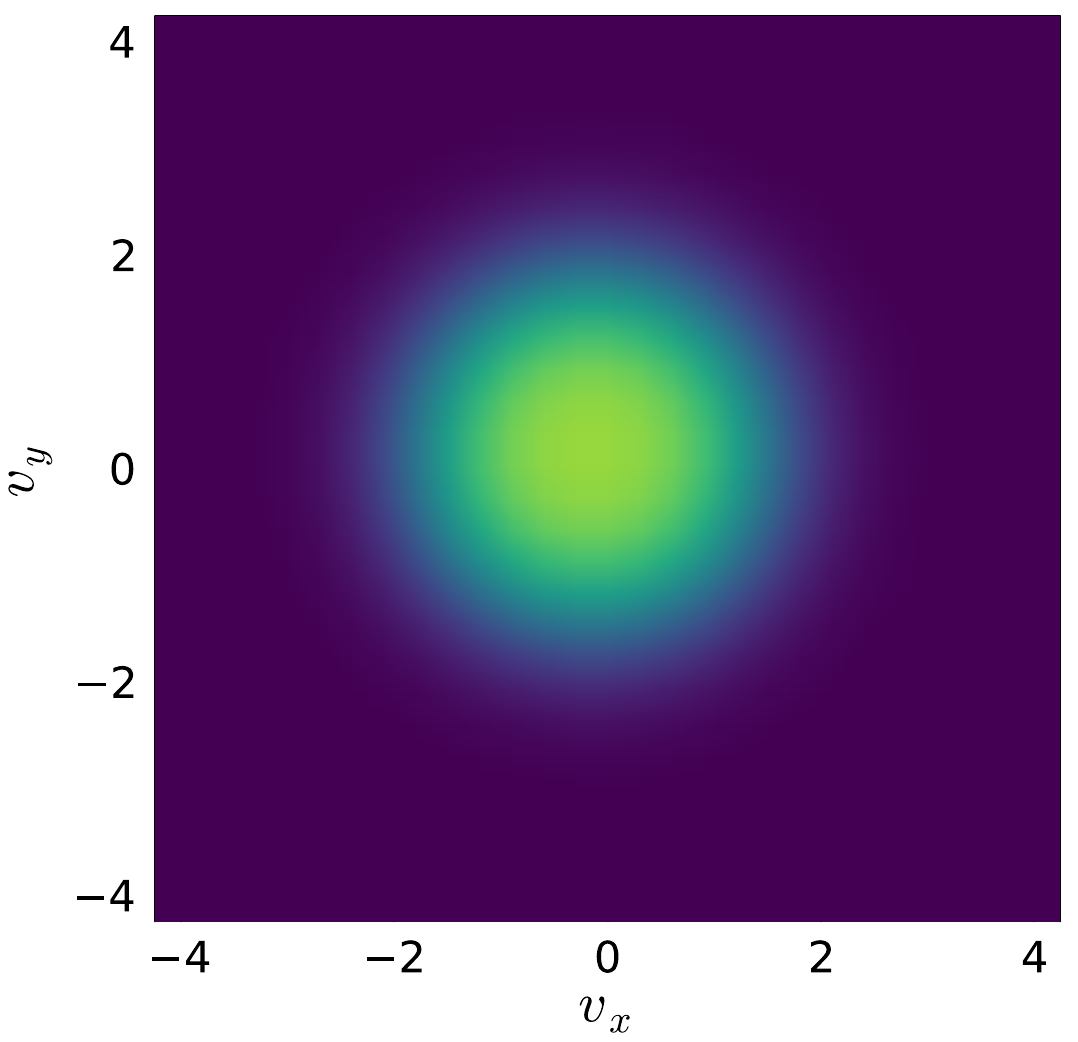} &
    \includegraphics[pagebox=mediabox,viewport=49 -21 573 503,clip,width=0.30\textwidth]{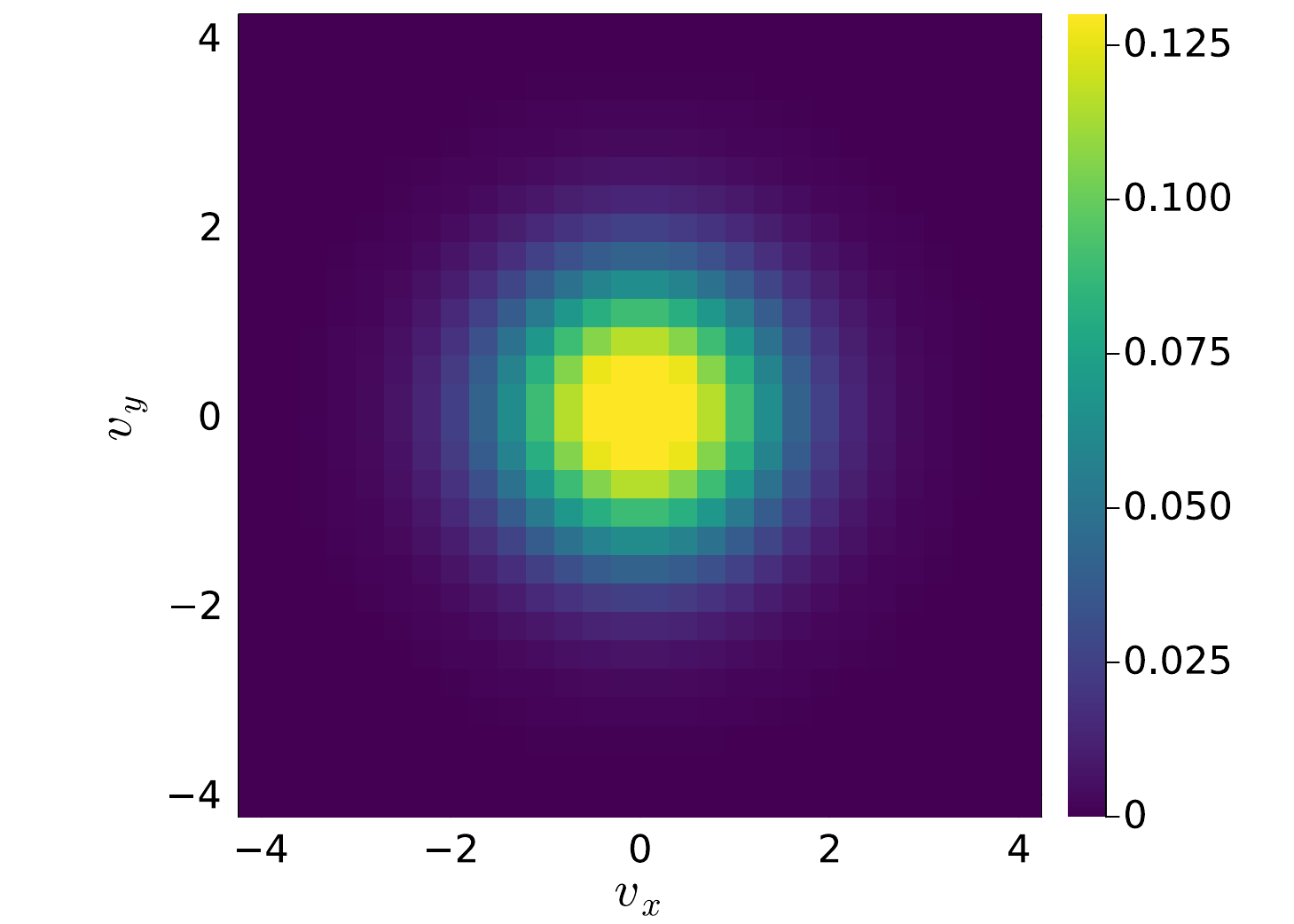} &
    \includegraphics[pagebox=mediabox,viewport=573 -21 675 503,clip,height=0.30\textwidth]{2dbkw_heatmap_final.pdf}
  \end{tabular}

  \vspace{-2mm}

 \caption{Evolution of the two-dimensional BKW initial state. Particle velocities (top) and reconstructed density (bottom) at $t=0$, $0.25$, and $5$ from left to right.}
  \label{fig:2dbkw_main}
\end{figure}

\begin{figure}[!ht]
  \centering

  % Optional tightening
  \vspace{-2mm}

  \begin{tabular}{@{}c@{\hspace{0.01\textwidth}}c@{\hspace{0.01\textwidth}}c@{}}
    % -------- Scatter row --------
    \includegraphics[trim=80bp 40bp 70bp 0bp,clip,width=0.32\textwidth]{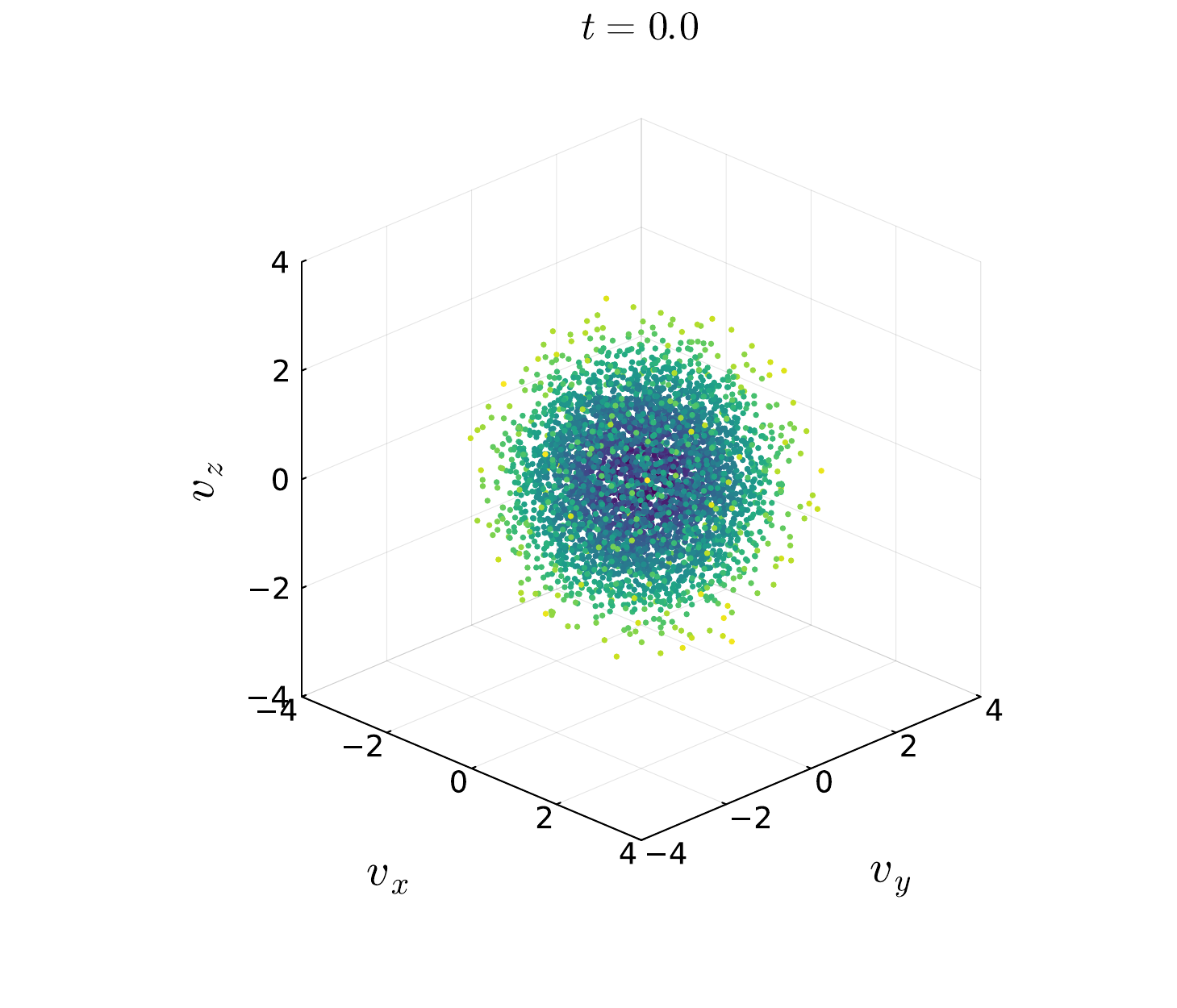} &
    \includegraphics[trim=80bp 40bp 70bp 0bp,clip,width=0.32\textwidth]{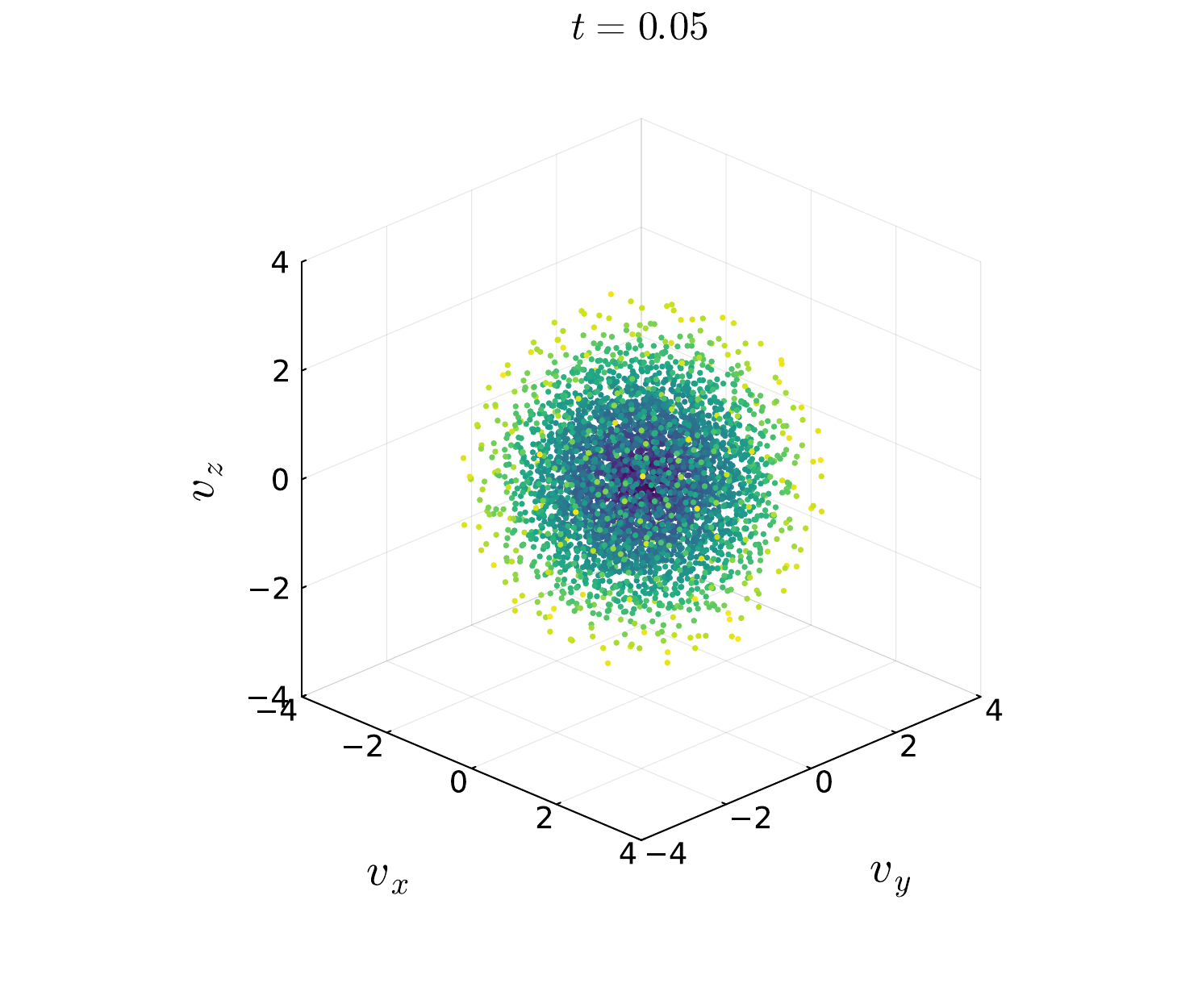} &
    \includegraphics[trim=80bp 40bp 70bp 0bp,clip,width=0.32\textwidth]{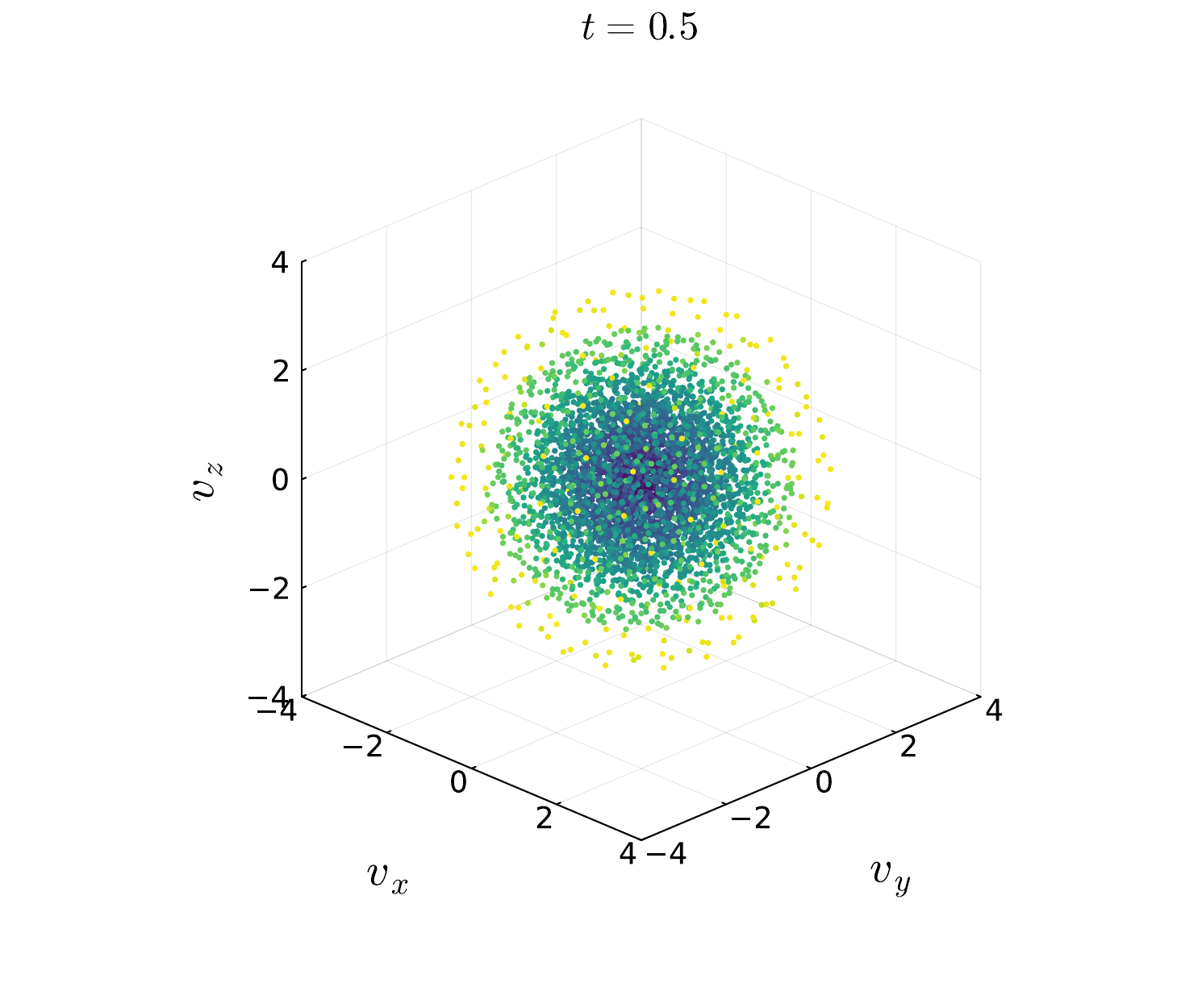}
  \end{tabular}

  \vspace{-2mm}

\caption{Evolution of the three-dimensional BKW initial state at $t=0$, $0.05$, and $0.5$ from left to right. Particles are coloured by speed.}

  \label{fig:3dbkw}
\end{figure}

Figures~\ref{fig:2dbkw_main} and \ref{fig:3dbkw} illustrate the time evolution of BKW type initial data in two and three velocity dimensions. The initial data are \eqref{eq:BKW} with $a=0$, $b=1$, and $T=0.5$ in $2$D, 
and \eqref{eq:BKW} with $a=0$, $b=5/9$, and $T=0.6$ in $3$D. This is the $t=0$ profile of the BKW solution \cite[Eq.~(35)]{MEDAGLIA2024}, with its equilibrium temperature set to one.

Both tests use the full dynamical dielectric response and stratified initial sampling as described in Subsection~\ref{subsec:initialization}. For the two-dimensional experiment, we use $n_v=40$ cells per velocity direction on $[-6,6]^2$, with an average $N_{\rm pc}=3$ particles per cell, giving $N=4800$. The time step is $\Delta t=10^{-2}$ and the final time is $T_{\rm end}=5$, with full particle interactions. The transverse Fourier quadrature uses $n_k=50$ nodes, $k_{\perp,\max}=5$, and $k_{\rm scale}=1$.

\begin{figure}[!ht]
    \centering
    \includegraphics[width=\textwidth]{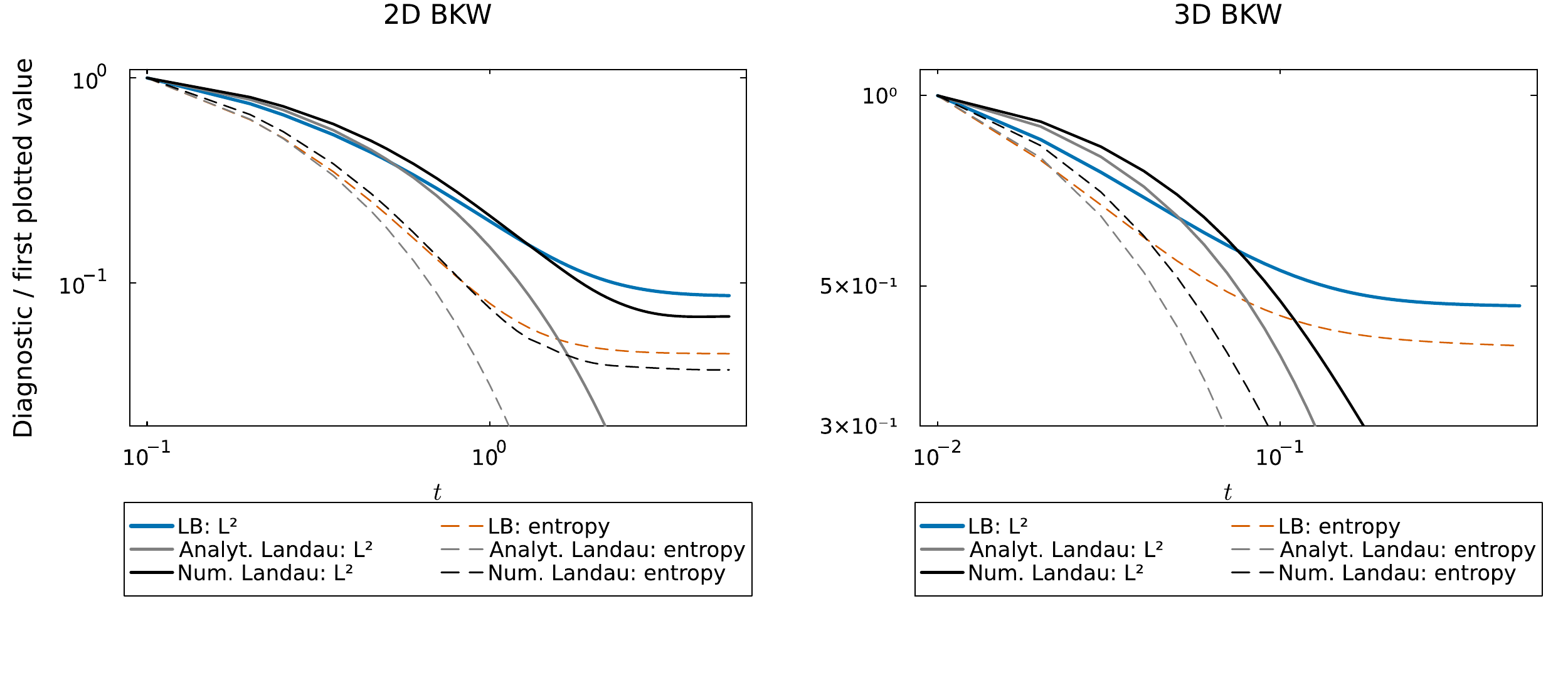}
    \caption{Relaxation of BKW initial states with $\widehat{V}=1$ interactions in two dimensions (left) and three dimensions (right). Each panel shows the relative $L^2$ distance to the matching Maxwellian (solid) and relative entropy (dashed) on log--log axes, omitting $t=0$. Each diagnostic is divided by its value at the first plotted positive time to compare decay rates despite the different initial particle and analytic reconstructions. Coloured curves show Lenard--Balescu simulations; black curves show the analytic Landau--BKW family~\cite{CARRILLO2020} with kernel $A(q)=B(|q|^2I-q\otimes q)$ with $B$ fitted to match the initial LB rate.}
    \label{fig:bkw_maxwell_relaxation}
\end{figure}

For the three-dimensional LB experiment, we use $n_v=13$ on $[-4,4]^3$, with $N_{\rm pc}=2$, giving $N=4394$ particles. The time step is $\Delta t=10^{-3}$ and the final time is $T_{\rm end}=0.5$. The transverse quadrature uses $n_k=30$ nodes per direction, $k_{\perp,\max}=5$, and $k_{\rm scale}=1$. This run uses random batching with target batch size $p=1000$.
Figure~\ref{fig:bkw_maxwell_relaxation} shows the convergence to the Maxwellian in terms of relative entropy as well as relative $L^2$ error. As it can be observed, the Lenard-Balescu relaxation is faster than the Landau equation with $\widehat{V}=1$ whenever the numerical results are reliable, i.e., only for small time before the plateau in the decay for Lenard-Balescu results.

Animations of these processes are available on Figshare~\cite{LenardBalescuAnimations}. Mass, energy and momentum remained approximately conserved throughout the simulations.

Our implementation further allows radial symmetry to be assumed, improving performance via radial bucketing of the particles.

The fitted three-dimensional Landau coefficient depends on the early-time window: the same procedure gives $B\approx1.15$ over $0.01\leq t\leq0.03$ and $B\approx0.56$ over $0.01\leq t\leq0.1$. It is a local rate match, not an identification of the Landau and Lenard--Balescu kernels.

\begin{figure}[!ht]
    \centering
    \includegraphics[width=\textwidth]{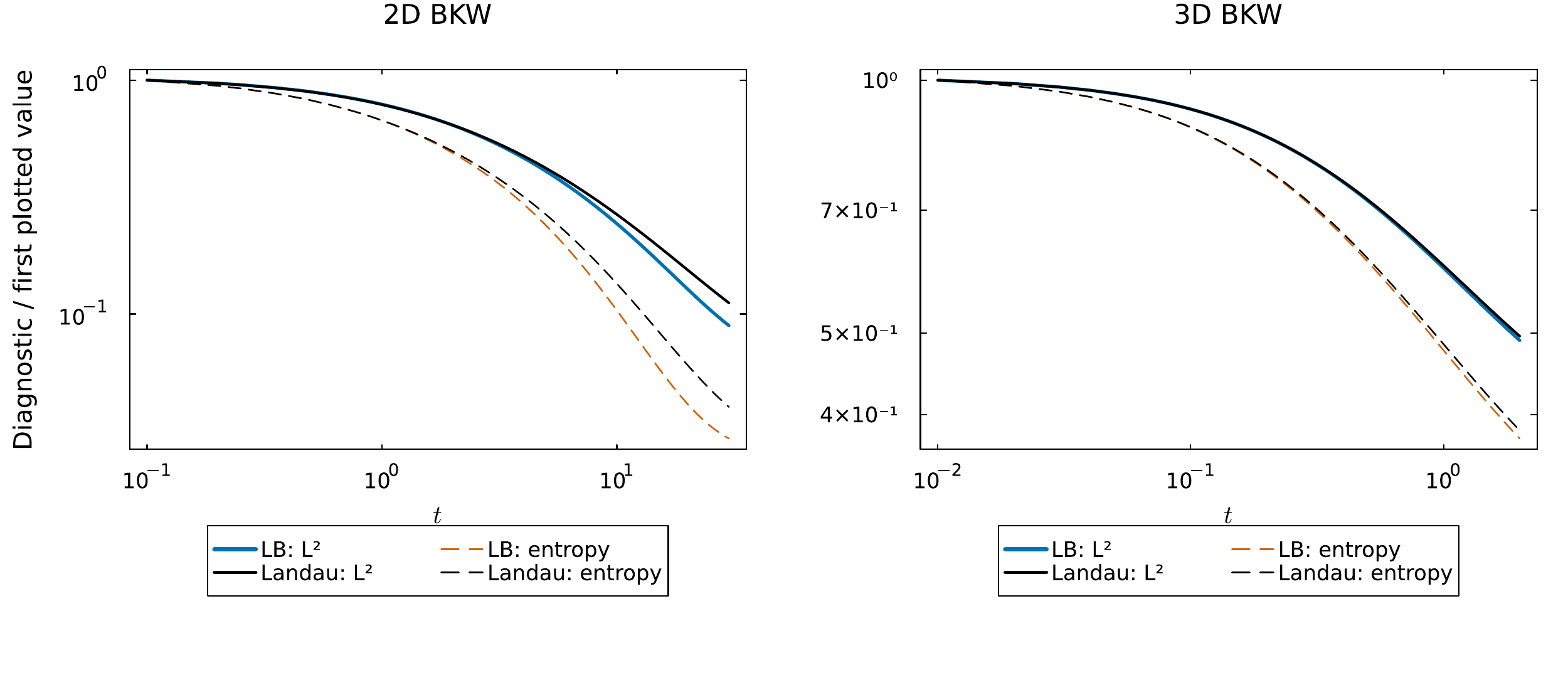}
    \caption{Relaxation of BKW initial states with Coulomb interactions in two dimensions (left, up to $t=30$) and three dimensions (right, up to $t=2$), showing the relative $L^2$ distance to the matching Maxwellian (solid) and relative entropy (dashed) on log--log axes, omitting $t=0$. Each diagnostic is divided by its value at the first plotted positive time. Coloured curves are Lenard--Balescu; black curves are numerical Coulomb Landau controls started from the same particles. Their coefficients $B$ are calibrated to the early decay rates of LB to make comparison possible.}
    \label{fig:bkw_coulomb_relaxation}
\end{figure}

We additionally consider Coulomb interactions, $\widehat V(k)=|k|^{-2}$, in Figure~\ref{fig:bkw_coulomb_relaxation}. The initial particles, velocity grids, time steps, and batching are matched to the runs with $\widehat{V}=1$ within each dimension; the Coulomb runs use $k_{\perp,\max}=100$ in two dimensions and $50$ in three dimensions, with the remaining quadrature parameters unchanged. Coulomb relaxation is slower over the simulated intervals. This comparison uses the stated finite Fourier cutoffs, which also affect the collision time scale.

\subsection{Convergence with respect to velocity-grid resolution and particle number}
\label{subsec:6.3}

We next examine convergence with respect to the number $n_v$ of velocity cells in each coordinate direction. Since the mollification scale is tied to the mesh width,
\[
\epsilon=h=\frac{2L_v}{n_v},
\]
this test measures the joint refinement of the velocity grid and the regularized particle reconstruction. We use the two-dimensional BKW initial condition from the previous subsection with $\widehat{V}=1$ interactions. The particle number is fixed at $N=5000$ and the time step at $\Delta t=0.01$, while
\[
n_v\in\{16,20,24,28,32,40,48\}.
\]
A full-interaction computation with $n_v^{\mathrm{ref}}=64$ is used as the reference solution. The velocity domain is $[-6,6]^2$, and the Fourier quadrature parameters are the same as in the two-dimensional BKW experiment: $n_k=50$, $k_{\perp,\max}=5$, and $k_{\rm scale}=1$.

To measure convergence during the non-equilibrium part of the evolution, we define
\begin{equation*}
\mathcal{E}_{\mathrm{self}}(n_v)
=
\frac{
\left\|f_{n_v}(T_*)-f_{\mathrm{ref}}(T_*)\right\|_{L^2}
}{
\left\|f_{\mathrm{ref}}(T_*)\right\|_{L^2}
},
\qquad T_*=0.2,
\end{equation*}
where $T_*$ is chosen approximately halfway through the relaxation towards equilibrium. We also report the finite-time distance from the Maxwellian,
\begin{equation*}
\mathcal{E}_{\mathrm{equil}}(n_v)
=
\frac{
\left\|f_{n_v}(T_{\mathrm{end}})-\mathcal{M}\right\|_{L^2}
}{
\left\|\mathcal{M}\right\|_{L^2}
},
\qquad T_{\mathrm{end}}=2.
\end{equation*}
All distributions are reconstructed on the reference grid before computing these errors.

\begin{figure}[!ht]
    \centering
    \includegraphics[width=0.78\textwidth]{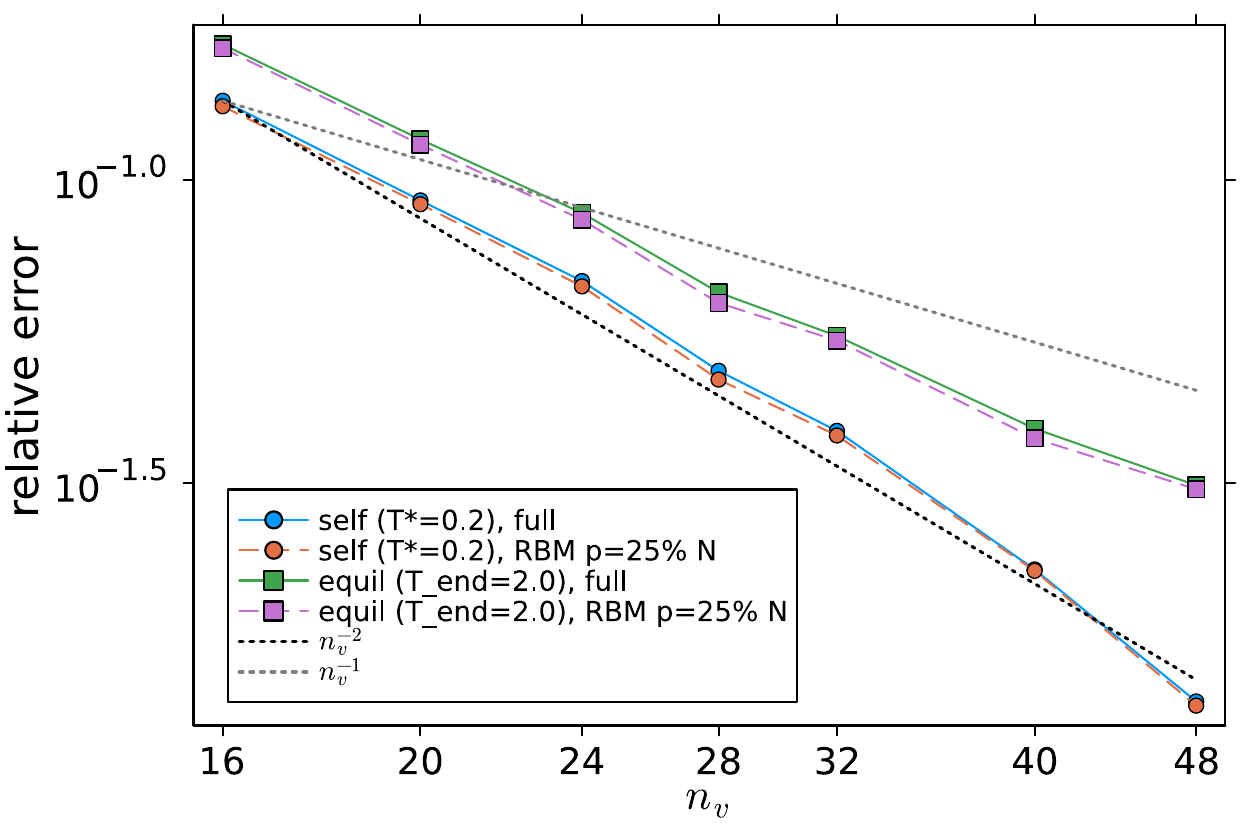}
    \caption{Velocity-grid convergence for the two-dimensional BKW test ($\log$-$\log$): relative errors versus $n_v$ at fixed $N=5000$, for full interactions and random batching, with first- and second-order reference slopes.}
    \label{fig:nv_convergence}
\end{figure}

Figure~\ref{fig:nv_convergence} shows approximately second-order self-convergence over the tested range. A least-squares fit gives rates $2.06$ for the full-interaction method and $2.05$ for the random-batch method. The latter uses batches of size $p=0.25N$ and produces errors almost indistinguishable from those of the full method, indicating that the random-batch error is smaller than the grid error in this regime. The error relative to the Maxwellian decreases more slowly, with an observed rate of approximately $1.55$, since at the finite final time it also contains the remaining relaxation error.

Since $h=2L_v/n_v$, the observed $n_v^{-2}$ decay corresponds to approximately $\cO(h^2)$ self-convergence. This is qualitatively consistent with the second-order convergence in $h$ observed for the regularized Landau particle method in \cite{CARRILLO2020}, although that study refines the particle number together with the mesh, whereas the present grid test keeps $N$ fixed.

\begin{figure}[!ht]
    \centering
    \includegraphics[width=0.78\textwidth]{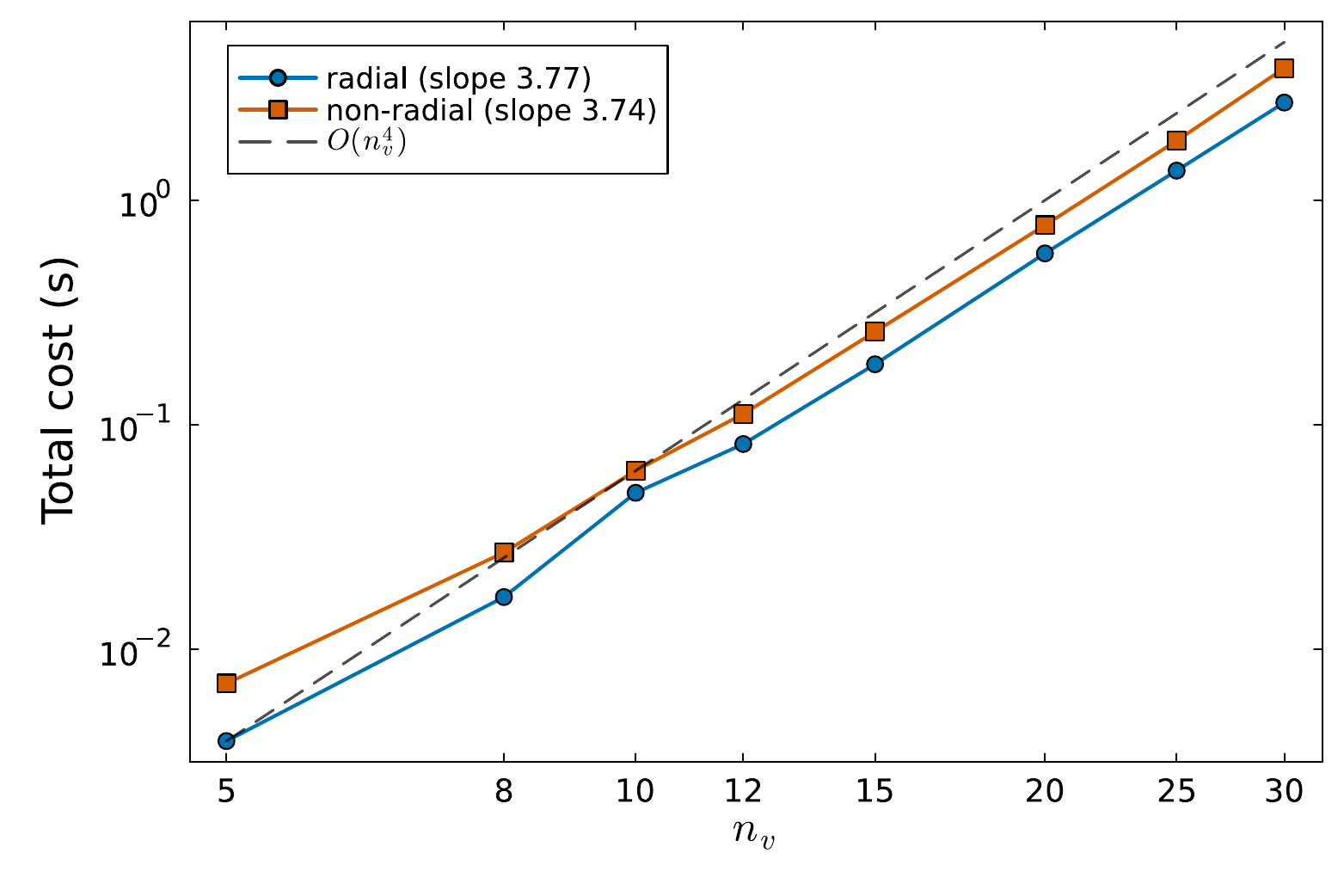}
    \caption{Computational cost ($\log$-$\log$) under joint grid and particle refinement for the two-dimensional BKW test. Radial and non-radial costs are shown with $N=3n_v^2$, including lookup-grid setup and one collision evaluation (median of three trials for each component, after warm-up). The non-radial/radial cost ratio is approximately $1.25$--$1.58$ for $8\le n_v\le30$ (about $20$--$37\%$ less time), and $1.80$ at $n_v=5$.}
    \label{fig:radial_scaling}
\end{figure}

Figure~\ref{fig:radial_scaling} compares the cost with and without radial reduction under joint grid and particle refinement, keeping $N_{\rm pc}=3$ so that $N=3n_v^2$. The observed cost is approximately proportional to $n_v^4$, consistent with quadratic pairwise-interaction cost in $N$; unlike the convergence test in Figure~\ref{fig:nv_convergence}, this benchmark does not keep $N$ fixed.

\begin{figure}[!ht]
    \centering
    \includegraphics[width=\textwidth]{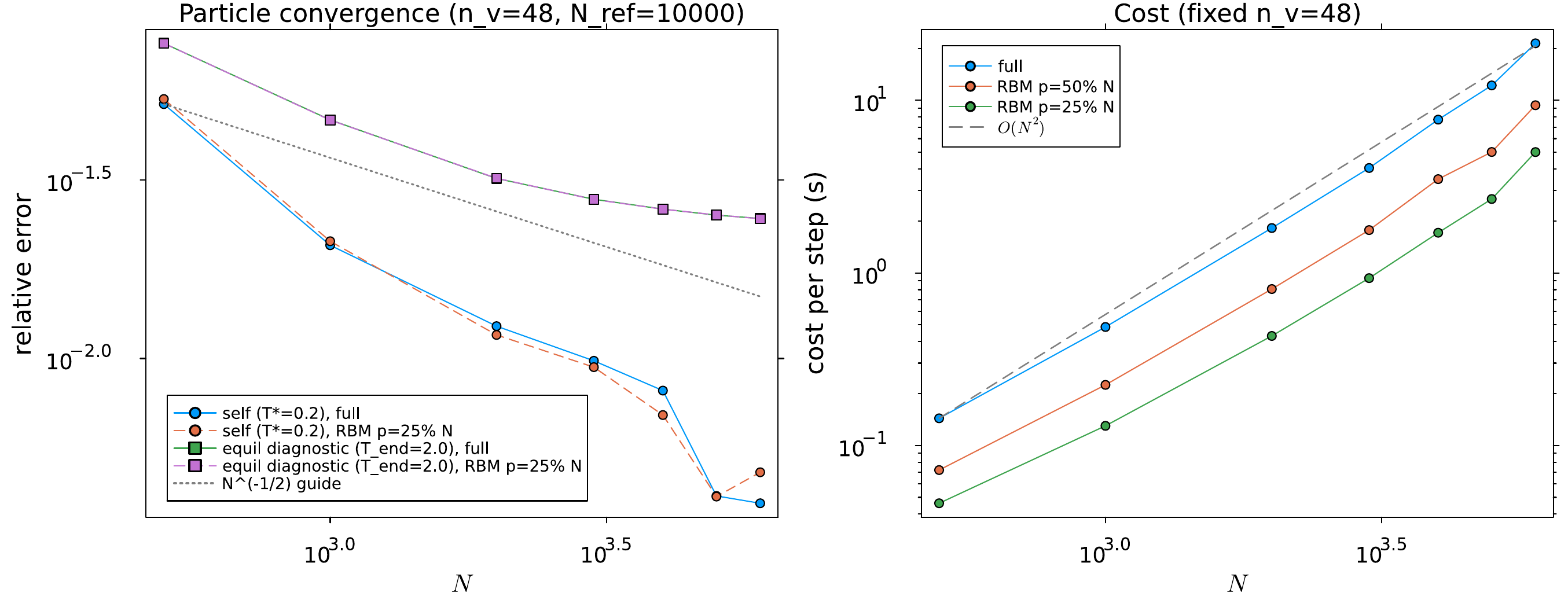}
    \caption{Particle-number convergence (left) and computational cost (right) for the two-dimensional BKW test at fixed $n_v=48$. Self-convergence errors use an independent $N^{\mathrm{ref}}=10000$ full-interaction reference on the same grid; equilibrium errors are measured at $T_{\mathrm{end}}=2$. The cost comparison includes full interactions and random batches of size $p=0.5N$ and $p=0.25N$.}
    \label{fig:N_convergence_cost}
\end{figure}

We also examine convergence in particle number at fixed $n_v=48$, keeping the mollification scale and all other discretization parameters unchanged. We vary
\[
N\in\{500,1000,2000,3000,4000,5000,6000\},
\]
and use an independent full-interaction reference with $N^{\mathrm{ref}}=10000$ on the same velocity grid. The relative self-convergence error is measured at $T_*=0.2$, and the distance to each run's matching Maxwellian at $T_{\mathrm{end}}=2$. Figure~\ref{fig:N_convergence_cost} shows decreasing self-convergence errors, with fitted rates $0.98$ for both full interactions and random batching with $p=0.25N$ over the tested range. The finite-time equilibrium error decreases more slowly and approaches a plateau. These rates are measured against a finite-particle reference. The same figure reports the median collision time per step over eight benchmark steps, discarding the first as warm-up. Full interactions exhibit approximately quadratic cost in $N$, while batches containing $50\%$ or $25\%$ of the particles reduce the cost without changing its quadratic scaling.

\subsection{Further out-of-equilibrium initial states}
\label{subsec:6.4}
\begin{figure}[!ht]
  \centering

  % Optional tightening
  \vspace{-2mm}

  \begin{tabular}{@{}c@{\hspace{0.01\textwidth}}c@{\hspace{0.01\textwidth}}c@{}c@{}}
    $t=0.0$ & $t=0.25$ & $t=7.5$ & \\[1mm]
    % -------- Scatter row --------
    \includegraphics[width=0.30\textwidth]{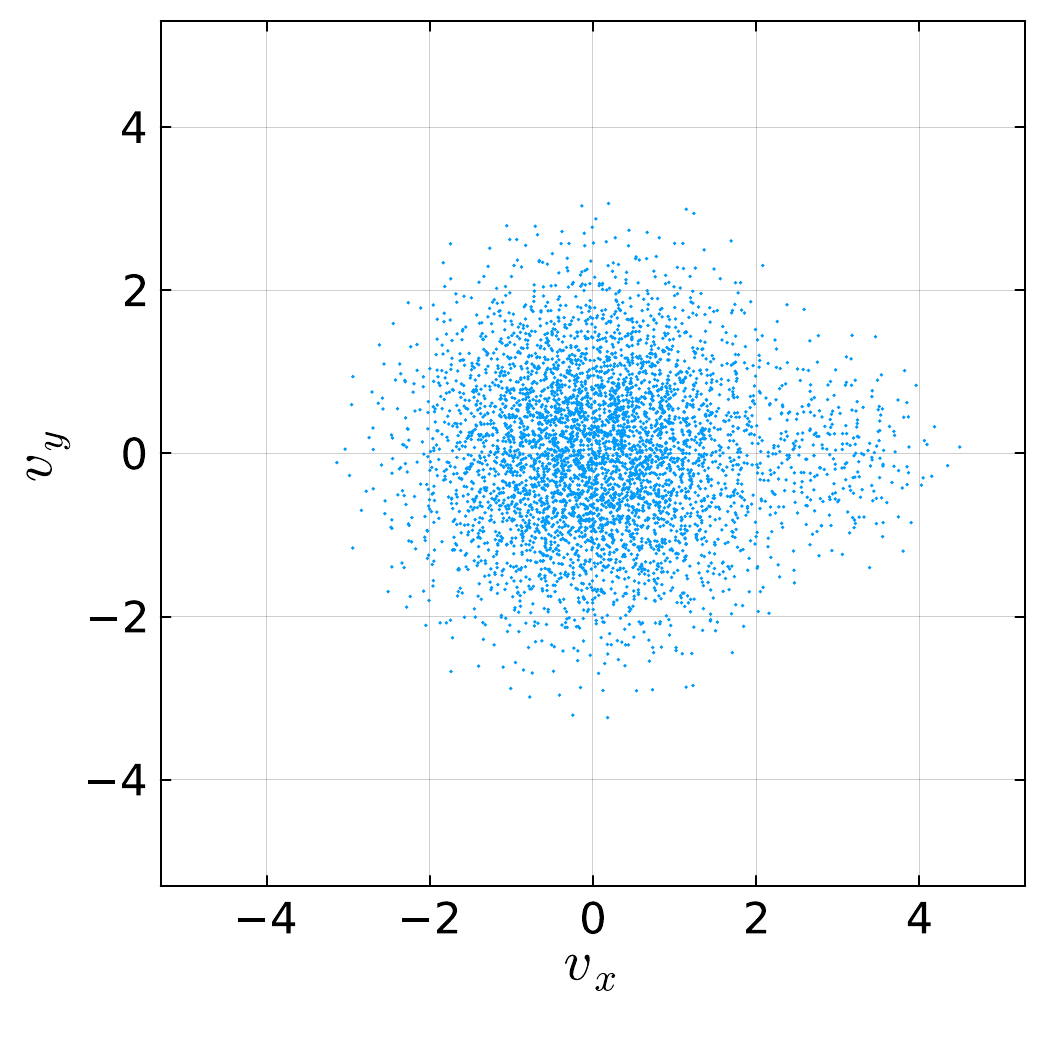} &
    \includegraphics[width=0.30\textwidth]{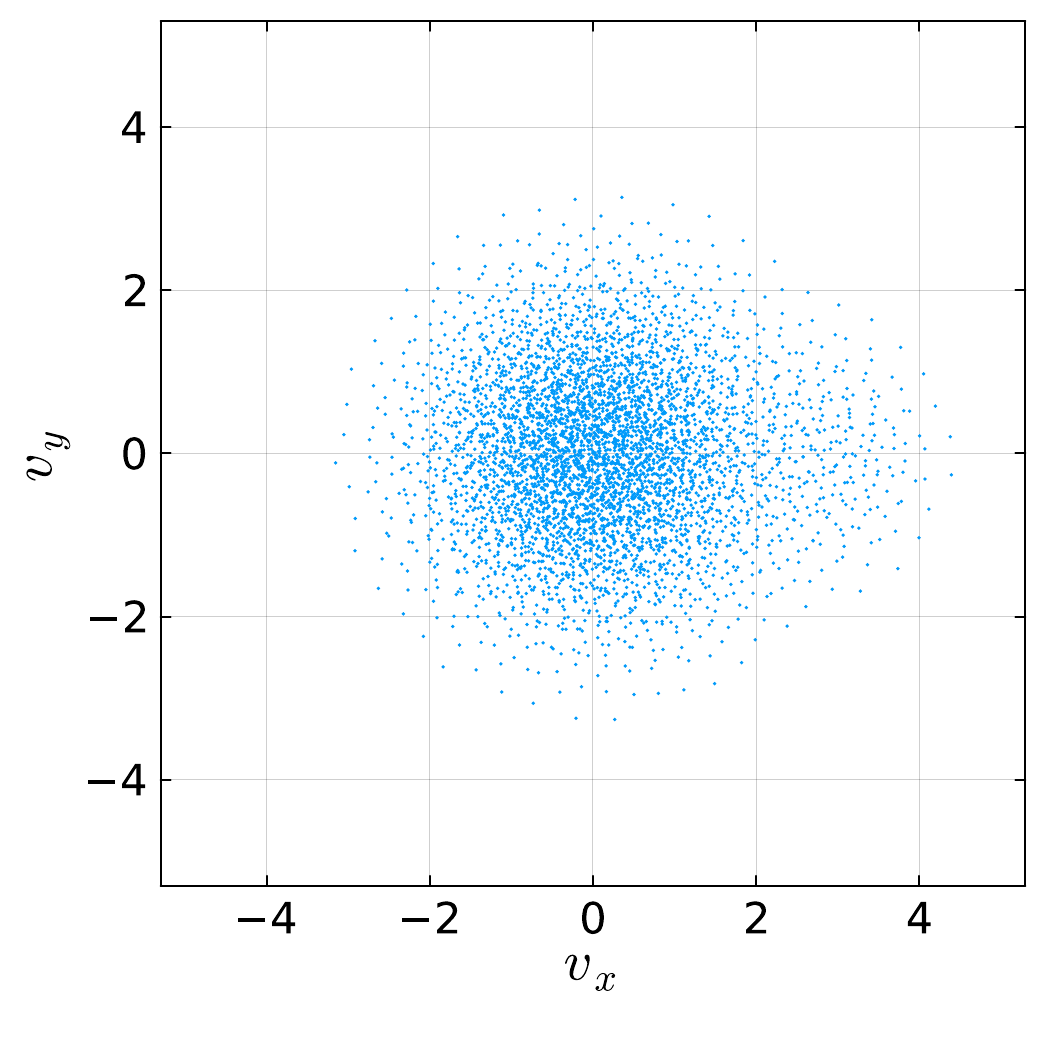} &
    \includegraphics[width=0.30\textwidth]{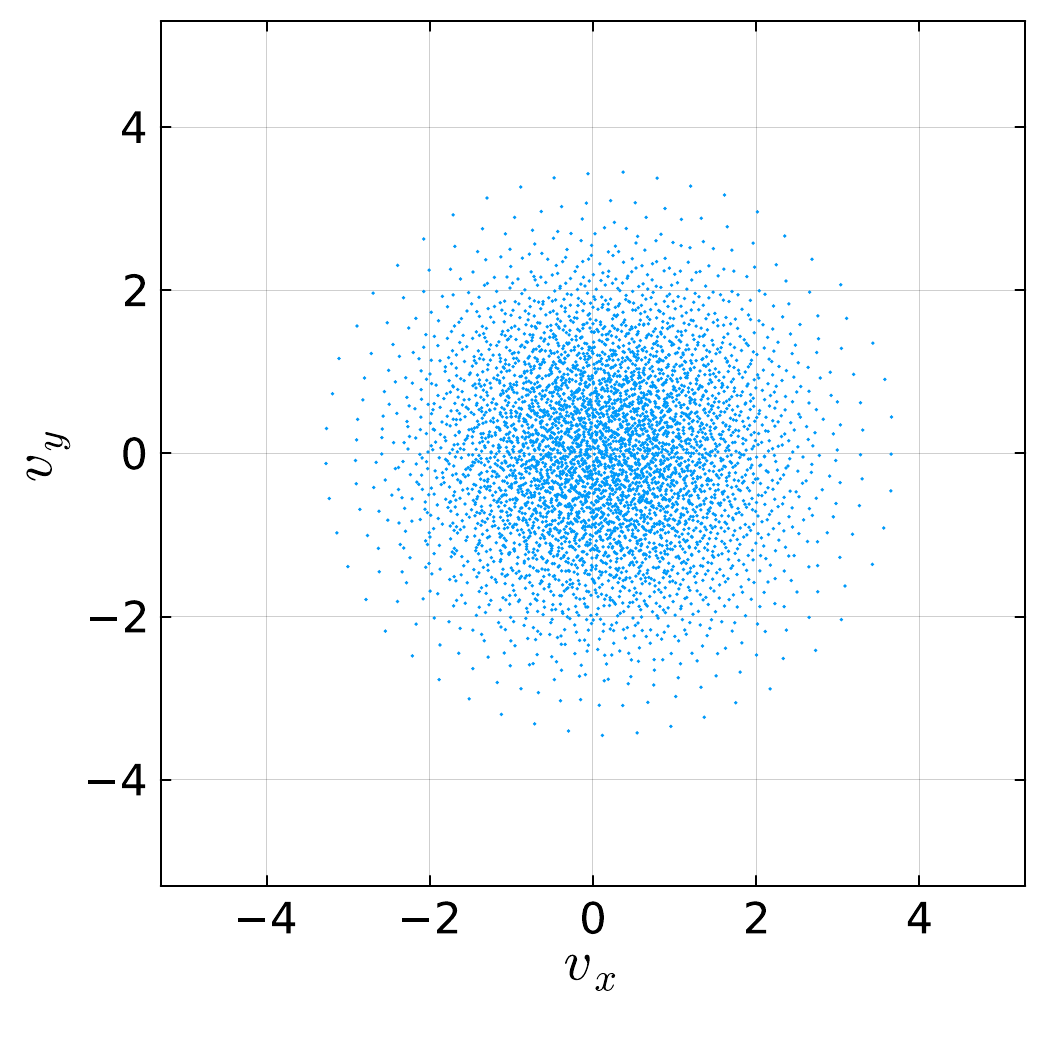} & \\[-2mm]

    % -------- Heatmap row --------
    % Match the axes boxes to the scatter panels, with a separate shared colorbar.
    \includegraphics[pagebox=mediabox,viewport=49 -21 573 503,clip,width=0.30\textwidth]{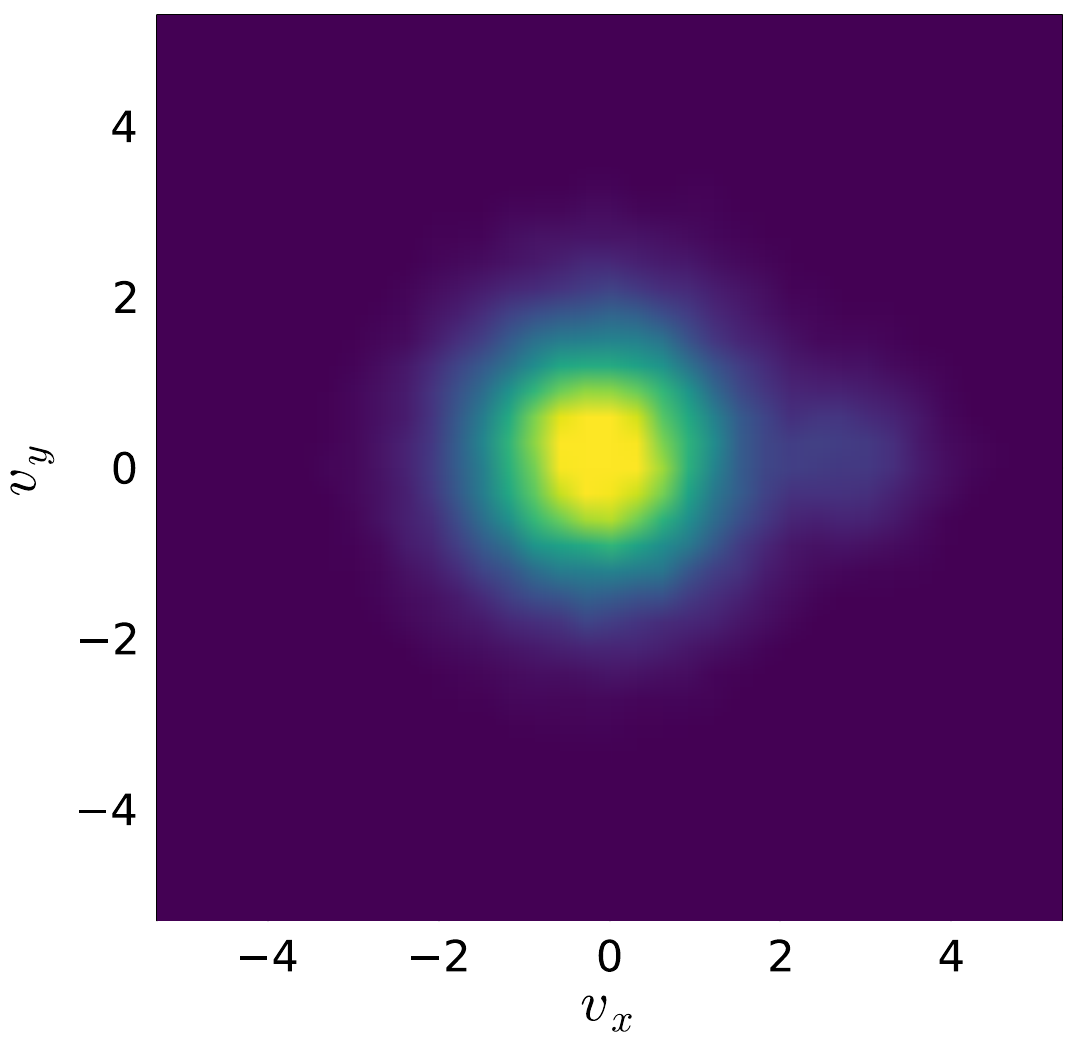} &
    \includegraphics[pagebox=mediabox,viewport=49 -21 573 503,clip,width=0.30\textwidth]{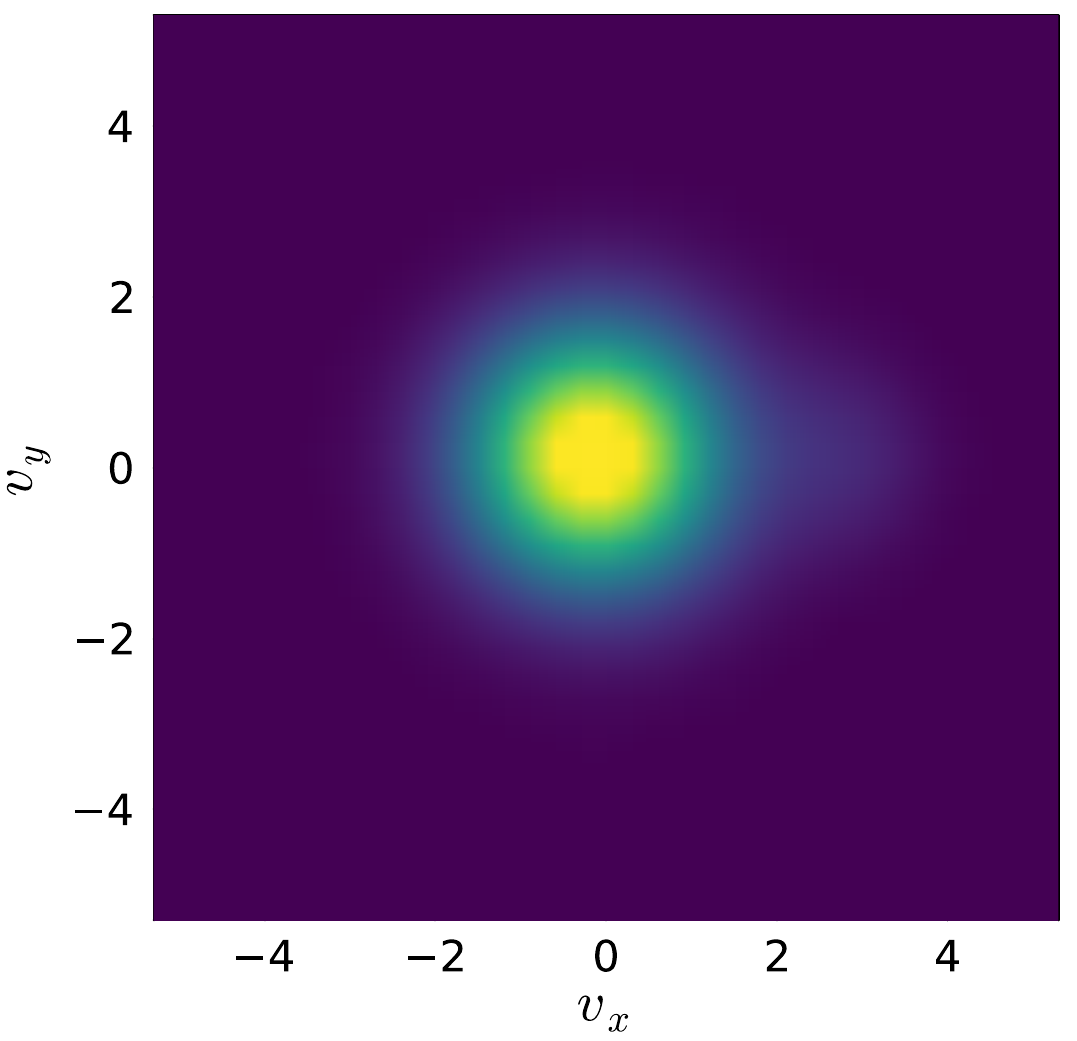} &
    \includegraphics[pagebox=mediabox,viewport=49 -21 573 503,clip,width=0.30\textwidth]{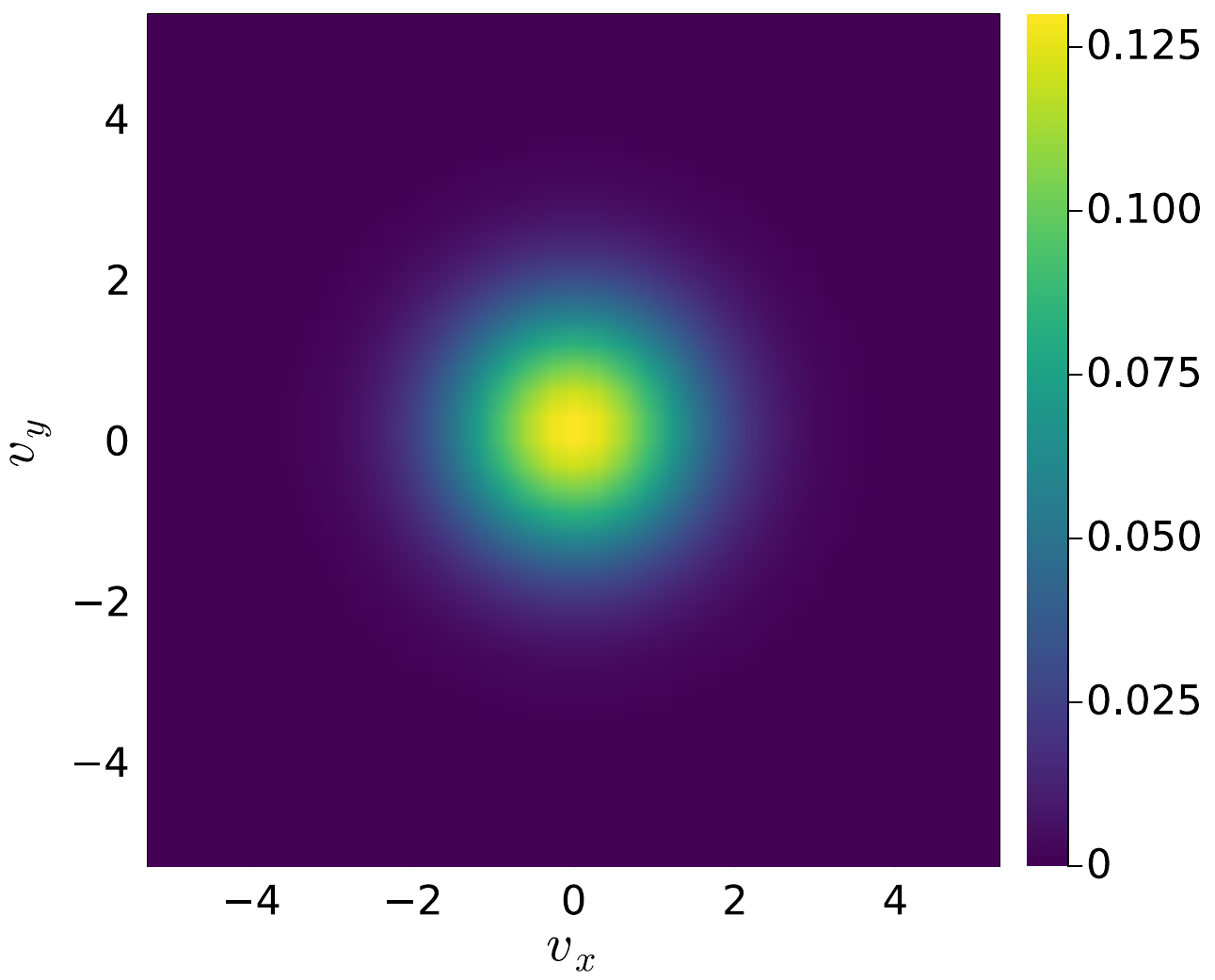} &
    \includegraphics[pagebox=mediabox,viewport=573 -21 675 503,clip,height=0.30\textwidth]{2dbump_frame_grid_0750.pdf}
  \end{tabular}

  \vspace{-2mm}

 \caption{Evolution of the two-dimensional bump-on-tail state. Particle velocities (top) and reconstructed density (bottom) at $t=0$, $0.25$, and $7.5$ from left to right.}

  \label{fig:2dbump_main}
\end{figure}

\begin{figure}[!ht]
  \centering

  \vspace{-2mm}

  \begin{tabular}{@{}cc@{}}
    \includegraphics[width=0.45\textwidth]{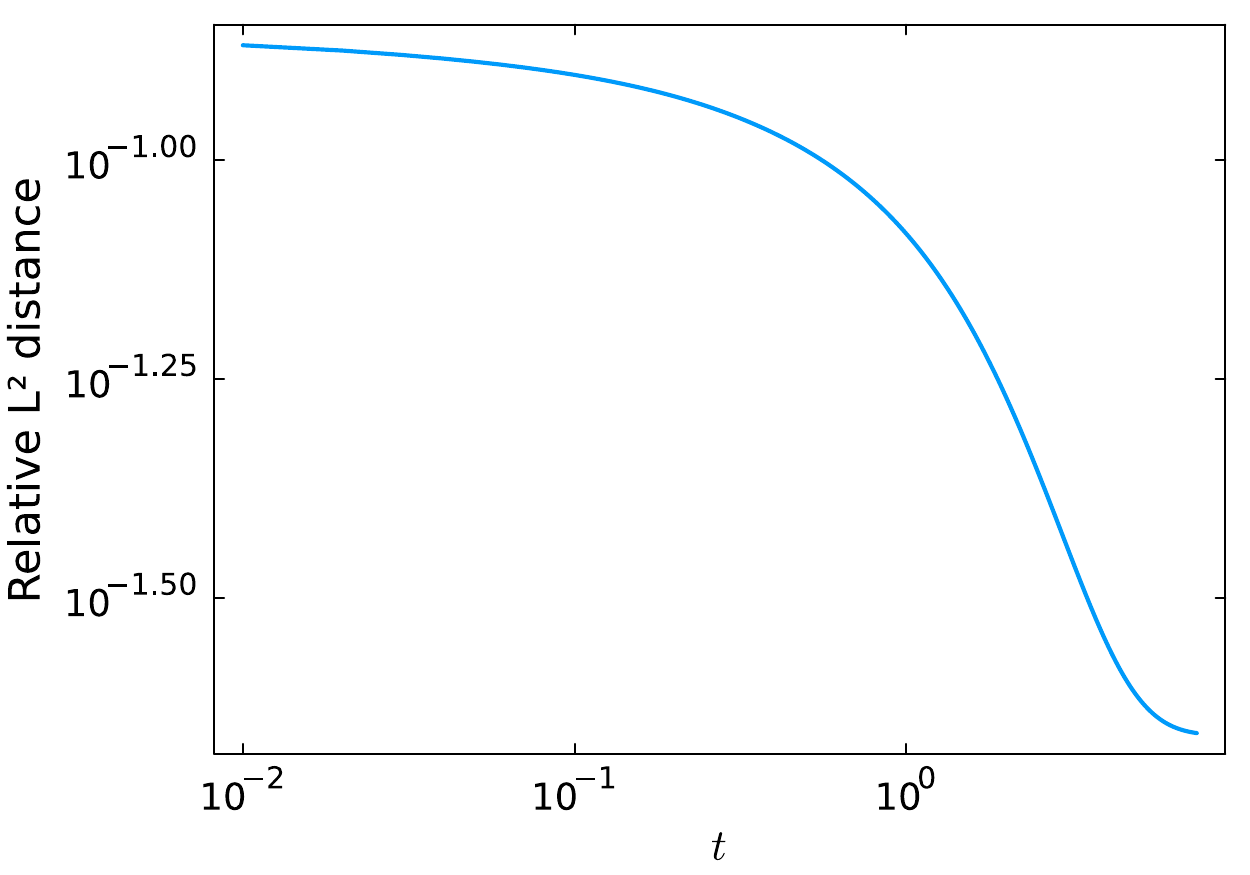} &
    \includegraphics[width=0.45\textwidth]{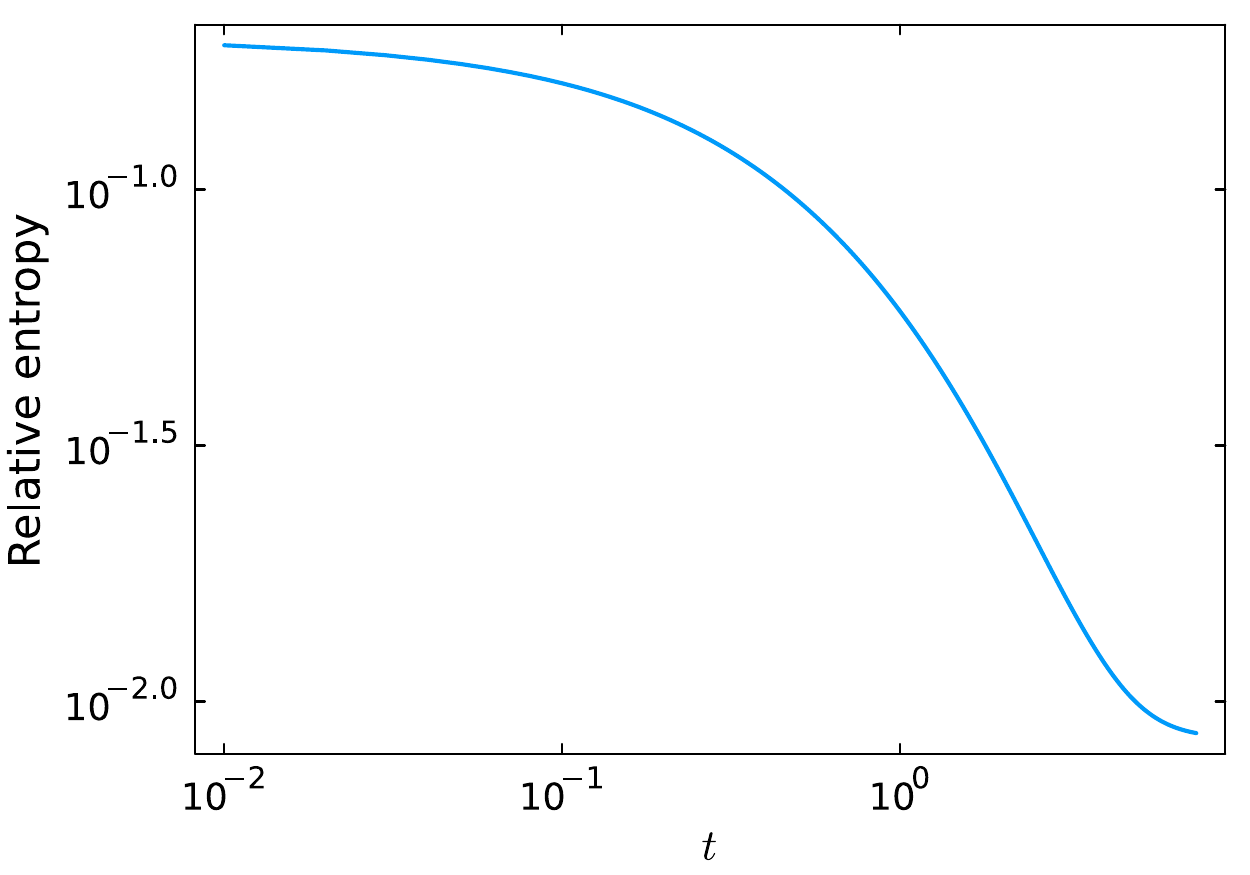}
  \end{tabular}

  \vspace{-2mm}

\caption{Relaxation of the bump-on-tail state towards the matching Maxwellian: relative $L^2$ error (left) and relative entropy (right).}

  \label{fig:2dbump_relaxation_metrics}
\end{figure}

\begin{figure}[!ht]
  \centering

  % Optional tightening
  \vspace{-2mm}

  \begin{tabular}{@{}c@{\hspace{0.01\textwidth}}c@{\hspace{0.01\textwidth}}c@{}}
    % -------- Scatter row --------
    \includegraphics[trim=80bp 40bp 70bp 0bp,clip,width=0.32\textwidth]{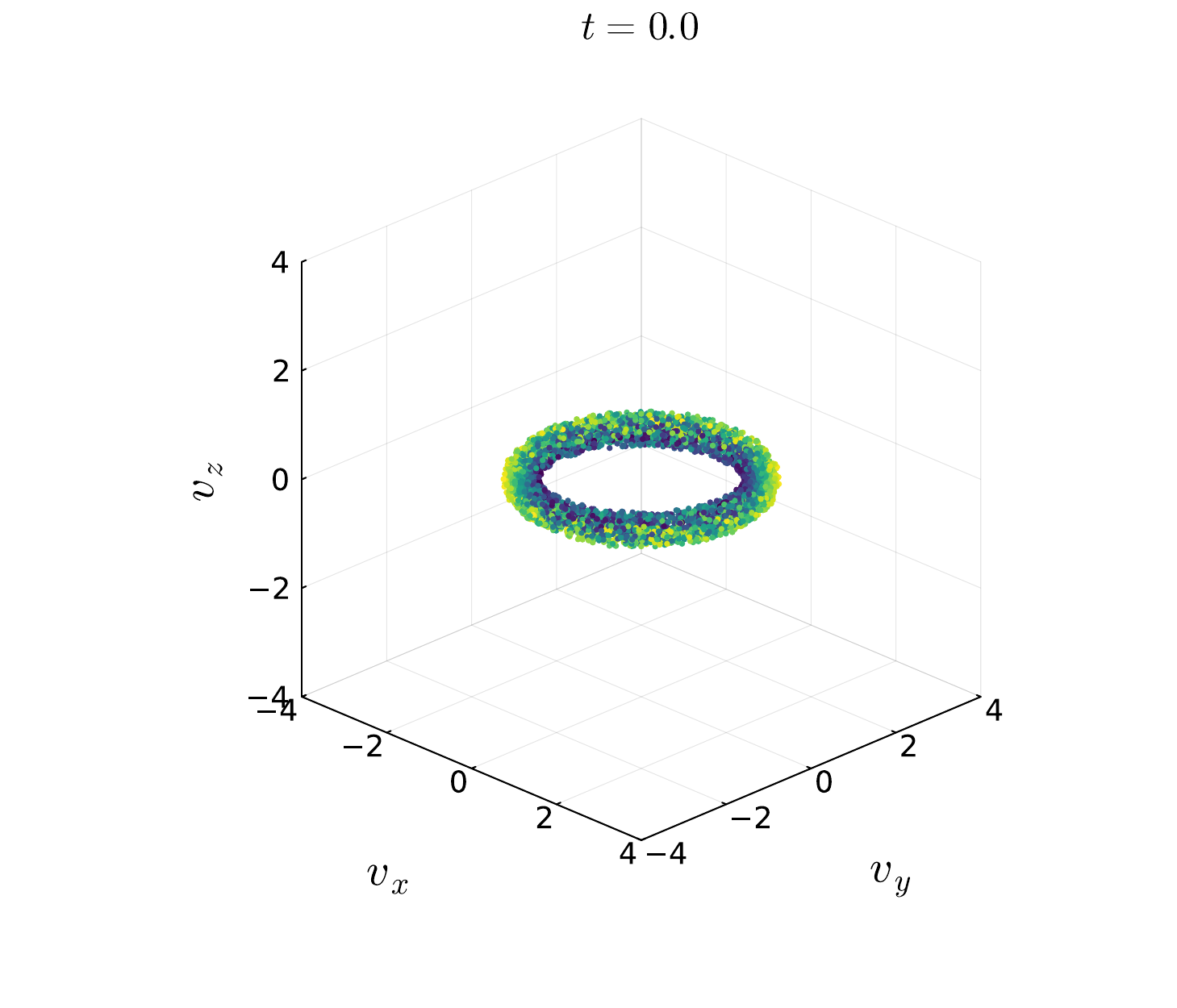} &
    \includegraphics[trim=80bp 40bp 70bp 0bp,clip,width=0.32\textwidth]{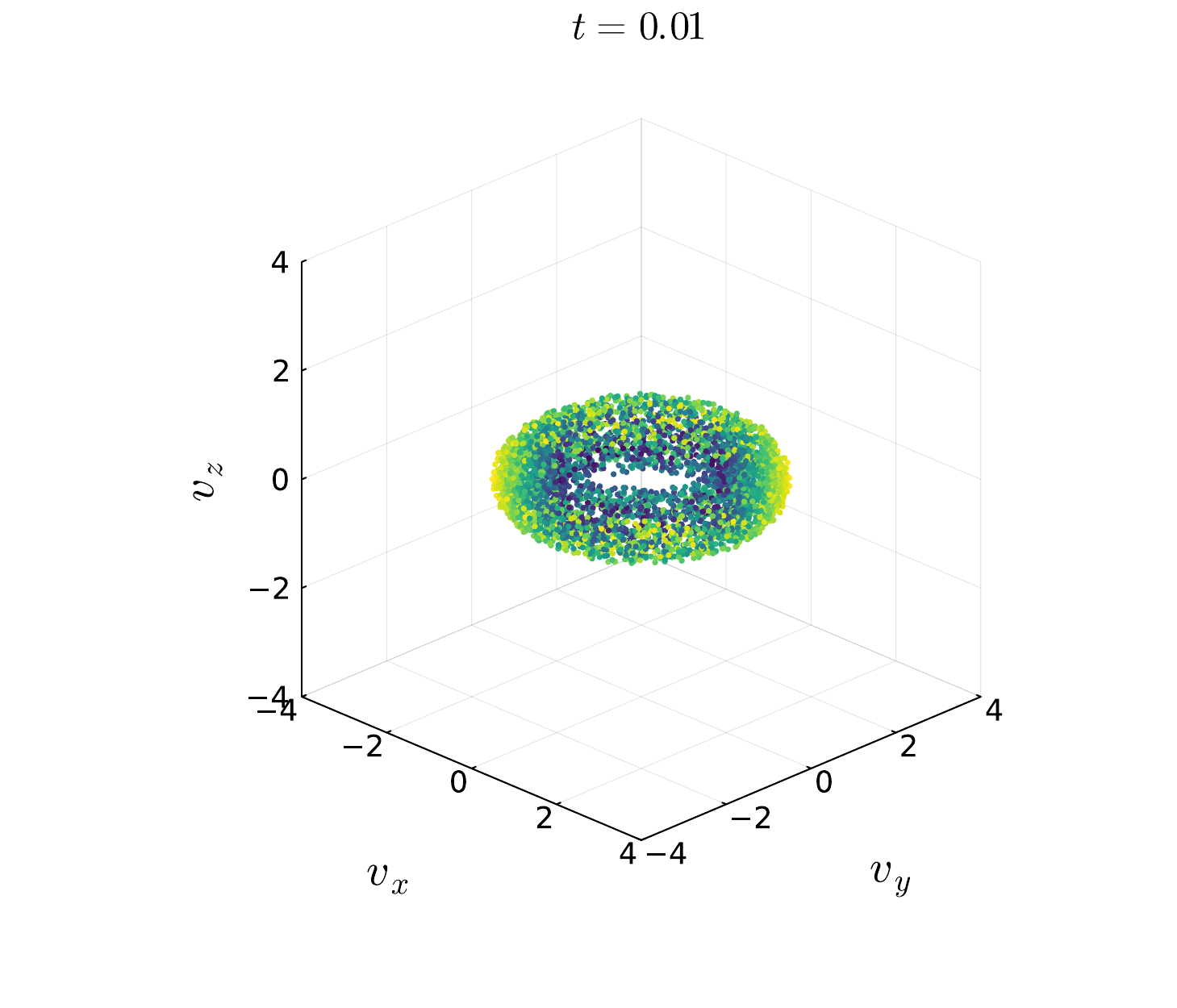} &
    \includegraphics[trim=80bp 40bp 70bp 0bp,clip,width=0.32\textwidth]{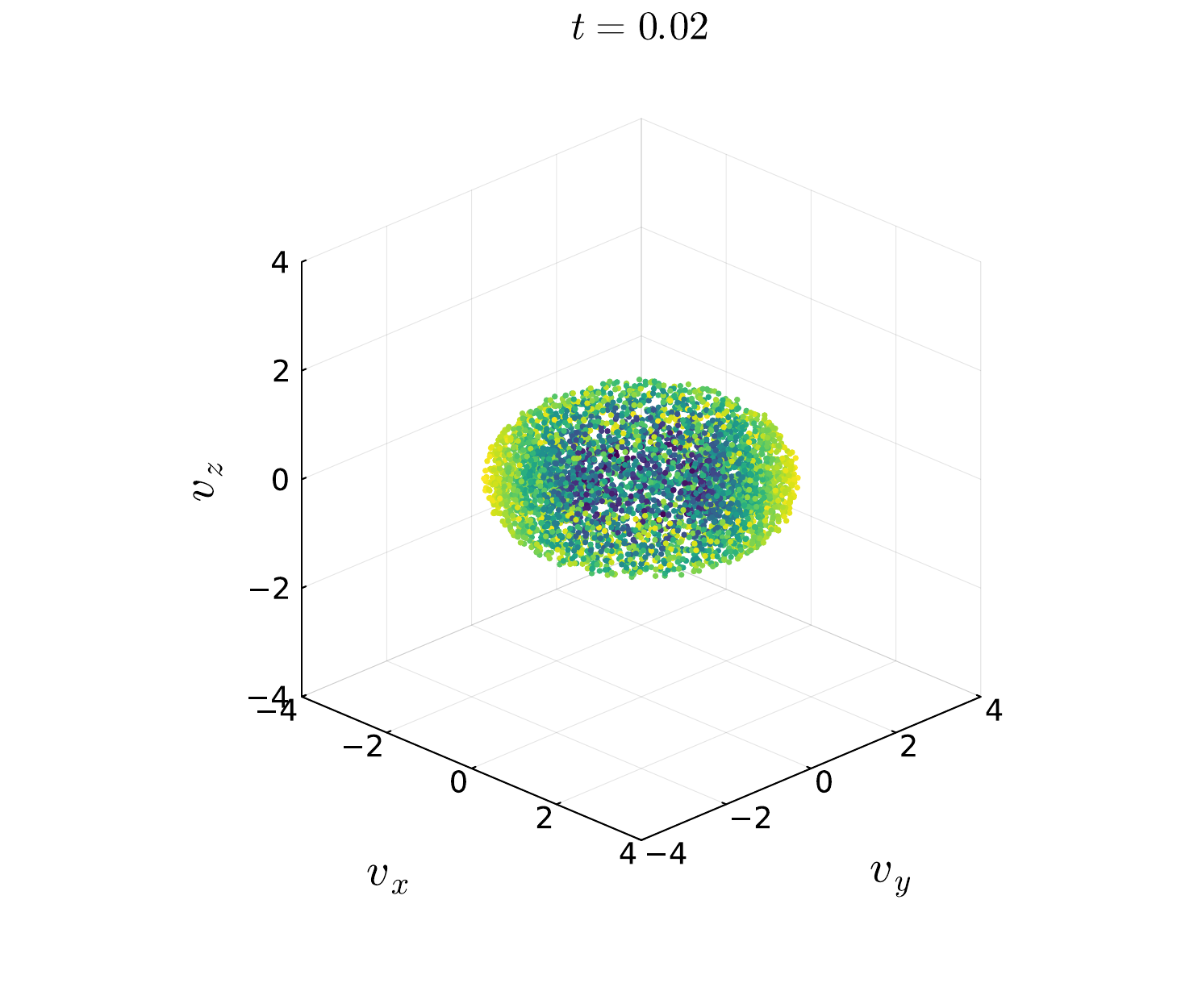} \\[2mm]

    % -------- Heatmap row --------
    \includegraphics[trim=80bp 40bp 70bp 0bp,clip,width=0.32\textwidth]{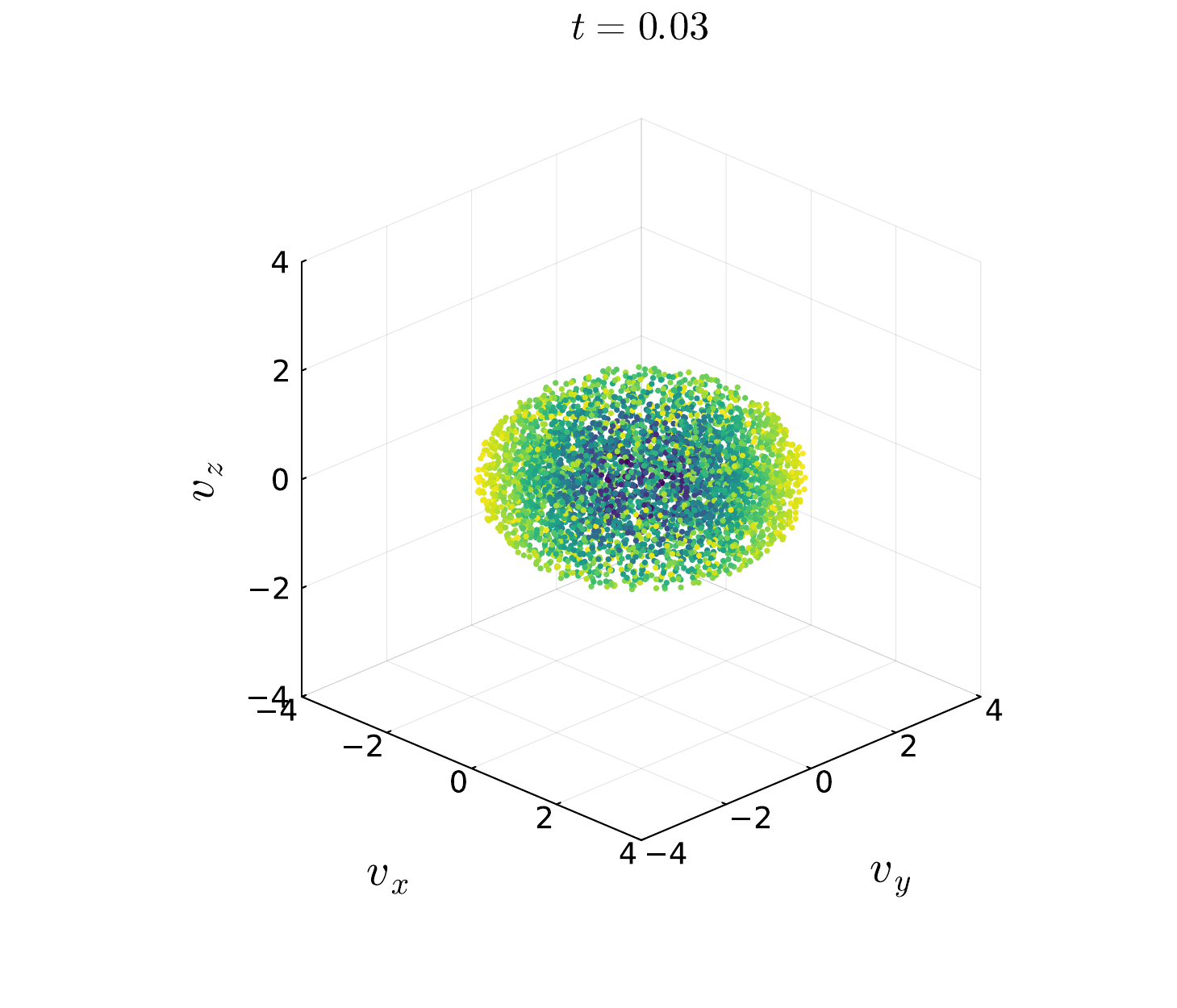} &
    \includegraphics[trim=80bp 40bp 70bp 0bp,clip,width=0.32\textwidth]{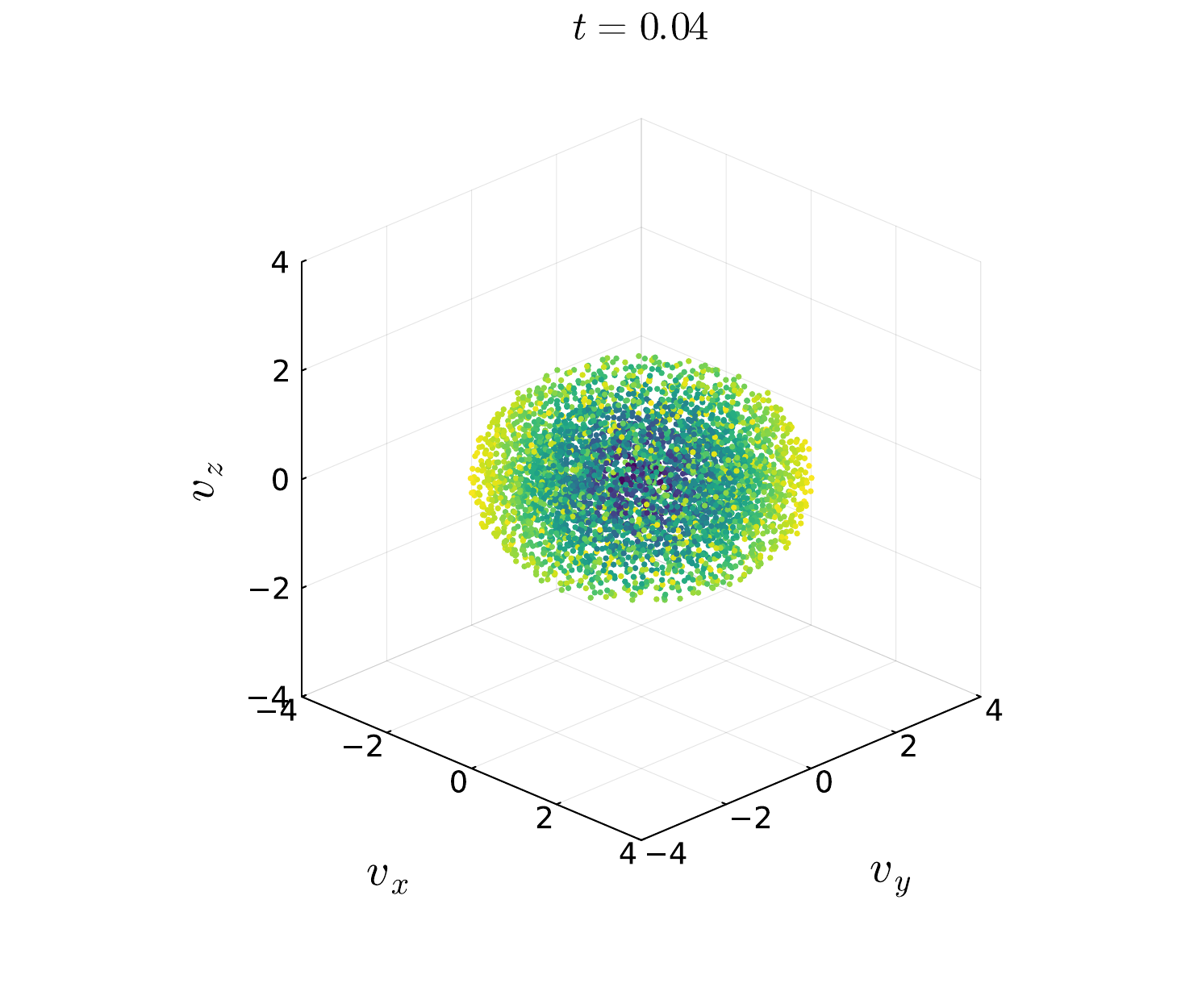} &
    \includegraphics[trim=80bp 40bp 70bp 0bp,clip,width=0.32\textwidth]{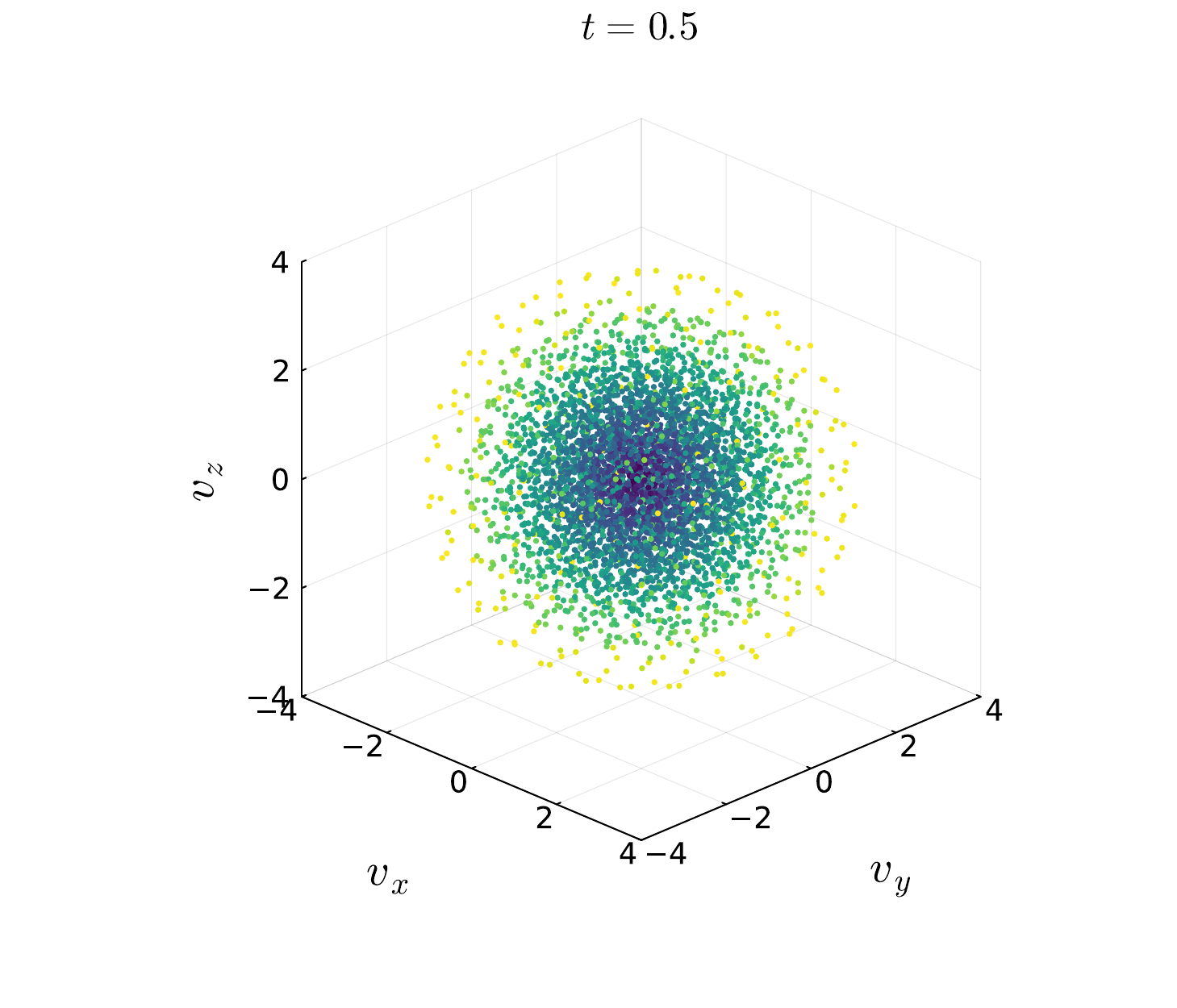}
  \end{tabular}

  \vspace{-2mm}

\caption{Evolution of the three-dimensional toroidal state at $t=0$, $0.01$, $0.02$, $0.03$, $0.04$, and $0.5$ from left to right and top to bottom. Particles are coloured by speed.}

  \label{fig:3dtoroidal}
\end{figure}

\begin{figure}[!ht]
  \centering

  \vspace{-2mm}

  \begin{tabular}{@{}cc@{}}
    \includegraphics[width=0.45\textwidth]{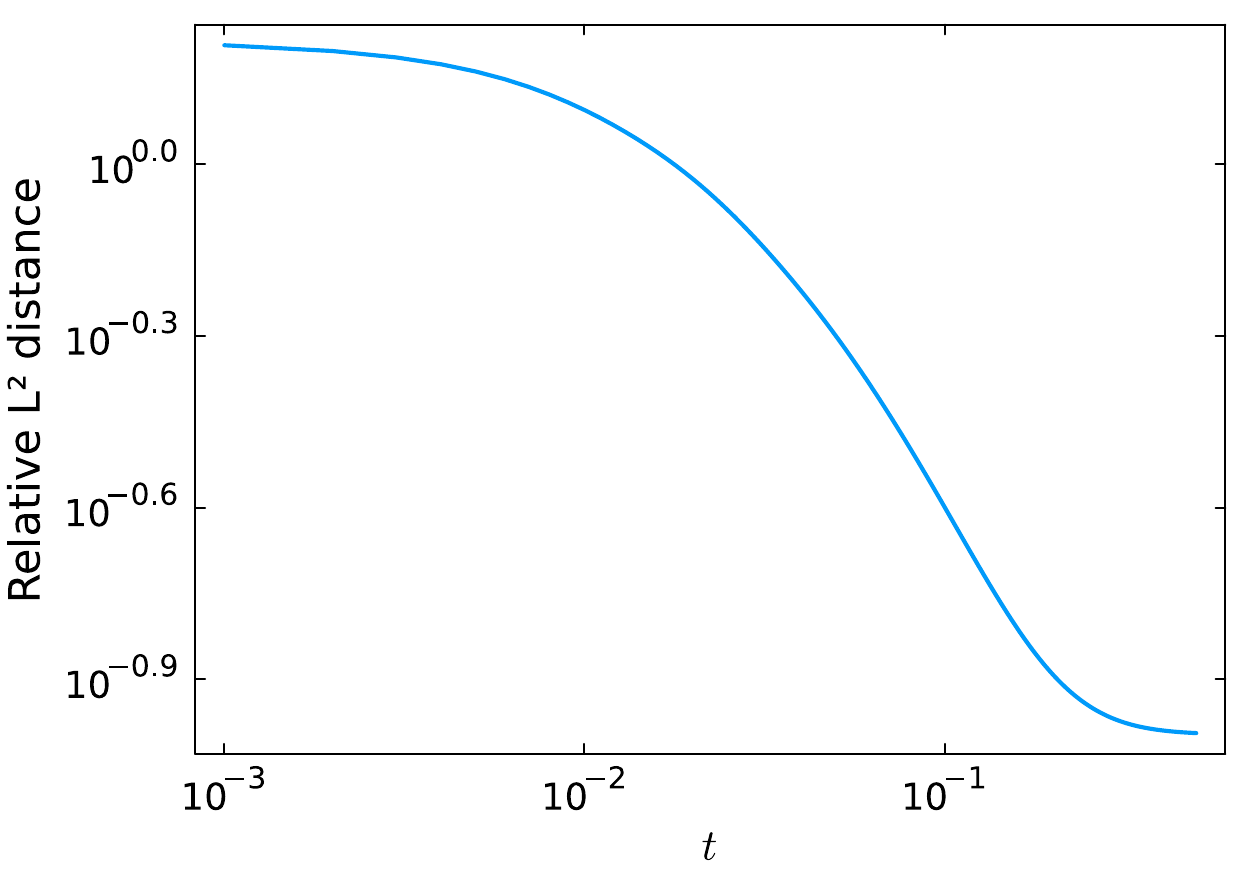} &
    \includegraphics[width=0.45\textwidth]{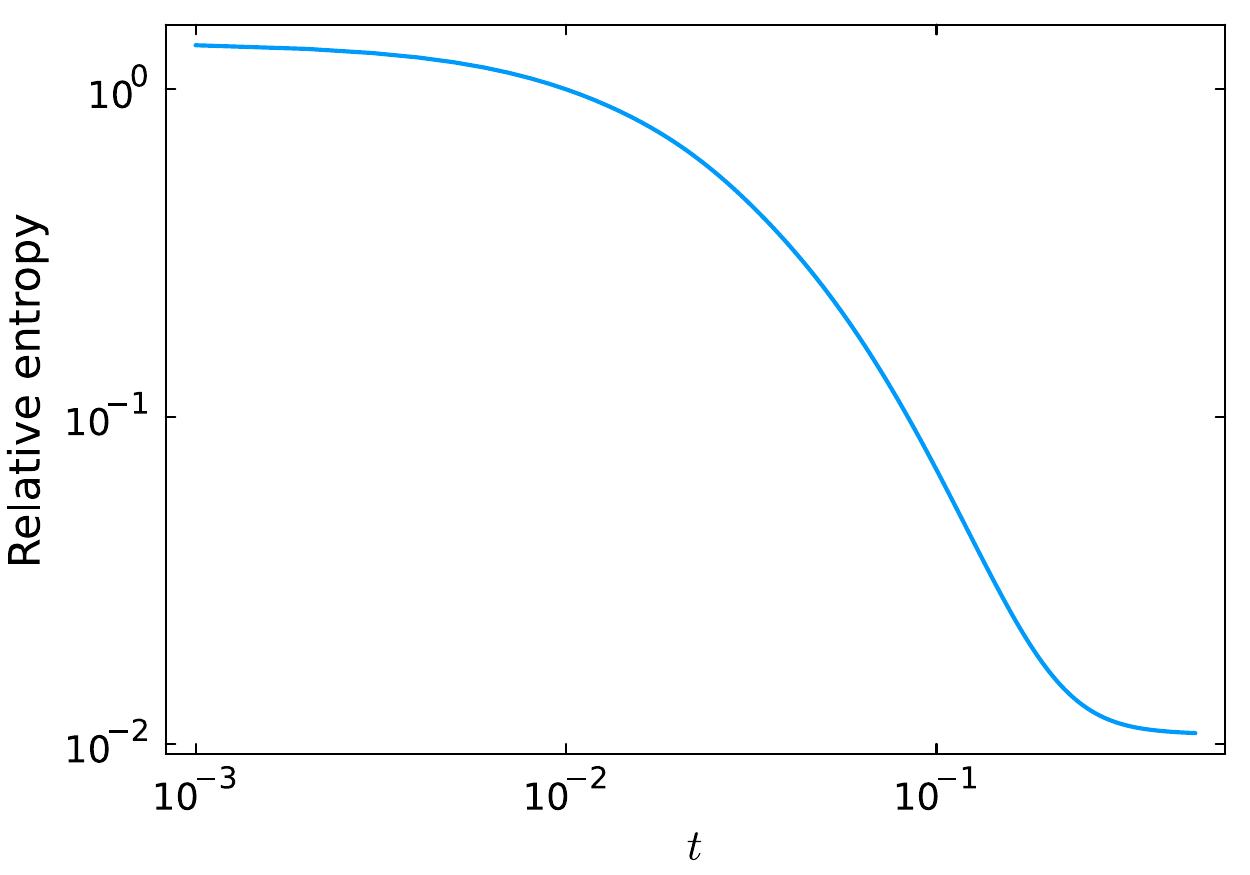}
  \end{tabular}

  \vspace{-2mm}

\caption{Relaxation of the toroidal state towards the matching Maxwellian: relative $L^2$ error (left) and relative entropy (right).}
  \label{fig:3dtoroidalrelaxation}
\end{figure}

The generality of the approach means that the evolution of essentially arbitrary initial distributions can be simulated. The evolution of a 2D bump-on-tail out-of-equilibrium initialisation is shown in Figure \ref{fig:2dbump_main}, while Figure \ref{fig:3dtoroidal} shows a 3D toroidal initial distribution evolving in time.
The bump-on-tail test uses $[-6,6]^2$, $n_v=40$, $N=4800$, $\Delta t=0.01$, full interactions, and transverse quadrature order $n_k=30$ with $k_{\perp,\max}=50$, $k_{\rm scale}=1$. The toroidal test uses $[-4,4]^3$, $n_v=13$, $N=4394$, $\Delta t=0.001$, requested batch size $p=1000$, and $n_k=30$ with $k_{\perp,\max}=5$, $k_{\rm scale}=1$. 
For the two-dimensional bump-on-tail test, we take 
\begin{equation*}
f_0(v)
=
0.95\,\mathcal{M}_{1}(v)
+
0.05\,\mathcal{M}_{0.3}\bigl(v-v_c\bigr),
\qquad v_c=(3,0)^T,\quad 
\mathcal{M}_{T}(v)
=
\frac{1}{2\pi T}
\exp\left(-\frac{|v|^2}{2T}\right).
\end{equation*}
For the three-dimensional toroidal test, we take
\begin{equation*}
f_0(v)
=
\frac{1
}{
2\pi^2a^2\sqrt{v_1^2+v_2^2}
} \mathbf{1}_{\left\{
(\sqrt{v_1^2+v_2^2}-R)^2+v_3^2\leq a^2
\right\}},
\qquad
R=2,\quad a=0.3.
\end{equation*}
The two-dimensional bump-on-tail experiment uses Coulomb interactions, $\widehat{V}(k) = |k|^{-2}$, while the three-dimensional toroidal experiment uses $\widehat{V}=1$ interactions, $\widehat{V}(k) = 1$. The corresponding convergence to equilibrium behaviour is plotted in Figures \ref{fig:2dbump_relaxation_metrics} and \ref{fig:3dtoroidalrelaxation}. Mass, energy and momentum remained approximately conserved throughout these simulations. Animations of these processes are also available on Figshare~\cite{LenardBalescuAnimations}.

\section*{Conclusions}
% \projectcomment{Andrea}{TO DO }
In this work, we propose a particle method to solve the spatially homogenous Lenard--Balescu equation. This method is based on the variational structure of the regularized collision operator, following the idea of the particle method for the classical Landau equation. Hence retaining the structural properties of the Lenard--Balescu equation, namely the conservation of mass, momentum and energy, as well as the decay of entropy. In addition, the Lenard--Balescu operator differs a lot with the Landau operator in the collision kernel, which is nonlocal and nonlinear, and describes the dynamic screening effect through a singular solution dependent dielectric response. To tackle the difficulty of the dielectric function, we rewrite the dielectric function into a explicit form with the help of a one-dimensional Hilbert transform of certain distribution dependent function. This function is an integral of the directional derivative of solution $f$ over the hyperplane orthogonal to $k$, hence favors a simplified form under the radial symmetric assumption.

This structure preserving particle method could be further extended to the non-homogenous case, the multi-species case and the quantum case. Moreover, theoretical justifications such as the gradient flow under some tailored metric, the convergence to the solution of original Lenard--Balescu as $\epsilon \to 0$ remain open. 

\section*{Acknowledgements}
JAC and AM acknowledge the support by the Advanced Grant Nonlocal-CPD,“Nonlocal PDEs for Complex Particle Dynamics: Phase Transitions, Patterns and Synchronization”, of the European Research Council Executive Agency (ERC) under the European Union’s Horizon 2020 research and innovation programme (grant 883363).
JAC was also partially supported by the EPSRC EP/V051121/1 and by the ``Maria de Maeztu'' Excellence Unit IMAG, reference CEX2020-001105-M, funded by MCIN/AEI/10.13039/501100011033/.
AM acknowledges the support by Fondo Italiano per la Scienza (FIS2023-01334) advanced grant “ADvanced numerical Approaches for MUltiscale Systems with uncertainties” - ADAMUS. YT acknowledges the support by the Deutsche Forschungsgemeinschaft (DFG) through the Priority Programme SPP2298 Theoretical Foundations of Deep Learning.

\bibliographystyle{abbrv}
\bibliography{References.bib}

\end{document}